\documentclass[10t]{article}          \usepackage{amsthm}
\usepackage{authblk}
\usepackage[margin=2cm]{geometry}
\usepackage{multicol}

\RequirePackage{calrsfs}
\RequirePackage{amsmath}
\RequirePackage{amsthm}

\RequirePackage{amssymb}
\RequirePackage{amsfonts}

\RequirePackage[most]{tcolorbox}

\RequirePackage{bbm}

\RequirePackage{graphicx}

\RequirePackage{booktabs}

\RequirePackage{url}

\RequirePackage{hyperref}
\hypersetup{
  colorlinks=true,
  linkcolor=blue,
  filecolor=magenta,
  urlcolor=cyan,
  citecolor=green!60!black,
  backref=page
}

\RequirePackage{xcolor}
\RequirePackage{color}

\RequirePackage{wasysym}
\RequirePackage{listings}
\RequirePackage{caption}
\RequirePackage{subcaption}

\RequirePackage{stmaryrd}

\RequirePackage{lineno}

\RequirePackage{xparse}

\RequirePackage{float}

\RequirePackage{multirow}

\RequirePackage{cellspace}

\RequirePackage{shellesc}

\RequirePackage[ruled,vlined,linesnumbered]{algorithm2e}

\SetCommentSty{xCommentSty}

\RequirePackage{tikz}
\usetikzlibrary{spy, calc, automata, positioning, arrows, external, matrix, shapes.misc}
\tikzset{cross/.style={cross out,thick,draw=black,minimum size=2*(#1-\pgflinewidth), inner sep=0pt, outer sep=0pt},
cross/.default={2.5pt}}

\RequirePackage{pgfplots}
\pgfplotsset{compat=1.16}
\usepgfplotslibrary{fillbetween}
\usepgfplotslibrary{groupplots}
\usepgfplotslibrary{external}

\DeclareMathAlphabet{\pazocal}{OMS}{zplm}{m}{n}
\DeclareMathAlphabet\bpazocal{OMS}{cmsy}{b}{n}

\providecommand{\meas}{\ensuremath{\text{meas}}}

\providecommand{\supp}{\ensuremath{\text{supp}}}

\providecommand{\bR}{\mathbb{R}}
\providecommand{\bs}{\boldsymbol}

\providecommand{\T}{\mathsf{T}}

\providecommand{\half}{\frac{1}{2}}

\providecommand{\vertiii}[1]{{\left\vert\kern-0.15ex\left\vert\kern-0.15ex\left\vert #1
    \right\vert\kern-0.15ex\right\vert\kern-0.15ex\right\vert}}

\NewDocumentCommand{\curlii}{sO{}m}
{
	\IfBooleanTF{#1}
    {\dgalext{#3}}
    {\dgalx[#2]{#3}}
}

\NewDocumentCommand{\dgalext}{m}{  \sbox0{    \mathsurround=0pt     $\left\{\vphantom{#1}\right.\kern-\nulldelimiterspace$  }  \sbox2{\{}  \ifdim\ht0=\ht2
    \{\kern-.625\wd2 \{#1\}\kern-.625\wd2 \}  \else
    \left\{\kern-.7\wd0\left\{#1\right\}\kern-.7\wd0\right\}  \fi
}

\NewDocumentCommand{\dgalx}{om}{  \sbox0{\mathsurround=0pt$#1\{$}  \sbox2{\{}  \ifdim\ht0=\ht2
    \{\kern-.625\wd2 \{#2\}\kern-.625\wd2 \}  \else
    \mathopen{#1\{\kern-.7\wd0 #1\{}
    #2
    \mathclose{#1\}\kern-.7\wd0 #1\}}
  \fi
}

\tikzset{
  partial ellipse/.style args={#1:#2:#3}{
    insert path={+ (#1:#3) arc (#1:#2:#3)}
  }
}

\newcommand{\overbar}[1]{\mkern 1.5mu\overline{\mkern-1.5mu#1\mkern-1.5mu}\mkern 1.5mu}

\providecommand{\pE}{\pazocal{E}}

\providecommand{\pG}{\pazocal{G}}

\providecommand{\pJ}{\pazocal{J}}

\providecommand{\pM}{\pazocal{M}}

\providecommand{\pO}{\pazocal{O}}

\providecommand{\pQ}{\pazocal{Q}}

\providecommand{\pT}{\pazocal{T}}

\providecommand{\bpK}{\bpazocal{K}}

\providecommand{\bpV}{\bpazocal{V}}

\providecommand{\bpX}{\bpazocal{X}}

\providecommand{\mA}{\mathcal{A}}

\providecommand{\mN}{\mathcal{N}}

\providecommand{\mV}{\mathcal{V}}

\providecommand{\matA}{\bs{A}}
\providecommand{\matB}{\bs{B}}

\providecommand{\matE}{\bs{E}}

\providecommand{\matI}{\bs{I}}

\providecommand{\matO}{\bs{O}}

\providecommand{\matQ}{\bs{Q}}
\providecommand{\matR}{\bs{R}}

\providecommand{\matT}{\bs{T}}

\providecommand{\matX}{\bs{X}}

\providecommand{\vecb}{\bs{b}}
\providecommand{\vecc}{\bs{c}}

\providecommand{\vece}{\bs{e}}
\providecommand{\vecf}{\bs{f}}
\providecommand{\vecg}{\bs{g}}

\providecommand{\vecn}{\bs{n}}

\providecommand{\vecr}{\bs{r}}

\providecommand{\vect}{\bs{t}}
\providecommand{\vecu}{\bs{u}}
\providecommand{\vecv}{\bs{v}}
\providecommand{\vecw}{\bs{w}}
\providecommand{\vecx}{\bs{x}}
\providecommand{\vecy}{\bs{y}}
\providecommand{\vecz}{\bs{z}}

\providecommand{\to}{\widetilde{o}}

\usepackage{color, colortbl}
\definecolor{Gray}{gray}{0.9}

\graphicspath{{./pngfigs/}}

\providecommand{\keywords}[1]
{
  \small
  \textbf{\textit{\ \ \ \ Keywords:}} #1
}

\makeatletter
\if@twocolumn

\else

\fi
\makeatother

\newtheorem{remark}{Remark}

\graphicspath{{./tikzpics/}}

\begin{document}
\title{A Generalized Multigrid Method for Solving Contact Problems in Lagrange Multiplier based Unfitted Finite Element Method}

\author[]{Hardik Kothari\thanks{hardik.kothari@usi.ch}}
\author[]{Rolf Krause\thanks{rolf.krause@usi.ch}}
\affil[]{Euler Institute, Universit\`{a} della Svizzera italiana, Lugano, Switzerland}

\maketitle
\begin{abstract}
  Internal interfaces in a domain could exist as a material defect or they can appear due to propagations of cracks.
Discretization of such geometries and solution of the contact problem on the internal interfaces can be computationally challenging.
We employ an unfitted Finite Element (FE) framework for the discretization of the domains and develop a tailored, globally convergent, and efficient multigrid method for solving contact problems on the internal interfaces.
In the unfitted FE methods, structured background meshes are used and only the underlying finite element space has to be modified to incorporate the discontinuities.
The non-penetration conditions on the embedded interfaces of the domains are discretized using the method of Lagrange multipliers.
We reformulate the arising variational inequality problem as a quadratic minimization problem with linear inequality constraints.
Our multigrid method can solve such problems by employing a tailored multilevel hierarchy of the FE spaces and a novel approach for tackling the discretized non-penetration conditions.
We employ pseudo-$L^2$ projection-based transfer operators to construct a hierarchy of nested FE spaces from the hierarchy of non-nested meshes.
The essential component of our multigrid method is a technique that decouples the linear constraints using an orthogonal transformation of the basis.
The decoupled constraints are handled by a modified variant of the projected Gauss-Seidel method, which we employ as a smoother in the multigrid method.
These components of the multigrid method allow us to enforce linear constraints locally and ensure the global convergence of our method.
We will demonstrate the robustness, efficiency, and level independent convergence property of the proposed method for Signorini's problem and two-body contact problems.

 \end{abstract}

\section{Introduction}

Contact problems are virtually ubiquitous in the field of mechanics and engineering.
An accurate and reliable simulation of the contact problem is important in many engineering applications.
From the numerical modeling point of view, contact problems are challenging to solve as the contact boundary is unknown a priori.
Hence, a special type of iterative scheme is needed to solve such problems, as the contact zone has to be identified during the solution process.
In this work, we present contact problems in the unfitted finite element (FE) framework and introduce the problem in terms of a variational inequality.
Here, we consider frictionless contact problems, where we neglect the tangential forces on the contact interfaces.
The main contribution of this work is a novel generalized multigrid method that is developed to solve the arising quadratic minimization problem with linear inequality constraints.

In the last two decades, unfitted FE methods have seen a rise in popularity and multiple frameworks for handling unfitted geometries have emerged.
Unlike traditional FE methods, these unfitted FE methods do not require a fitted mesh that describes the computational domain explicitly.
The unfitted FE methods, generally, require a background mesh that encapsulates the computational domain, and the FE spaces associated with the background meshes are modified to capture the information of the domain.
These methods are ideal for solving problems with complex computational domains, interface problems with discontinuous coefficients, or moving interfaces.
As the background mesh and the computational domain are created independently, the interfaces/boundaries of the domain are, generally, embedded in the background mesh.
For this reason, it becomes essential to enforce boundary conditions or interface conditions in a weak sense.
To this end, the penalty method, the method of Lagrange multipliers, and Nitsche's method are used to impose the interface/boundary conditions.
In practice, Nitsche's method is significantly more popular than the method of Lagrange multipliers and the penalty method.
This is due to the fact that the penalty method is variationally inconsistent and the method does not produce optimal convergence rates of the discretization error in absence of a sufficiently large penalty parameter.
Nitsche's method can be regarded as a variationally consistent penalty method and due to its robustness, it is widely used in the unfitted FE methods.
Nitsche's method also requires a penalty/stabilization parameter, where the parameter has to be chosen such that the coercivity of the bilinear form is ensured.
The method of Lagrange multipliers gives rise to mixed FE formulations, the linear systems of equations stemming from this discretization scheme have a saddle point structure and the method is not stable if the FE spaces do not satisfy discrete inf-sup condition.
The eXtended finite element method (XFEM) was introduced as a partition of unity method to enrich the underlying FE spaces to tackle the problems in fracture mechanics with crack propagation~\cite{sukumar_extended_2000, moes_finite_1999}.
A similar enrichment scheme in combination with Nitsche's method was introduced to handle the interface problems~\cite{hansbo_unfitted_2002,hansbo_finite_2004}.
This method evolved into the CutFEM method~\cite{burman_cutfem:_2015} which in addition to Nitsche's method is also equipped with a ghost penalty stabilization scheme~\cite{burman_ghost_2010}.

The modeling of contact problems in the context of fitted FE methods has been studied from both numerical and theoretical points of view in detail, for example, in~\cite{peterwriggers2006-06-28, wohlmuth_variationally_2011, toda.laursen2013-03-13}.
In the unfitted FE framework, the initial work regarding the contact problem was carried out by Dolbow et al.~in the context of XFEM to tackle the frictional sliding contact on the crack faces~\cite{dolbow_extended_2001}.
In the work of Dolbow et al., the fictional contact constraints on the crack faces were handled using the Large Time Increment (LaTIn) method~\cite{ladeveze2012nonlinear}.
Since the initial work, the penalty method, Nitsche's method, and the method of Lagrange multipliers have been pursued for solving the contact problems in unfitted frameworks.
The penalty formulation for the contact between the open crack faces was utilized in a few works~\cite{khoei_contact_2006,liu_contact_2008,mueller-hoeppe_crack_2012}.
Nitsche's method has also been proposed for solving the frictional contact problem in the unfitted framework~\cite{coon_nitsche-eXtended_2011}.
A formal theoretical framework of Nitsche's method for solving the contact problem was later given by Chouly et al.~\cite{chouly_nitsche-based_2013,chouly_symmetric_2015} and later it was extended to the fictitious domain methods~\cite{fabre_fictitious_2016}.
As employing Nitsche's method for enforcing the contact conditions gives rise to non-smooth energy functional and in order to solve such a problem the generalized Newton's method was used as a solution strategy~\cite{renard_generalized_2013}.
In the context of the CutFEM solver, a LaTIn-based solution scheme was proposed for solving the contact problems, where the contact condition is handled with Nitsche's method~\cite{claus_stable_2018}.
Lagrange multiplier based approaches for contact problems in unfitted methods have also been pursued in several works~\cite{kim_mortared_2007,bechet_stable_2009,nistor_x-fem_2009,akula_mortex_2019}.
In these approaches, multiple techniques for constructing a stable multiplier space have been considered.
In our work, we utilize the method of Lagrange multipliers to discretize the contact condition.
In addition, we utilize a ghost penalty stabilization term in our bilinear form to bound the condition number of the arising linear system of equation by means of the mesh size of the background meshes.

As noted earlier, the method of Lagrange multipliers gives rise to mixed FE formulations, and the stability of the mixed formulation is ensured only if the discrete inf-sup condition is satisfied.
In the unfitted FE framework, it is shown that the most convenient approaches to construct the multiplier spaces give rise to instabilities~\cite{belytschko_structured_2003,ji_strategies_2004}.
To circumvent the strict requirement of satisfying the inf-sup condition, a different approach was introduced by Barbosa and Hughes~\cite{barbosa_finite_1991}.
In the Barbosa-Hughes approach, the restriction over the choices for FE spaces is dropped and the stability of the formulation is ensured using a stabilization term~\cite{barbosa_circumventing_1992}, which penalizes the jump between the multiplier and its physical interpretation.
The Barbosa-Hughes approach was extended by Haslinger and Renard to the fictitious domain method in the XFEM framework~\cite{haslinger_new_2009}.
In the unfitted FEM framework, a different type of a stabilization method was introduced by Burman and Hansbo~\cite{burman_fictitious_2010}, where the multiplier is chosen as a piecewise constant function and the stability of the saddle-point formulation is achieved by penalizing the jump of the multiplier over the elemental faces~\cite{burman_fictitious_2010}.
Also, a primal space with bubble-stabilized basis functions is considered such that they satisfy the discrete inf-sup condition~\cite{mourad_bubble-stabilized_2007,dolbow_residual-free_2008}.
Other approaches, where the primal space was kept untouched and coarser multiplier spaces were considered in several works~\cite{moes_imposing_2006,kim_mortared_2007,bechet_stable_2009}.
B\'echet et al.~developed a stable Lagrange multiplier space based on a vital vertex algorithm~\cite{bechet_stable_2009}, which was later extended by Hautefeuille et al.~\cite{hautefeuille_robust_2012}.
This method does not require any stabilization terms or modification of the primal space, only the multiplier space is designed carefully such that it satisfies the inf-sup condition and ensures optimal convergence of discretization error.
In this work, we employ the vital vertex algorithm for constructing a stable Lagrange multiplier space.

In the unfitted FE framework, the background mesh and the computational domains are independent entities, hence the elements associated with the background mesh are allowed to intersect arbitrarily.
Due to this reason, the linear system of equations arising from the unfitted FE discretization can be highly ill-conditioned.
Additionally, the system of equations arising from the method of Lagrange multipliers can also be formulated as a quadratic minimization problem with linear inequality constraints.
To this end, we propose a tailored generalized multigrid method for solving contact problems in the unfitted discretization method, where the non-penetration condition is discretized with the method of Lagrange multipliers.
Our generalized multigrid method utilizes the pseudo-$L^2$-projection to compute the transfer operator which was originally proposed in~\cite{kothari_multigrid_2019} and utilized in~\cite{kothari_multigrid_2021}.
The generalized multigrid method is motivated by the monotone multigrid method proposed in~\cite{kornhuber_monotone_1994,kornhuber_monotone_1996}, which was developed for solving a quadratic minimization problem with pointwise constraints arising from the variational inequalities.
In the monotone multigrid method, the energy functional is minimized successively such that each iterate satisfies the constraints.
An important component of the monotone multigrid is projected Gauss-Seidel (PGS) smoother, which can simultaneously minimize the energy functional and project the current iterate onto a feasible set in each local iteration.
For the linear constraints, which are represented by a linear combination of several variables, the traditional PGS method is unusable.
To overcome this difficulty, we introduce the orthogonal transformation of the linear constraints and propose a novel variant of the PGS method that can handle the linear constraints.

The outline of this paper is given as follows.
We introduce the two-body contact problem as a model problem and discuss the unfitted FE discretization in detail in Section~\ref{sec:model}.
In Section~\ref{sec:QRMG}, we introduce our generalized multigrid method and explain each component of the multigrid method in detail.
We discuss the orthogonalization strategy to decouple the linear constraints and introduce the modified PGS method to tackle the decouple constraints.
Lastly, in Section~\ref{sec:results} we present the results of numerical experiments.
We study the discretization error and the performance of the multigrid method with respect to several parameters.
We show the robustness of the proposed generalized multigrid method for Signorini's problem and the two-body contact problem and also compare the performance of our multigrid method with other solution strategies.
 \section{Two-body Contact Problem in XFEM Framework} \label{sec:model}
In this section, we introduce the two-body contact problem within the unfitted finite element framework.
Here, we assume that the contact between the two bodies takes place on an embedded interface.
The Dirichlet boundary and Neumann boundaries can also be assumed to be embedded, but in order to simplify the presentation of the problem, these boundaries are assumed to be fitted with the background mesh.

\subsection{Problem Description}
We assume two elastic bodies $\Omega^1,\Omega^2 \in \bR^d$, $d\in\{2,3\}$, with the Lipschitz continuous boundaries $\Gamma^1,\Gamma^2$.
The bodies are assumed to be subjected to volume forces $\vecf:\Omega \to \bR^d$ and the traction/surface forces on the Neumann boundary ${\vect_N:\Gamma_N \to \bR^d}$.
Both bodies undergo deformation due to the influence of these external forces.
A material point $\matX \in \Omega$ in the undeformed state moves to the location $\matX + \vecu$ after the deformation.
Here, the vector-valued quantity $\vecu:\Omega \to \bR^d$ denotes the displacement field of the material point $\matX$, denoted as $\vecu := \vecu(\matX)$.
The boundary $\Gamma$ is decomposed into three parts: the Dirichlet boundary $\Gamma_D$, the Neumann boundary $\Gamma_N$, and apriori unknown contact boundary $\Gamma_c$.

In elastostatics, the displacement field $\vecu :=(\vecu^1,\vecu^2)$ can be given as a solution of the following boundary value problem:
\begin{equation}
  \begin{aligned}
    -\nabla \cdot \bs{\sigma} & = \vecf   \quad &  & \text{in } \Omega^i ,    \\
    \vecu                     & = \vecu_D       &  & \text{on } \Gamma_D^i, \\
    \bs{\sigma} \cdot \vecn^i & = \vect_N       &  & \text{on } \Gamma_N^i,
  \end{aligned}
  \label{eq:two_body}
\end{equation}
where $\vecn^i$ denotes the outward normal on the Neumann boundary $\Gamma_N^i$.
We note that unlike Dirichlet and Neumann boundaries the contact boundary is shared between both bodies, $\Gamma_c = \Gamma_c^1 = \Gamma_c^2$.
These parts of the boundaries are assumed to be disjoint and the contact boundary is assumed to have a positive measure, i.e., $\meas(\Gamma_c) > 0$.
In \eqref{eq:two_body}, we denote the Cauchy stress tensor as $\bs{\sigma}:= \bs{\sigma}(\vecu)$.
Here, we assume the bodies $\Omega^1,\Omega^2$ to be linear elastic, where the constitutive law is provided by Hooke's law
\[
  \bs{\sigma} = \lambda \text{tr}(\bs{\varepsilon}) \matI + 2 \mu \bs{\varepsilon},
  \label{eq:contitutive_law}
\]
where $\lambda$ and $\mu$ are the Lam\'e parameters, $\text{tr}(\cdot)$ denotes the trace operator, $\matI$ is second order identity tensor, and the linearized strain tensor $\bs{\varepsilon} := \bs{\varepsilon}(\vecu)$ is defined as
\(
  \bs{\varepsilon}(\vecu) := \frac{1}{2}\big( \nabla \vecu +(\nabla \vecu)^\T \big).
  \label{eq:strain_tensor}
\)

We assume that the initial gap function $g_c:\Gamma_c \to \bR^+$ is given between two bodies in the direction of outward normal $\vecn$, where the outward normal is defined as $\vecn = \vecn^1 = -\vecn^2$.
The point-wise gap in the displacement fields for both domains is defined as
\[
  \llbracket \vecu \cdot \vecn \rrbracket := \vecu^1\cdot\vecn^1 + \vecu^2\cdot \vecn^2 = (\vecu^1-\vecu^2)\cdot \vecn^1.
\]
We define a gap function $g_c$ as the distance from the possible contact boundary of one body to the other body in the normal direction.
The non-penetration condition on the possible contact boundary $\Gamma_c$ is given as in \eqref{eq:1}.
The contact pressure or stress developed in the normal direction on $\Gamma_c$ is compressive \eqref{eq:2}. We also decompose the traction vector at the contact boundary into the normal and tangential components, given as ${\bs{\sigma}\vecn = \sigma_n\cdot\vecn + \bs{\sigma}_t}$, where $\sigma_n = \vecn\cdot \bs{\sigma} \vecn$.
The third contact condition is given as complementarity condition, given as in \eqref{eq:3}, which ensures that the gap between the body and the rigid obstacle is zero in presence of non-zero contact pressure and the contact pressure is zero in absence of contact.
As we are considering the frictionless contact problem, the body is allowed to move freely in the tangential direction and the induced tangential stresses are given as in \eqref{eq:4}.
The frictionless linearized contact conditions are given as follows:
\begin{linenomath}
  \begin{subequations}
    \begin{align}
      \llbracket \vecu\cdot \vecn \rrbracket - g_c           & \leqslant 0   \quad \text{ on } \Gamma_c, \label{eq:1} \\
      \sigma_n                                               & \leqslant 0   \quad \text{ on } \Gamma_c, \label{eq:2} \\
      (\llbracket \vecu\cdot \vecn \rrbracket - g_c)\sigma_n & = 0           \quad \text{ on } \Gamma_c, \label{eq:3} \\
      \bs{\sigma}_t                                          & = 0           \quad \text{ on } \Gamma_c. \label{eq:4}\end{align}\label{eq:contact_two_body_cond}\end{subequations}\end{linenomath}In contact mechanics these conditions are known as Hertz–Signorini–Moreau conditions for frictionless contact, while in optimization literature they are known as Karush–Kuhn–Tucker (KKT) conditions of the constraints.

\begin{remark}
  In this work, we also consider Signorini's contact problem in the unfitted FE framework.
  In Signorini's problem, a contact between a linear elastic body and a rigid foundation is considered and the gap function is computed as a distance from apriori unknown contact boundary to the rigid foundation.
  The non-penetration condition for this problem is given as, $ \vecu\cdot \vecn - g_c \leqslant 0 $ on $\Gamma_c$.
\end{remark}

\subsection{Unfitted FE Discretization}
In this section, we discuss the discretization of the two-body contact problem. For simplicity, we assume that only the contact boundary is not fitted with the mesh, while Dirichlet and Neumann boundaries are fitted.

We assume a shape regular, quasi-uniform, conforming quadrilateral mesh $\widetilde{\pT}_h$.
The domain $\Omega=\Omega^1\cup\Omega^2$ is encapsulated by the mesh, $\Omega \subset \widetilde{\pT}_h$. The contact boundary $\Gamma_c$ is assumed to be resolved sufficiently well by the mesh $\widetilde{\pT}_h$ and the curvature of the boundary is assumed to be bounded.
Let $h_K$ be the diameter of the element $K$, and mesh size is defined as $h=\max_{K\in \widetilde{\pT}_h} h_K$.
We define an active mesh, which is strictly intersected by the domains $\Omega^1,\Omega^2$ as
\[
  \pT^i_{h} = \{ K\in \widetilde{\pT}_h: K\cap \Omega^i  \neq \emptyset\}, \quad  i \in \{1,2\}.
\]
Here, each domain $\Omega^i$ is captured by the respective active mesh, $\Omega^i \subset \pT^i_{h}$, for $i=\{1,2\}$.
The active meshes exclude all the elements that are neither intersected by the boundary $\Gamma_c$ nor are in the interior of the domain.
We define a set of elements that are intersected by the contact boundary $\Gamma_c$ as
\[
  \pT_{h,\Gamma_c}=\{K \in \widetilde{\pT}_h: K\cap \Gamma_c \neq \emptyset \}.
\]
For all elements $K \in \pT_{h,\Gamma_c}$, let $K_{\Omega} := K \cap \Omega $ be part of $K$ in domain $\Omega$.
The elements $K \in  \pT_h \setminus \pT_{h,\Gamma_c}$ are strictly in the interior of domain $\Omega$.
For all $K \in \pT_{h,\Gamma_c}$, let $\Gamma_K:=\Gamma_c \cap K$ be part of $\Gamma_c$ in $K$.

We define a continuous FE space over the mesh $\widetilde{\pT}_h$ as
\begin{equation}
  \label{eq:background_FEspace}
  \widetilde{\bpV}_h = \{\vecv \in [H^1 (\widetilde{\pT}_h)]^d: \vecv\vert_K \in \pQ_1(K), \vecv\vert_{(\partial \widetilde{\pT}_h)_D} = 0, \forall K\in \widetilde{\pT}_h \},
\end{equation}
where $\pQ_1$ denotes the space of piecewise bilinear functions.
Following the original XFEM literature~\cite{moes_finite_1999}, we define a characteristic function of the computational domains $\Omega^i$ for $i=\{1,2\}$, as
\begin{equation}
  \label{eq:Heaviside}
  \chi_{\Omega^i} : \bR^d \to \bR, \qquad
  \chi_{\Omega^i}(\matX) = \begin{cases}
    1 & \quad \forall \matX \in \overbar{\Omega}^i, \\
    0 & \quad \text{otherwise.}
  \end{cases}
\end{equation}
The function space $\widetilde{\bpV}_h$ is spanned by the nodal basis functions $\widetilde{\Phi}_h = (\widetilde{\phi}^p_h)_{p\in \widetilde{\mathcal{N}}_h}$, where $\widetilde{\mathcal{N}}_h$ denotes the set of nodes of the background mesh $\widetilde{\pT}_h$.
The characteristic function is used to restrict the support of the finite element space $\widetilde{\bpV}_h$ to the respective domain $\Omega^i$ thus
\({
  \bpV^i_{h} = \chi_{\Omega^i}(\matX)\widetilde{\bpV}_h.
}\)
We seek the approximation $\vecu_h = (\vecu^1_h \oplus \vecu^2_h)$ in space $\bpV_h = \bpV^1_{h} \oplus \bpV^2_{h}$.
We define the set of nodes on the active mesh $\pT_{h, i}$ associated with domain $\Omega_i$ as
\[
  \mathcal{N}_{h,i} := \{ p \in \widetilde{\mathcal{N}}_h: \supp(\widetilde{\phi}^p_h) \cap \Omega_i \neq \emptyset\} \quad i=\{1,2\}.
\]
We now define the ``cut'' basis function associated with a node $p$ as
\[
  \phi^p_h = \chi_{\Omega_i}(\matX) \widetilde{\phi}^p_h \qquad \forall p \in \mathcal{N}^i_{h},\ i=\{1,2\}.
\]
The function space $\bpV^i_{h}$ is spanned by the nodal basis functions $\Phi^i_{h} = ({\phi}^p_h)_{p\in \mathcal{N}^i_h}$.
We define the span of nodal basis function on $\bpV_h$ as $\Phi_h = \Phi^1_{h}\oplus\Phi^2_{h}$, and the set of nodes associated with the mesh $\pT_h$ is given by $\mathcal{N}_h = \mathcal{N}^1_{h} \oplus \mathcal{N}^2_{h}.$

\begin{figure*}
  \begin{subfigure}[t]{.23\textwidth}
    \centering
    \includegraphics{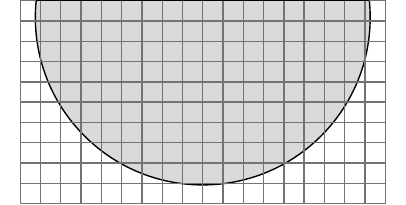}
\caption{A quadrilateral mesh $\widetilde{\pT}_h$ is used as a background mesh to capture a circular domain $\Omega$.}
\end{subfigure}\hfill
  \begin{subfigure}[t]{.23\textwidth}
    \centering
    \includegraphics{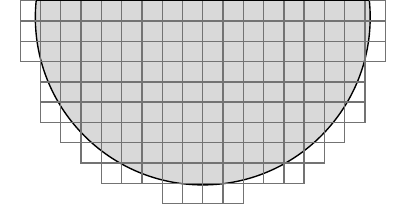}
\caption{The mesh $\pT_h$ is strictly intersected by the domain $\Omega$.}
\end{subfigure}
  \hfill
  \begin{subfigure}[t]{.23\textwidth}
    \centering
    \includegraphics{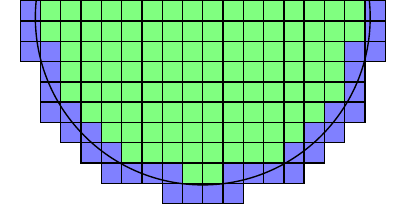}
\caption{The mesh $\pT_{h,\Gamma_c}$ is shaded in blue, while the interior mesh $\pT_h\backslash \pT_{h,\Gamma_c}$ is shaded in green.}
\end{subfigure}\hfill
  \begin{subfigure}[t]{.23\textwidth}
    \centering
    \includegraphics{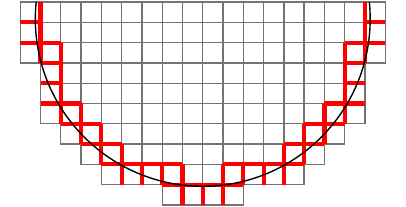}
\caption{ The set of faces $\pG_{h,\Gamma_c}$ is shown in the red.}
    \label{fig:EdgesForGhostPenalty}
  \end{subfigure}
  \caption{An example of a domain $\Omega$ with a background mesh $\widetilde{\pT}_h$.}
\label{fig:mesh_fictitious_domain}
\end{figure*}

\begin{remark}
  For Signorini's problem, the characteristic function \eqref{eq:Heaviside} is defined only on the domain $\Omega^1$.
  An example of a domain $\Omega$ embedded in the background mesh $\widetilde{\pT}_h$, the active mesh $\pT_h$, the interior mesh, and the cut mesh are shown in Figure~\ref{fig:mesh_fictitious_domain}.
\end{remark}

\subsection{Variational Formulation}
The variational formulation of the two-body contact problem using the principle of virtual work is given as:
\begin{equation}
  \text{find } \vecu_h \in \bpK_h \text{ such that, } \quad  a(\vecu_h,\vecv_h-\vecu_h) \geqslant F(\vecv_h-\vecu_h) \quad \forall \vecv_h \in \bpK_h,
  \label{eq:abstract_variational}
\end{equation}
where $a(\cdot,\cdot): \bpV_h \times \bpV_h \to \bR$ is a symmetric, continuous and coercive bilinear form, and ${F(\cdot):\bpV_h \to \bR}$ denotes continuous and bounded linear form.
The bilinear and the linear forms are defined as
\begin{equation}
  \begin{aligned}
    \label{eq:2bodyelast_bilinear_linear}
    a(\vecu_h,\vecv_h) = \sum_{i=\{1,2\}} \int_\Omega \bs{\sigma}(\vecu_h^i):\bs{\varepsilon}({\vecv_h^i}) d\Omega,  \qquad \qquad  F(\vecv_h) = \sum_{i=\{1,2\}} \Bigg( \int_\Omega \vecf \vecv_h^i d\Omega + \int_{\Gamma_N} \vect_N \vecv_h^i d\Gamma \Bigg).
  \end{aligned}
\end{equation}
Traditionally, in the fitted finite element framework, we can use the method of Lagrange multipliers, the penalty method, Nitsche's method, the regularization methods, the augmented Lagrangian methods, etc., to impose the contact conditions~\cite{wohlmuth_variationally_2011}.
Here, in the unfitted FE framework, we employ the method of Lagrange multipliers to enforce the non-penetration contact conditions \eqref{eq:1}.
This is due to the fact that the Lagrange multiplier formulation does not require modification of the primal formulation, and the contact condition can be handled by the multipliers implicitly.
Whereas, Nitsche's formulation for the contact problem is more complex as we have to handle the non-penetration conditions in the primal formulation.

We impose the non-penetration contact condition using the method of multipliers, where the multiplier space is constructed by employing the vital vertex algorithm~\cite{bechet_stable_2009,hautefeuille_robust_2012}.
We introduce the multiplier space $\pM_h \subseteq H^{-\half}(\Gamma_c)$ and define the bilinear form $b(\cdot,\cdot):\pM_h\times \bpV_h \to \bR$, given as
\begin{equation}
  b(\mu_h,\vecu_h):= \sum_{K \in \pT_{h,\Gamma_c}} \int_{\Gamma_K} {\mu_h} \llbracket \vecu_h \cdot \vecn \rrbracket \ d \Gamma \qquad \forall \mu_h \in \pM_h,\forall \vecu_h \in \bpV_h,
  \label{eq:two_body_b}
\end{equation}
and the linear form $G(\cdot):\pM_h \to \bR$, given as
\[
  G(\mu_h) :=  \sum_{K \in \pT_{h,\Gamma_c}} \int_{\Gamma_K} \mu_h g_c \ d {\Gamma} \qquad \forall \mu_h \in \pM_h.
\]
Finally, we can define the space of admissible displacements, such that they satisfy the contact conditions as
\begin{equation}
  \bpK_h:=\{\vecv_h \in \bpV_h : b(\mu_h, \vecv_h) \leqslant G(\mu_h), \forall \mu_h\in \pM_h \}.
  \label{eq:admissible_space}
\end{equation}
The space of admissible displacements $\bpK_h$ is a closed convex subset of FE space $\bpV_h$.
Due to the inequality condition in \eqref{eq:abstract_variational}, the contact problem is inherently non-linear.

\subsubsection{Ghost Penalty Stabilization}\label{sec:ghost_penalty}
In unfitted methods, a background mesh captures the computational domain of arbitrary shape, hence the elements are allowed to be cut arbitrarily by the boundary/interface.
In general, this flexibility can result in disproportionally cut elements, which might not be shape regular anymore.
For this reason, the bound on the gradient of a function can become arbitrarily weak for the unfortunately cut elements.
By adding the ghost penalty term~\cite{burman_ghost_2010}, we regain control over the gradients of the function on the cut elements with very small support, and by extension, we can overcome the issue of ill-conditioning.
This stabilization term has to be chosen in such a way that it provides sufficient stability and stays weakly consistent with the original formulation for smooth functions.
We define the set of faces $\pG_{h,\Gamma_c}^i$ for each subdomain $\Omega^i$ as
\[
  \pG_{h,\Gamma_c}^i =\{ G \subset \partial K \mid K\in \pT^i_{h,\Gamma_c},\ \partial K \cap \partial \pT_i = \emptyset\} \quad i=\{1,2\}.
\]
An example of the set of faces $\pG_{h,\Gamma_c}$ in the context of Signorini's problem can be seen in Figure~\ref{fig:EdgesForGhostPenalty}.
The ghost penalty term is enforced on the set of edges $\pG_{h,\Gamma_c}$, and it is defined as
\begin{equation}
  j(\vecu_h,\vecv_h) = \sum_{i=\{1,2\}}\sum_{G\in \pG^i_{h,\Gamma_c}} \int_G \epsilon_G {h_G} (2\mu^i+\lambda^i) \llbracket \nabla \pE_h\vecu_h \cdot \vecn_G \rrbracket \llbracket \nabla\pE_{h} \vecv_h \cdot \vecn_G \rrbracket\ dG,
\end{equation}
where $h_G$ is the diameter of the face $G$, $\vecn_G$ denotes unit normal to face $G$, $\epsilon_G$ is a positive constant and $\lambda^i, \mu^i$ denote the Lam\'{e} parameters associated with the either domain~\cite{sticko_high-order_2020}.
Here, $\pE_h$ denotes the canonical extension of the function from the domain to the background mesh, which is defined as $\pE_h:\bpV_h\vert_{K_\Omega} \to \widetilde{\bpV}_h\vert_K$.
This ghost penalty term is enforced in the normal derivatives of the displacement field.
A different approach is also considered in~\cite{claus_stable_2018}, where the ghost penalty term is enforced in the normal derivatives of the stress field.

Now, we modify the variational formulation of the two-body contact problem by adding the ghost penalty stabilization term.
The updated variational problem is defined as:
\begin{equation}
  \text{find } \vecu_h \in \bpK_h \text{ such that }\quad  a_j(\vecu_h,\vecv_h-\vecu_h) \geqslant F(\vecv_h-\vecu_h) \quad \forall \vecv \in \bpK_h,
  \label{eq:abstract_variational_2}
\end{equation}
with the bilinear form $a_j(\vecu_h,\vecv_h) := a(\vecu_h,\vecv_h) + j(\vecu_h,\vecv_h)$.

\subsubsection{Discretization of Non-penetration Condition}
We have used the method of Lagrange multipliers to discretize and enforce the non-penetration contact condition.
In order to achieve optimal convergence rates of the discretization method, the choice of the FE spaces for primal variable $\vecu_h \in \bpV_h$ and the dual variable $\lambda_h \in \pM_h$ is crucial.
In the unfitted FE framework, the most convenient options for $\bpV_h$ and $\pM_h$ are very rarely stable.
As, the method of multipliers is stable only if the following discrete inf-sup condition is satisfied
\begin{equation}
  \label{eq:infsup}
  \inf_{\mu_h \in \pM_h} \sup_{\vecv_h \in \bpV_h} \frac{b(\mu_h,\vecv_h)}{\|\vecv_h\|_{\bpV_h} \|\mu_h\|_{\pM_h} } \geqslant \beta > 0,
\end{equation}
where the constant $\beta$ does not dependent on mesh-size $h$.
If the inf-sup condition~\eqref{eq:infsup} is not satisfied, it can give rise to spurious modes in the discrete Lagrange multiplier space.
The effect of the spurious modes in the solution can be observed as locking phenomena on the interface.
In the unfitted FE framework, it is not trivial to create an optimal multiplier space, as the interface is not resolved by the background mesh.

We employ the vital vertex algorithm to create a stable multiplier space, where the primal space is kept the same and a coarser multiplier space is chosen.
Such approach of constructing a coarser multiplier space can be found in several works in the context of unfitted FE framework~\cite{bechet_stable_2009,kim_mortared_2007,moes_imposing_2006}.
B\'echet et al.~developed a stable multiplier space based on a vital vertex algorithm~\cite{bechet_stable_2009}. This method does not require any stabilization terms and also the primal space $\bpV_h$ is not modified.
Only the multiplier space $\pM_h$ is designed carefully such that it satisfies the inf-sup condition.
We denote the list of vital vertices on the interfaces as $\mV_{h,\Gamma_c}$, and the dimension of the multiplier space $\pM_h$ is given as $\lvert\mV_{h,\Gamma_c}\rvert$.
For each vital vertex $p\in \mV_{h,\Gamma_c}$, we define the associated basis function $\phi_h^q, \forall q \in \mN_h$ restricted to the interface $\Gamma_c$.
Also, we introduce a set of nodes $\mN_{h,\Gamma_c}$, given as
\[
  \mN_{h,\Gamma_c}:=\{q \in \mN_h : \phi_h^q\vert_{\Gamma_c} \neq 0\},
\]
where $\mN_{h,\Gamma_c}$ includes all nodes that are endpoints of the cut-edges.
The set of nodes $\mN_{h,\Gamma_c}$ is later divided into a set of active nodes $\mN^A_{h,\Gamma_c}$ and inactive nodes $\mN^I_{h,\Gamma_c}$.
Here, the set of active nodes $\mN^A_{h,\Gamma_c}$ are defined as the endpoints of the edges on which the vital vertices are located and the inactive nodes are given as $ \mN^I_{h,\Gamma_c} =\mN_{h,\Gamma_c} \setminus \mN^A_{h,\Gamma_c}$.

The vital vertex method is ideal for the unfitted FEM as basis functions for the multiplier space are still constructed as a trace of the higher dimensional basis functions.
The nodal basis functions for each vital vertex $p\in \mV_{h,\Gamma_c}$ is given as
\[
  \psi_h^p := \sum_{p\in \mN_{h,\Gamma_c}} w_{pq} \phi_h^q\vert_{\Gamma_c} \quad \forall q\in \mN_{h,\Gamma_c},
\]
where $w_{pq}$ are coefficients of the linear combination.
\begin{figure*}
  \begin{subfigure}[t]{.23\textwidth}
    \centering
    \includegraphics{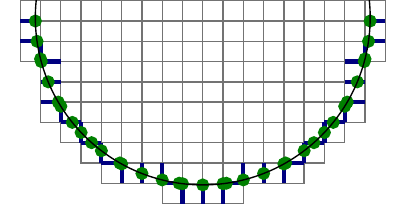}
\caption{All vertices are marked with a green circle, all cut edges which are intersected by the boundary are marked with blue line.}
    \label{fig:allvertices}
  \end{subfigure}\hfill
  \begin{subfigure}[t]{.23\textwidth}
    \centering
    \includegraphics{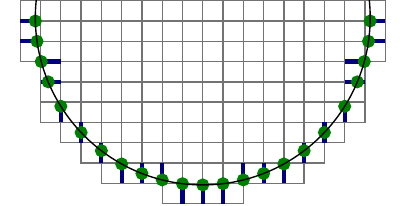}
\caption{Vital vertices are marked with a green circle, the cut edges on which the vital vertices lie are marked with blue line.}
    \label{fig:vitalvertices}
  \end{subfigure}
\hfill
  \begin{subfigure}[t]{.23\textwidth}
    \centering
    \includegraphics{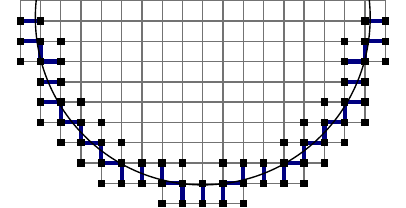}
\caption{The nodes $\mathcal{N}_{h,\Gamma_c}$ that belong to the sub-mesh $\pT_{h,\Gamma_c}$ are marked with black squares}
    \label{fig:allnodes}
  \end{subfigure}\hfill
  \begin{subfigure}[t]{.23\textwidth}
    \centering
    \includegraphics{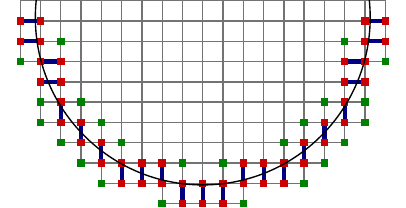}
\caption{All active nodes $\mathcal{N}^A_{h,\Gamma_c}$ are marked with red squares, and the inactive nodes $\mathcal{N}^I_{h,\Gamma_c}$ are marked with green squares}
    \label{fig:activeinactivenodes}
  \end{subfigure}
  \caption{Different type of nodes characterized by the vital-vertex algorithm.}
  \label{fig:vitalvertex_info}
\end{figure*}
In Figure~\ref{fig:vitalvertex_info}, we can see the set of all vertices (Figure~\ref{fig:allvertices}), a set of vital vertices $\mV_{h,\Gamma_c}$ (Figure~\ref{fig:vitalvertices}), the set of all nodes associated with the cut elements $\mN_{h,\Gamma_c}$ (Figure~\ref{fig:allnodes}), and an identified set of active nodes $\mN^A_{h,\Gamma_c}$ and inactive nodes $\mN^I_{h,\Gamma_c}$ (Figure~\ref{fig:activeinactivenodes}) for Signorini's problem.

\subsection{Local Basis Transformation}\label{sec:householder_reflection}
We note, the non-penetration constraints in the contact problem are given by the relative displacement of the bodies in the normal direction.
Thus, the constraint at any node is given by the coupling of the degrees of freedom (DoFs) on the node.
To create the constraint matrix such that the non-penetration condition is enforced only on one DoF per node, we transform the system into a new basis.

Let $\{\matE_i\}_{i=1,\ldots,d}$ be the Euclidean basis of $\bR^d$ and $\vecn_p$ be the outward normal on the node $p$.
On each node $p\in \mN_{h,\Gamma_c}$, we define a new basis $\vece_1(p) = \vecn_p$ and also change $\{\vece_i\}_{i=2,\ldots,d}$ such that these redefined basis are also orthonormal, while for all $q \in \mN_h \setminus \mN_{h,\Gamma_c}$ the definition of the Euclidean basis remain same.
This approach was introduced for Signorini's problems~\cite{krause_monotone_2001} and has been effectively applied to multi-body contact problems~\cite{dickopf_efficient_2009,krause_nonsmooth_2009,wohlmuth_monotone_2003}.
The transformed basis are constructed using a local Householder transformation on $\bR^d$, given as
\[
  \matO_{pp} = \matI - 2(\vecw_p \otimes \vecw_p) \qquad \forall p \in \mN_{h,\Gamma_c},
\]
where the vector $\vecw_p$ is computed by $\vecw_p = (\vecn_p - \matE_1)/\| \vecn_p - \matE_1\|_2$.
Now, due to the Householder transformation for all $p \in \mN_{h,\Gamma_c}$ we can uniquely define the local unit vectors as $\vece_i(p) = \matO_{pp} \matE_i(p)$.
While, we define $\matO_{qq} = \matI$ for all $q \in \mN_h \setminus \mN_{h,\Gamma_c}$, which ensure that the basis system on those nodes remain unchanged, i.e., $\vece_i(q) = \matE_i(q)$.
Thus, by using these local transformation matrices, we can construct the global matrix $\matO\in \bR^{nd\times nd}$ where $nd = |\mN_h|\cdot d$ and $\matO = \oplus_{p\in \mN_h} \matO_{pp}$, which is an orthonormal matrix with the properties, $\matO \matO^\T = \matO^\T \matO = \matI$ and $\matO = \matO^\T$.

This transformation decouples and locally modifies the constraints and it is only applicable in the normal direction.
The bilinear form \eqref{eq:two_body_b} for the two-body contact problem  can be reformulated as
\begin{equation}
  \label{eq:two_body_b_hh}
  b(\mu_h,\vecu_h) := \sum_{K\in \pT_{h,\Gamma_c}}  \langle \mu_h, \llbracket \vecu_h \cdot \matE_1 \rrbracket  \rangle_{\Gamma_K}  \quad \forall \mu_h\in \pM_h, \forall \vecu_h \in \bpV_h.
\end{equation}

In the next section, we discuss the algebraic formulation of the contact problem \eqref{eq:abstract_variational_2} and also discuss the effect of the local basis transformation algebraically.

\subsection{Algebraic Formulation}
The abstract variational problem for the contact problem can be reformulated as an optimization problem with inequality constraints, given as
\begin{equation}
  \begin{aligned}
    \label{eq:min_problem}
     & \min_{\vecu_h\in \bpV_h} \pJ(\vecu) = \half a_j(\vecu_h,\vecu_h) - F(\vecu_h)                 \\
     & \text{subject to } \quad b(\mu_h,\vecu_h) \leqslant G(\mu_h) \quad \forall \mu_h \in \pM_h.
  \end{aligned}
\end{equation}
The above minimization problem can be written in an algebraic formulation using the local basis transformation introduced in the previous section.
We denote the local entries of the stiffness matrix and the right-hand side as
\begin{equation}
  \label{eq:stiff_b}
  A_{pq} = (a_j(\phi^p_h \matE_i, \phi^q_h \matE_k))_{i,k=1,\ldots,d}, \quad b_{p} = (F(\phi^p_h \matE_i))_{i=1,\ldots,d}.
\end{equation}
The global stiffness matrix and the right-hand side vector can be assembled as
\[
  \matA = (A_{pq})_{p,q \in \mN_h}, \quad \vecb = (b_p)_{p\in \mN_h}.
\]
The bilinear form $b(\cdot,\cdot)$ can be decomposed into two parts, ${b^1(\cdot,\cdot): \pM_h \times \bpV^1_{h} \to \bR }$  and ${b^2(\cdot,\cdot): \pM_h \times \bpV^2_{h} \to \bR}$ associated with domain $\Omega^1$ and $\Omega^2$, respectively, written as
\[
  b(\mu_h,\vecu_h) = b^1(\mu_h,\vecu^1_h) - b^2(\mu_h,\vecu^2_h).
\]
This decomposition allows us to write the local entries of the constraint matrix as
\begin{equation}
  \label{eq:local_b}
  B^1_{rk} = (b^1(\psi^r_h, \phi^k_h \matE_i))_{i=1,\cdots,d} \quad B^2_{rl} = (b^2(\psi^r_h, \phi_h^l \matE_i))_{i=1,\cdots,d}, \quad g_r = G(\psi^r_h),
\end{equation}
where $\psi_h^r, \phi^k_h$ and $\phi_h^l$ denote the basis functions associated with nodes $r,k,l$ in the FE spaces $\pM_h$, $\bpV^1_h$ and $\bpV^2_h$, respectively.
Thus, the entries of the constraint matrix $\matB$ and the gap vector $\vecg$ are given as
\[
  \begin{aligned}
    \matB^1 = (B^1_{rk})_{r\in \mV_{h,\Gamma_c},k\in\mN^1_h}, \quad  \matB^2 = (B^2_{rl})_{r\in \mV_{h,\Gamma_c}, l\in \mN^2_h}, \quad \vecg = (g_r)_{r\in \mV_{h,\Gamma_c}},
  \end{aligned}
\]
where $\mN_h^1$ and $\mN_h^2$ denote the set of nodes of the active meshes associated with each body, respectively, and $\mV_{h,\Gamma_c}$ denotes the set of vital vertices.
The matrix $\matB$ can be constructed as $\matB = [\matB^1 -\matB^2]$.
Now, we can write the algebraic formulation of the constraint minimization problem~\eqref{eq:min_problem} as
\begin{equation}
  \begin{aligned}
     & \min_{\vecx \in \bR^{nd}} J(\vecx) = \half \vecx^\T\matA\vecx  - \vecx^\T \vecb \\
     & \text{subject to } \quad \matB \vecx  \leqslant \vecg,
  \end{aligned}
  \label{eq:const_min_mat}
\end{equation}
where \(\vecx,\vecb \in \bR^{nd}\), \(\matA \in \bR^{nd \times nd}\), \(\matB \in \bR^{m \times nd}, \vecg \in \bR^m, m \ll nd \) and $\text{rank}(\matB)=m$.
Here, $\vecx$ denotes the unknown displacements, and $\vecg$ denotes the gap between two bodies on the contact boundary.

As we have changed the definition of the Euclidean basis (in Section~\ref{sec:householder_reflection}), we have to also modify the problem algebraically.
The matrix $\matO$ can be used to transform the variables into the new basis as $\overbar{\vecx} = \matO\vecx$ and $\vecx = \matO \overbar{\vecx}$.
Similarly, the stiffness matrix $\matA$ and the right hand side $\vecb$ in the new basis are given as $\overbar{\matA} = \matO \matA \matO$ and $\overbar{\vecb} = \matO \vecb$.
The constraint matrix can be written in the new basis as $\overbar{\matB} = \matB \matO$, the matrix $\matB$ is obtained by discretization of \eqref{eq:two_body_b} and the matrix $\overbar{\matB}$ is constructed by discretization of \eqref{eq:two_body_b_hh}.
The algebraic formulation of the contact problem~\eqref{eq:const_min_mat} in the new basis system is given as following minimization problem:
\begin{equation}
  \begin{aligned}
    \label{eq:contact_problem_bar}
     & \min_{\overbar{\vecx} \in \bR^{nd}} J(\overbar{\vecx}) = \half \overbar{\vecx}^\T \overbar{\matA} \overbar{\vecx} - \overbar{\vecx}^\T \overbar{\vecb} \\
     & \text{subject to } \quad \overbar{\matB} \overbar{\vecx} \leqslant \vecg.
  \end{aligned}
\end{equation}

\begin{remark}
  The local basis transformation for this problem can be carried out directly during the assembly process.
  We can compute the stiffness matrix, constraint matrix and the right-hand side with locally transformed basis by directly utilizing the new basis $\vece_i$ instead of the Euclidean basis $\matE_i$ in \eqref{eq:stiff_b} and \eqref{eq:local_b}.
\end{remark}

 \section{A Generalized Multigrid Method} \label{sec:QRMG}
In this section, we introduce a new generalized multigrid method for solving a quadratic minimization problem with linear constraints \eqref{eq:contact_problem_bar}.
This multigrid method is motivated by the monotone multigrid method~\cite{kornhuber_monotone_1994,kornhuber_monotone_1996}, which was originally developed to solve a quadratic minimization problem with point-wise inequality constraints.
Here, we present an extension of this method for solving a quadratic minimization problem with linear inequality constraints.

The monotone multigrid method is an iterative method, where within each iteration the energy functional is minimized successively such that the current iterate satisfies the constraints.
This task is carried out using the PGS method, which simultaneously minimizes the energy functional and projects the current iterate onto a feasible set.
The traditional PGS method is unable to tackle the linearly constrained minimization problem, which represents a linear combination of several variables.
To overcome this difficulty, we introduce an orthogonal transformation and a variant of the PGS method that can handle the linear inequality constraints.
In addition, in this multigrid method, we employ the transfer operators constructed using the pseudo-$L^2$-projections.
Here, we introduce the necessary ingredients used in our generalized multigrid method. \subsection{Standard Multigrid Method}

The multigrid method is an ideal iterative method for solving many large-scale linear systems of equations that arise from the discretization of elliptic differential equations~\cite{wolfganghackbusch1986-01-14}.
This method obtains optimal convergence rates by exploiting discretizations with different mesh sizes.
The multigrid method is considered to have optimal complexity as its convergence rate is bounded from above and it does not depend on the size of the problem.
Even though the convergence rate of the multigrid method does not depend on the problem size, the number of arithmetic operations grows proportionally with the problem size.
Hence, the complexity of the multigrid method is given as $\pO(n)$.
The robustness of the multigrid method depends on a sophisticated combination of smoothing iterations and coarse-level corrections.
These components are complementary to each other and reduce the error in a different part of the spectrum.

The multigrid method requires a hierarchy of nested FE spaces and appropriate transfer operators to pass the information between the FE spaces on subsequent levels.
For the geometric multigrid method, a hierarchy of nested meshes is employed, which naturally gives rise to a hierarchy of nested FE spaces.
For nested meshes, we can utilize the standard interpolation operator as a prolongation operator and its adjoint is, traditionally, used as a restriction operator.
As the remaining component of the multigrid method, the basic iterative methods, such as Jacobi method or Gauss-Seidel method, are usually employed as smoothers.
These basic iterative methods can eliminate high-frequency or oscillatory components of the error efficiently, but they are relatively slow in eliminating the low-frequency or smooth components of the error.
Thus, after the application of a few smoothing steps, the oscillatory components of the errors are removed, only smooth components of the error remain.
However, when the low-frequency components of the error are mapped onto a coarse mesh, they appear oscillatory on larger resolution of the coarse mesh.
By employing the recursion process, we employ a smoother on each level.
This allows us to remove the corresponding high-frequency components of the error associated with the respective levels.
While the remaining error components, that can not be removed with the smoothers, are annihilated by a direct solver on the coarsest level.

We define the levels as $\ell \in \{0,\ldots, L\}$, where $\ell=0$ denote the coarsest level and $\ell=L$ denotes the finest level.
In the context of this subsection, we assume that the multigrid method is devised to solve the linear system of the equation such as \({ \matA_L \vecx_L = \vecb_L}\). The corresponding residual on the finest level is given as \({\vecr_L = \vecb_L - \matA_L \vecx_L}\).
The standard multigrid method is described in Algorithm~\ref{algo:mmg_coarse}, where \({\nu_1,\nu_2}\) are the number of pre-smoothing and post-smoothing steps, respectively.
The values of \({\gamma=1}\) and \({\gamma=2}\), in the multigrid algorithm transform the multigrid method to a \({V(\nu_1,\nu_2)}\)-cycle and \({W(\nu_1,\nu_2)}\)-cycle, respectively.
The matrix \({\matT_{\ell-1}^{\ell}}\) denotes the prolongation operator, and its adjoint denotes the restriction operator.
Algorithm~\ref{algo:mmg_coarse} is written in an abstract way, such that it returns the correction \({\vecc_L}\) rather than the iterate explicitly, as the method can also be used as a preconditioner.

\begin{algorithm*}[t]
    \caption{Standard multigrid cycle}
    \label{algo:mmg_coarse}
    \SetKwComment{Comment}{$\triangleright$\ }{}
    \SetKwInOut{Input}{Input}\SetKwInOut{Output}{Output}
    \Input{$(\matA_\ell)_{\ell=0,\dots,L}, \vecr_{L}, L, \nu_1, \nu_2, (\matT_{\ell-1}^\ell)_{\ell=1,\dots,L},\gamma$}
    \Output{$\vecc_{L}$}
    \BlankLine
    Function: MG($\matA_{\ell}, \vecr_\ell, \ell, \nu_1, \nu_2, \matT_{\ell-1}^\ell$,$\gamma$)\\
    \eIf{$\ell \neq 0$}
    {
      $\vecc_{\ell}   \mapsfrom \bs{0}$\Comment*[r]{initialize correction}
      $\vecc_{\ell}   \mapsfrom $ Smoother($\matA_\ell, \vecc_\ell, \vecr_\ell, \nu_1$)\Comment*[r]{pre-smoothing}
      $\vecr_{\ell-1} \mapsfrom (\matT_{\ell-1}^{\ell})^\T (\vecr_\ell - \matA_\ell \vecc_\ell)$\Comment*[r]{restriction}
      $\matA_{\ell-1} \mapsfrom (\matT^\ell_{\ell-1})^\T \matA_\ell \matT^\ell_{\ell-1}$\Comment*[r]{Galerkin projection}
      $\vecc_{\ell-1} \mapsfrom \bs{0} $ \Comment*[r]{initialize coarse level correction}
      \For{$i=1,\ldots,\gamma$}
      {
        $\vecc_{\ell-1} \mapsfrom \vecc_{\ell-1} $+ MG($\matA_{\ell-1}, \vecr_{\ell-1}, \ell-1, \nu_1, \nu_2, \matT_{\ell-2}^{\ell-1},\gamma$)\Comment*[r]{coarse level cycle}
      }
      $\vecc_{\ell}   \mapsfrom \vecc_\ell + \matT_{\ell-1}^{\ell} \vecc_{\ell-1}$\Comment*[r]{prolongation}
      $\vecc_{\ell}   \mapsfrom $ Smoother($\matA_\ell, \vecc_\ell, \vecr_\ell, \nu_2$)\Comment*[r]{post-smoothing}
    }
    {
      $\vecc_0 \mapsfrom \matA_0^{-1} \vecr_0 $\Comment*[r]{direct solver}
    }
\end{algorithm*}

\subsection{Transfer operators for Unfitted FEM}\label{sec:transfer}
It is well-known that the efficiency of the multigrid method depends heavily on the underlying hierarchy of meshes and FE spaces.
In the multigrid method, the multilevel decomposition of the FE space is performed in such a way that the FE spaces associated with coarser levels are subspaces of the FE space associated with the finest level.
This does not necessarily hold for the hierarchy of FE spaces in the unfitted FEM framework.
Here, we briefly present a strategy for constructing a nested hierarchy of FE spaces from a hierarchy of non-nested meshes.

We define a sequence of background meshes, denoted as $\{\widetilde{\pT}_\ell\}_{\ell=0,\ldots,L}$.
We associate the original background mesh on which the problem is defined as the mesh on the finest level, give as, $\widetilde{\pT}_L:=\widetilde{\pT}_h$.
The sequence of meshes is created by either choosing a mesh on the coarsest level $\widetilde{\pT}_0$ and uniformly refining this mesh or uniformly coarsening the mesh on the finest level $\widetilde{\pT}_L$.
Here, we assume that the mesh on each level encapsulates the domain $\Omega$.
Now, we can associate FE spaces $\widetilde{\bpV}_\ell$ to the meshes on each level, in the same way as given in \eqref{eq:background_FEspace}.
If the background meshes $\{\widetilde{\pT}_\ell\}_{\ell=0,\ldots,L}$ are nested then the associate FE spaces are also subspaces of the FE space on finest level, given as \({\widetilde{\bpV}_{\ell-1} \subset \widetilde{\bpV}_{\ell} }\), for all \({\ell=1,\ldots,L}\).
In order to create a hierarchy of meshes for unfitted FEM, the background meshes are enriched, decomposed and they are associated with either of the domains.

\begin{figure*}
  \begin{subfigure}[t]{.32\textwidth}
    \centering
    \includegraphics{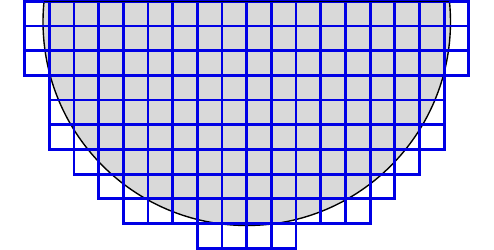}
\caption{Quadrilateral mesh ${\pT}_{\ell-1,i}$.}
  \end{subfigure}\hfill
  \begin{subfigure}[t]{.32\textwidth}
    \centering
    \includegraphics{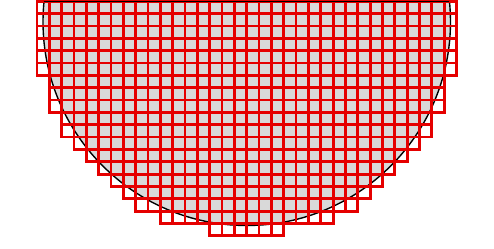}
\caption{Quadrilateral mesh $\pT_{\ell,i}$.}
  \end{subfigure}
  \hfill
  \begin{subfigure}[t]{.32\textwidth}
    \centering
    \includegraphics{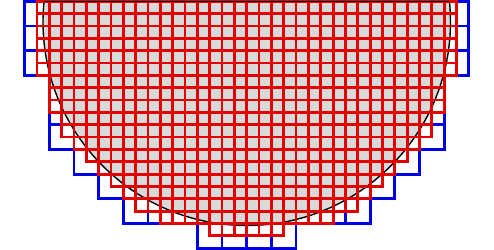}
\caption{Superimposition of $\pT_{\ell,i}$ on $\pT_{\ell-1,i}$.}
  \end{subfigure}\caption{2D Triangular meshes on different levels encapsulating the domain $\Omega_i$, (domain $\Omega_i$ is shaded in gray).}
  \label{fig:multilevel_decomp}
\end{figure*}

In Figure~\ref{fig:multilevel_decomp}, we can see that even though the background meshes are nested, the captured meshes are not necessarily nested.
The nestedness of the captured meshes depends heavily on the embedded interfaces.
Now, utilizing the characteristic function \eqref{eq:Heaviside}, we restrict the support of the FE spaces $\widetilde{\bpV}_\ell$ to the domains, given as $\bpV_\ell^1$ and $\bpV_\ell^2$, respectively.
Hence, the enriched FE spaces associated with a domain are also not nested, i.e., \({\bpV_{\ell-1}^i \not\subset \bpV_{\ell}^i }\), for \(i\in\{1,2\}\) and $\ell\in\{1,\ldots,L\}$.
To create a hierarchy of nested FE spaces from the hierarchy of non-nested meshes, we adopt the variational transfer approach introduced for the unfitted FEM~\cite{kothari_multigrid_2019}.
We define a prolongation operator which projects quantities from a FE space associated with a coarse level to a FE space associated with a fine level, thus as
\[
  \bs{\Pi}_{\ell-1,i}^\ell:\bpV_{\ell-1}^i \to \bpV_{\ell}^i \quad \forall\ \ell \in \{1,\ldots,L\},\ i=\{1,2\},
\]
such that \({ \bs{\Pi}_{\ell-1,i}^\ell \bpV_{\ell-1}^i \subset \bpV_{\ell}^i}\).
By employing this prolongation operator, a FE space associated with an enriched mesh \({\pT_\ell^i}\) is constructed by  composition of a sequence of prolongation operators,
\[
  \bpX_\ell^i := \bs{\Pi}_{L-1,i}^L \cdots \bs{\Pi}_{\ell,i}^{\ell+1} \bpV_\ell^i, \quad \forall \ \ell \in \{1,\ldots,L-1\},\ i = \{1,2\}.
\]
We borrow the definition of the FE space on the finest level as $\bpX_L^i:= \bpV_L^i$, while the coarse levels, a hierarchy of nested FE spaces associated with each domain is created using such prolongation operator, $\bs{\Pi}_{\ell, i}^{\ell+1}$.
The nested FE spaces are given as
\[
  \bpX_{0}^i \subset \bpX_{1}^i \subset \cdots \subset \bpX_{\ell-1}^i \subset \bpX_{\ell}^i\subset \bpX_{\ell+1}^i\subset \cdots \subset \bpX_{L-1}^i \subset \bpX_{L}^i, \qquad  \forall i\in\{1,2\}.
\]
Following the previous section, we can construct the prolongation operator for the domain $\Omega$ as a direct sum of the prolongation operators on each domain, i.e., $\bs{\Pi}_{\ell-1}^\ell:=\bs{\Pi}_{\ell-1,1}^\ell \oplus \bs{\Pi}_{\ell-1,2}^\ell$, for all $\ell \in \{1,\ldots,L-1\}$.
Thus, we can create the enriched FE spaces $\bpX_\ell := \bpX_\ell^1 \oplus \bpX_\ell^2$ associated with each level $\ell$.
In addition, we can create a hierarchy of nested FE spaces $\{\bpX_\ell\}_{\ell\in\{0,\ldots,L\}}$ for the domain $\Omega$ by using the prolongation operator $\bs{\Pi}_{\ell-1}^\ell$.

The computation of the prolongation operator $\bs{\Pi}_{\ell-1}^\ell$ can be carried out using the $L^2$-projection or pseudo-$L^2$-projection \cite{dickopf_evaluating_2014, kothari_multigrid_2019}.
In this work, we compute the algebraic representation of the prolongation operator $\matT_{\ell-1}^\ell$, with the pseudo-$L^2$-projection.

\subsection{Orthogonal Transformation}
In this section, we introduce the orthogonal transformation for the contact problem \eqref{eq:contact_problem_bar}.
This transformation is necessary to decouple the linear constraints, which in turn allows us to utilize the  modified PGS method.
In order to decouple the constraints, we perform a QR decomposition of the constraint matrix $\overbar{\matB}^\T$
\[
  \overbar{\matB}^\T =  \matQ\matR\ \text{  and   }\  \overbar{\matB} = \matR^\T \matQ^\T,
\]
where $\matQ \in \bR^{nd \times nd}$ is an orthonormal matrix. Thus, we have \( \matQ\matQ^\T = \matQ^\T\matQ = \matI \), where \(\matI\in \bR^{nd \times nd}\), represents the identity matrix.
The decomposition of the matrix \(\matR \in  \bR^{nd \times m}\) is given by $\matR = [\matR_1 \ \matO_1]^\T$, where $\matR_1 \in \bR^{m \times m}$ is an upper triangular matrix and $\matO_1 \in \bR^{(nd-m)\times m}$ is a matrix with all zero entries.
The matrix $\matQ$ simply provides a change of basis, and on this new basis system, the representation of the constraint is modified.
It is clear from the structure of the new constraint matrix $\matR_1$ that in the modified basis system the constraints are sequentially dependent on the previous linear constraint.

We use the matrix $\matQ$ and project the problem on a different basis system.
The matrix $\matQ$ is used to define the variables in the new basis system, given as $\widehat{\vecx} = \matQ^\T \overbar{\vecx}$ and $\overbar{\vecx} = \matQ\widehat{\vecx}$.
Moreover, we can observe that $\matQ^\T \overbar{\matB}^\T = \matR $ and $\overbar{\matB}\matQ = \matR^\T$.
By incorporating the transformed matrices and the vectors, we can reformulate the constrained minimization problem~\eqref{eq:contact_problem_bar} as follows:
\begin{equation}
  \begin{aligned}
     & \min_{\widehat{\vecx}\in \bR^{nd}} J(\widehat{\vecx}) = \half  \widehat{\vecx}^\T\widehat{\matA} \widehat{\vecx} - \widehat{\vecx}^\T \widehat{\vecb} \\
     & \text{subject to }\quad \matR^\T \widehat{\vecx}  \leqslant \vecg,
  \end{aligned}
\end{equation}
where $\widehat{\matA} = \matQ^\T\overbar{\matA}\matQ$ and $\widehat{\vecb} = \matQ^\T\overbar{\vecb}$.
The constraints of the above optimization problem can be written algebraically as,
\begin{equation}
  \begin{pmatrix}
    \matR_{11} & 0          & 0          & \cdots & 0          & 0      & \cdots & 0      \\
    \matR_{11} & \matR_{22} & 0          & \cdots & 0          & 0      & \cdots & 0      \\
    \matR_{13} & \matR_{23} & \matR_{33} & \cdots & 0          & 0      & \cdots & 0      \\
    \vdots     & \vdots     &            & \ddots & \vdots     & \vdots & \ddots & \vdots \\
    \matR_{1m} & \matR_{2m} & \matR_{3m} & \cdots & \matR_{mm} & 0      & \cdots & 0      \\
  \end{pmatrix}
  \begin{pmatrix}
    \widehat{\vecx}_1     \\
    \widehat{\vecx}_2     \\
    \widehat{\vecx}_3     \\
    \vdots                \\
    \widehat{\vecx}_m     \\
    \widehat{\vecx}_{m+1} \\
    \vdots                \\
    \widehat{\vecx}_{nd}  \\
  \end{pmatrix}\leqslant
  \begin{pmatrix}
    \vecg_1 \\
    \vecg_2 \\
    \vecg_3 \\
    \vdots      \\
    \vecg_m
  \end{pmatrix}.
  \label{eq:eqq1}
\end{equation}
As $\matQ$ is an orthonormal matrix, the spectral properties of the $\widehat{\matA}$ and $\matA$ are equivalent.
But, the sparsity pattern of the original matrix $\overbar{\matA}$ and its rotated variant $\widehat{\matA}$ are quite different.
In practice, the matrix $\widehat{\matA}$ is denser than the original matrix, which in turn increases the computational cost of the matrix-vector products in the algorithm.
The new constraint matrix $\matR^\T$ has a lower triangular structure, which can be handled easily by forward substitution.
It is important to note, this type of constraint can be handled easily by the PGS method, due to its inherent sequential nature.

Now, we define a constrained subspace or a feasible set as
\[
  \widehat{\mathbb{K}} = \{\widehat{\vecx} \in \bR^{nd} : \matR^\T\widehat{\vecx} \leqslant \vecg \}.
\]
We pose our problem as an energy minimization problem in the following algebraic formulation:
\begin{equation}
  \text{find }\widehat{\vecx} \in \widehat{\mathbb{K}} \text{ such that } \qquad J(\widehat{\vecx}) \leqslant J(\widehat{\vecy})  \qquad \forall \widehat{\vecy} \in \widehat{\mathbb{K}}.
  \label{eq:min_energy_Q}
\end{equation}

\subsection{Modified Projected Gauss-Seidel Method}
Here, we introduce a modified PGS method for solving the problem~\eqref{eq:min_energy_Q}.
The Gauss-Seidel method is known to minimize the energy functional $J(\cdot)$ in each local iteration step.
The energy minimization takes place in the direction of the nodal basis functions that span the FE space.
The Gauss-Seidel method can be written as a subspace correction method, where the subspace decomposition is achieved by a direct splitting of the underlying FE space into one-dimensional subspaces spanned by the nodal basis functions.
The PGS method is used widely to solve various forms of obstacle problems, and it is known to be globally convergent~\cite{kornhuber_monotone_1994,kornhuber_adaptive_2001}.
We remark that decoupling of the constraints with respect to the nodal basis function is essential for the global convergence of the PGS method~\cite{kornhuber_monotone_1994,rolandglowinski1984-02-14}.
The original linear contact condition, $\overbar{\matB}\overbar{\vecx} \leqslant \vecg$, does not satisfy this property, as the constraints are represented by the linear combination of basis functions.
The QR decomposition allows us to decouple the constraints by expressing them in new basis as \(\matR^\T \widehat{\vecx} \leqslant \vecg\).
In order to discuss this method in generic way, we introduce abstract upper bound $\bs{ub}\in \bR^m$ and lower bound $\bs{lb} \in \bR^m$.
In the context of the contact problem, the lower bound and upper bound are defined as $\bs{lb} = \{-\infty\}$ and $\bs{ub} = \vecg$, respectively.
In addition, we define the set of all active set as a set of all the DoFs where the constraints are binding, thus as
\[
  \mA :=  \{p : (\matR^\T \widehat{\vecx})_p = \vecg_p \}.
\]
The matrix \(\matR^\T\) is a lower triangular matrix, which allows us to write the constraints as a linear combination of the current nodal basis function and previously constrained basis.
This key idea allows us to use the PGS method to solve the problem \eqref{eq:min_energy_Q}.

The iterative process is given as follows.
For a given $k$-th iterate $\widehat{\vecx}^{(k)} \in \widehat{\mathbb{K}}$, we compute a sequence of local intermediate iterates, $\vecz^{(0)},\vecz^{(1)},\ldots,\vecz^{(nd)}$.
We begin with the first local iterate $\vecz^{(0)} := \widehat{\vecx}^{(k)}$, and the next local iterates are given by $\vecz^{(i)} = \vecz^{(i-1)} + \vecc^{(i)}$, for $i = 1,\ldots, nd$.
Once all local intermediate iterates are computed, the new global iterate of the Gauss-Seidel step is given by $\widehat{\vecx}^{(k+1)} := \vecz^{(n)}$.
The corrections $\vecc^{(i)}$ are obtained as the unique solution of the following local subproblems, given as,
\[
  \text{find } \vecc^{(i)} \in \mathbb{D}^{(i)} \text{ such that } \quad J( \vecz^{(i-1)} + \vecc^{(i)} ) \leqslant J(\vecz^{(i-1)} + \vecy ) \qquad \forall \vecy \in\mathbb{D}^{(i)},
\]
with closed, convex set $ \mathbb{D}^{(i)} $, defined for abstract upper bound $\bs{ub}$ and lower bound $\bs{lb}$ as
\begin{equation}
  \label{eq:convex_set_GS}
  \mathbb{D}^{(i)} = \{ \vecc^{(i)} \in \bR^n : \bs{lb} - \matR^\T\vecz^{(i-1)}\leqslant \matR^\T \vecc^{(i)} \leqslant \bs{ub} - \matR^\T \vecz^{(i-1)} \}.
\end{equation}
Each intermediate step ensures that the iterate does not violate the constraints.
If the current iterate violates the constraints, it is projected to the admissible space.
The PGS method for a generic linear inequality constrained minimization problem is summarized in Algorithm~\ref{algo:projectedGS}.

\begin{algorithm*}[t]\small
  \caption{Modified Projected Gauss-Seidel method}
  \label{algo:projectedGS}
  \SetKwData{Left}{left}\SetKwData{This}{this}\SetKwData{Up}{up}
  \SetKwComment{Comment}{$\triangleright$\ }{}
  \SetKwFunction{Union}{Union}\SetKwFunction{FindCompress}{FindCompress}
  \SetKwInOut{Input}{Input}\SetKwInOut{Output}{Output}
  \Input{$\widehat{\matA},\widehat{\vecb},\matR,\widehat{\vecx}^{(0)},\bs{lb},\bs{ub},\nu_\ast$}
  \Output{$\widehat{\vecx}^{(\nu_\ast)}, \mA$}
  \BlankLine
  Function: Projected GS($\widehat{\matA},\widehat{\vecb},\matR,\widehat{\vecx}^{(0)},\bs{lb},\bs{ub},\nu_\ast$)\\
  \For{ $k = 1,2,\ldots,\nu_\ast$}
  {
    $\mA \mapsfrom \emptyset $  \Comment*[r]{initialize empty active set}
    \For{$i = 1,2,\ldots,n$}
    {
      $\widehat{\vecx}^{(k)}_i = \dfrac{1}{\widehat{\matA}}_{ii} (\widehat{\vecb}_i - \sum_{j<i} \widehat{\matA}_{ij} \widehat{\vecx}^{(k)}_j - \sum_{j>i} \widehat{\matA}_{ij} \vecx^{(k-1)}_j)$ \Comment*[]{update the iterate}
      \If{$i \leqslant m$}
      {
        $ {lb}_t  = \dfrac{1}{\matR_{ii}}\big( \bs{lb}_i - \sum_{j=1}^{i-1} \matR_{ji}\widehat{\vecx}^{(k)}_j \big)$;\ ${ub}_t  = \dfrac{1}{\matR_{ii}}\big( \bs{ub}_i - \sum_{j=1}^{i-1} \matR_{ji}\widehat{\vecx}^{(k)}_j \big) $ \Comment*[r]{updated local bounds}
\If{${lb}_t <\widehat{\vecx}^{(k)}_i  < {ub}_t$}
        {
          $\widehat{\vecx}^{(k)}_i =  \max({lb}_t, \min(\widehat{\vecx}^{(k)}_i, {ub}_t)) $ \Comment*[r]{project onto feasible set}
          $\mA \mapsfrom \mA \cup \{i\} $ \Comment*[r]{add current index to the active set}
        }
      }
    }
  }
\end{algorithm*}
Thus, we have a globally convergent PGS method that can be used to solve the problem \eqref{eq:min_energy_Q}.
But the convergence rate of the Gauss-Seidel method is known to deteriorate as the size of the problem increases.
Hence, we employ the modified PGS method as a smoother in our multigrid method. \begin{remark}
  In the Algorithm~\ref{algo:projectedGS}, we have assumed that the values of diagonal entries of the matrix $\matR_1$ (where $\matR = [\matR_1 \ \matO_1]^\T$) are positive, which does not hold in general.
  It is necessary to pay attention to the sign of diagonal entries of the matrix $\matR_1$, as the sign of the diagonal values may change the inequality bounds. \end{remark}
\subsection{Multigrid Method}
In this section, we summarize the generalized multigrid method, which includes all the components introduced in the previous section.
In particular, we have a sequence of non-nested finite element spaces $\{\bpV_\ell\}_{\ell=0,\ldots, L}$ associated with the hierarchy of meshes $\{\pT_\ell\}_{\ell=0,\ldots, L}$.
Following Section~\ref{sec:transfer}, we have the transfer operators $\{\bs{\Pi}_{\ell-1}^\ell\}_{\ell=1,\ldots,L}$ which are computed using the pseudo-$L^2$-projections.
By means of these transfer operators, we create a hierarchy of nested finite element spaces $\{\bpX_\ell\}_{\ell=0, \ldots, L}$ from the hierarchy of background meshes.
The prolongation matrices associated with the transfer operators are given as $\{\matT_{\ell-1}^\ell\}_{\ell=1,\ldots,L} $.
\begin{algorithm*}[t]\small
  \caption{Generalized Multigrid algorithm}
  \SetKwComment{Comment}{$\triangleright$\ }{}
  \label{algo:mmg}
  \SetKwInOut{Input}{Input}\SetKwInOut{Output}{Output}
  \Input{$\matA_L, \vecb_L, L, \nu_1, \nu_2, (\matT_{\ell-1}^\ell)_{\ell=1,\dots,L}, \matB,\matO, \bs{lb},\bs{ub} ,\gamma$}
  \Output{$\vecx_L \mapsfrom \matO \overbar{\vecx}_L$ ,$\overbar{\vecx}_L \mapsfrom \matQ \widehat{\vecx}_L$}
  \BlankLine
  Function: GMG($\matA_{L}, \vecb_L, L, \nu_1, \nu_2, (\matT_{\ell-1}^\ell)_{\ell=1,\dots,L}, \matB, \matO, \bs{lb},\bs{ub}, \gamma$)\\
  $\widehat{\vecx}_{L}       \mapsfrom \bs{0}$\Comment*[r]{initialize solution}
  $\overbar{\matT}^L_{L-1}   \mapsfrom \matO \matT_{L-1}^L;\ \overbar{\matA}_L \mapsfrom \matO \matA_L \matO;\ \overbar{\vecb}_L \mapsfrom \matO \vecb_L; \overbar{\matB} \mapsfrom \matB\matO $\Comment*[r]{local basis transformation}
  $\matQ,\matR \mapsfrom $ QR Transformation($\overbar{\matB}^\T$)\Comment*[r]{QR decomposition}
  $\widehat{\matT}^L_{L-1}   \mapsfrom \matQ^\T\overbar{\matT}_{L-1}^L;\ \widehat{\matA}_L \mapsfrom \matQ^\T \overbar{\matA}_L \matQ;\ \widehat{\vecb}_L \mapsfrom \matQ^\T \overbar{\vecb}_L $\Comment*[r]{orthogonal rotation}
  \BlankLine
  \While{not converged}{
  $\widehat{\vecx}_L,\mA_L \mapsfrom \widehat{\vecx}_L + $ Projected GS($ \widehat{\matA}_L , \widehat{\vecb}_L, \matR^\T,\widehat{\vecx}_L,\bs{lb},\bs{ub}, \nu_1$)\Comment*[r]{$\nu_1$ pre-smoothing steps}
  $\widehat{\vecr}_{L} \mapsfrom \widehat{\vecb}_L - \widehat{\matA}_L \widehat{\vecx}_L$\Comment*[r]{residual}
  $\widehat{\vecr}_{trc} \mapsfrom \text{trc}(\widehat{\vecr}_L,\mA_L);\ \widehat{\matA}_{trc} \mapsfrom \text{trc}(\widehat{\matA}_L,\mA_L) $\Comment*[r]{truncation}
  $\vecr_{L-1} \mapsfrom (\widehat{\matT}^L_{L-1})^\T \widehat{\vecr}_{trc}$ \Comment*[r]{restriction}
  $\matA_{L-1} \mapsfrom (\widehat{\matT}^L_{L-1})^\T \widehat{\matA}_{trc} \widehat{\matT}^L_{L-1} $\Comment*[r]{Galerkin projection}
  $\vecc_{L-1} \mapsfrom \bs{0} $ \Comment*[r]{initialize coarse level correction}
  \For{$i=1,\ldots,\gamma$}
  {
    $\vecc_{L-1} \mapsfrom \vecc_{L-1} +$ MG($\matA_{L-1},\vecr_{L-1},L-1,\nu_1,\nu_2,\matT_{L-2}^{L-1},\gamma$)\Comment*[r]{coarse level cycle}
  }
  $\widehat{\vecc}_L \mapsfrom \widehat{\matT}_{L-1}^{L} \vecc_{L-1} $\Comment*[r]{prolongation}
  $\widehat{\vecc}_{\text{trc}} \mapsfrom \text{trc}(\widehat{\vecc}_{L},\mA_L)$\Comment*[r]{truncation}
  $\widehat{\vecx}_L \mapsfrom \widehat{\vecx}_L + \widehat{\vecc}_{\text{trc}} $\Comment*[r]{update iterate}
  $\widehat{\vecx}_L,\mA_L \mapsfrom \widehat{\vecx}_L + $ Projected GS($ \widehat{\matA}_L , \widehat{\vecb}_L, \matR^\T, \widehat{\vecx}_L, \bs{lb},\bs{ub},\nu_2$)\Comment*[r]{$\nu_2$ post-smoothing steps}
  }
\end{algorithm*}

The orthogonal transformation of the matrix $\overbar{\matB}^\T$ plays a vital role in our multigrid method.
We recall, on the finest level the definition of the FE space is kept the same, as $\bpX_L=\bpV_L$, and hence also the nodal basis functions defined on these FE spaces are given as $\zeta_L^p=\phi_L^p$, for all $p \in \mN_L$.
These nodal basis functions are modified or rotated after the orthogonal transformation, which can be written as
\[
  \widehat{\zeta}_L^q \vece_i(q) := \sum_{p\in \mN_L} \matQ_{pq} \zeta_L^p \vece_i(p) \qquad \forall q\in \mN_L.
\]
The transfer operators are computed using the nodal basis functions that span the FE space on a coarse level and a fine level.
With the modified nodal basis functions on the finest level, it becomes essential to compute the transfer operator associated with the finest level such that the vector and the matrix quantities are projected on the FE space spanned by a modified basis system.
Thus, the prolongation matrix $\matT_{L-1}^L$ is also modified in two stages.
The first stage is necessary because of the local basis transformation, which is carried out locally to modify the basis system by means of the Householder rotation matrix $\matO$ such that contact conditions are only applicable in the normal direction.
The second transformation is carried out using the orthogonal transformation matrix $\matQ$, which is computed to decouple the linear contact constraints.
The updated transfer operator is defined as
\[
  \widehat{\matT}_{L-1}^L  = \matQ^\T \overbar{ \matT}_{L-1}^L = \matQ^\T \matO \matT_{L-1}^L.
\]
Now, the nodal basis function associated with the FE space $\bpX_{L-1}$ are given as
\[
  \widehat{\zeta}_{L-1}^q \vece_i(q):= \sum_{p\in \mN_L} (\widehat{\matT}_{L-1}^L )_{pq} \zeta_L^p \vece_i(p) = \sum_{p\in \mN_L} (\matQ^\T \overbar{\matT}_{L-1}^L )_{pq} \zeta_L^p \vece_i(p), \qquad \forall q\in \mN_{L-1}.
\]
This modification of the transfer operator is only required on the finest level, while all other transfer operators on the coarser levels $\{\matT_{\ell-1}^\ell\}_{\ell=0,\ldots, L-1}$ remain the same.

The modified PGS method is employed as a smoother in the generalized multigrid method only on the finest level.
It minimizes the energy functional in each local iteration in each smoothing step.
At the end of the smoothing iterations, we obtain a list of an active set where the constraints are binding.
The most crucial feature of this multigrid method is that the coarse level corrections do not violate the fine level constraints.
As a consequence, we solve the constrained optimization problem only on the finest level, while on the coarse levels, we solve the unconstrained linear problem.
This is also very convenient, as in this algorithm the representation of the contact constraints is only required on the finest level.

To ensure that the coarse level corrections do not violate the constraint on the finest level, we modify the restriction of the residual and the stiffness matrix, and the prolongation of the coarse level correction.
Following the discussion in Section~\ref{sec:transfer}, we know the nodal basis functions associated with the coarse level FE space are computed as a linear combination of the nodal basis function defined on the FE space on the finest level.
If the value of a nodal basis function on the finest level is set to zero, the basis function constructed on the coarse levels is represented by truncated basis functions.
For all DoFs that are in the active set, we set the corresponding entries of the residual or the prolongated correction to zero.
While for the stiffness matrix, we set the rows and columns associated with the active set to be zero.
This is equivalent to removing the nodal basis function associated with all DoFs in the active set.

As we are employing transfer operators constructed by the pseudo-$L^2$-projection, this multigrid method including the truncation process can be carried out algebraically.
In comparison with the standard multigrid method used for solving the linear system, this algorithm is computationally more expensive.
This can be attributed to the cost of computing the orthogonal transformation of the matrix $\overbar{\matB}^\T$ and then projecting the problem onto a new basis system.
Even though the generalized multigrid method is computationally more expensive, it has optimal convergence properties.
Additionally, it is significantly cheaper in comparison with the other iterative methods, e.g., interior-point method or semi-smooth Newton method.
If we are solving an optimization problem with inequality constraints, the active set changes in a few initial multigrid iterations.
However, once the active set of the solution is identified, the algorithm converges linearly.

In Algorithm~\ref{algo:mmg}, we can see the detailed generalized multigrid algorithm with the modified PGS method as a smoother.
On the coarse levels, we employ the standard MG methods as described in the Algorithm~\ref{algo:mmg_coarse}, with any regular smoothers.
Here, we note that Algorithm~\ref{algo:mmg} is given in an abstract setting for inequality constraints with upper bounds and lower bounds, assuming that the active set may change in each multigrid iteration.

 \section{Numerical Results}\label{sec:results}
In this section, we evaluate the performance of the proposed generalized multigrid method for Signorini’s problem and the two-body contact problem.
The non-penetration condition at the contact interface is discretized using the vital vertex algorithm.

We utilize Givens rotation to perform the orthogonal transformation of the matrix $\matB^\T$, as this method produces a sparser matrix $\matQ$ than the other methods.
The performance analysis of the multigrid method is carried out in the Utopia library~\cite{utopiagit}.
For these experiments, we choose correction in energy norm as a termination criterion, given as
\begin{equation}
  {\| \vecx^{(k+1)} - \vecx^{(k)} \|_A} < 10^{-10}.
  \label{eq:termination_cont}
\end{equation}
Also, we define the convergence rate of an iterative solution scheme as
\[
  \rho^{k+1} := \dfrac{\|\vecx^{(k+1)} - \vecx^{(k)} \|_A }{\|\vecx^{(k)} -\vecx^{(k-1)} \|_A}.
\]
We denote the asymptotic convergence rate as $\rho^\ast$, where the iterate $\vecx^{(k+1)}$ satisfies the termination criterion.

\subsection{Signorini's Problem}
In this section, we describe the problem setup for Signorini's problem for two different types of rigid obstacles.

\subsubsection{Problem Description}
All experiments in this section are carried out on a structured background mesh with the quadrilateral elements.
On the coarsest level, the background mesh $\widetilde{\pT}_0$ is given on a rectangle of dimension $[-1.09,1.09] \times [0, 1.09]$, with $100$ elements in $X$-direction and $50$ in $Y$-direction, which we denoted as mesh on level $L0$.
By uniformly refining the mesh $\widetilde{\pT}_0$, we obtain a hierarchy of meshes $\{\widetilde{\pT}_\ell\}_{\ell=1}^5$ associated with the levels $L1,L2,\ldots,L5$.
The Dirichlet boundary condition is defined as $\vecu = (0,0)$ on $x=[-1.09, 1.09]$ and $y=0$.
The body force for this example is considered to be zero.
In these experiments, the material parameters are chosen as Young's modulus $E=10\,\text{MPa}$ and Poisson's ratio $\nu=0.3$.
We can compute Lam\'e parameters $\lambda$ and $\mu$ using the following relation:
\[
  \lambda = \frac{E \nu}{(1+\nu)(1-2\nu)} \quad \text{and} \quad \mu = \frac{E}{2(1+\nu)}.
\]

\begin{figure*}[t]
  \begin{subfigure}{0.49\textwidth}
    \includegraphics[width=.9\linewidth]{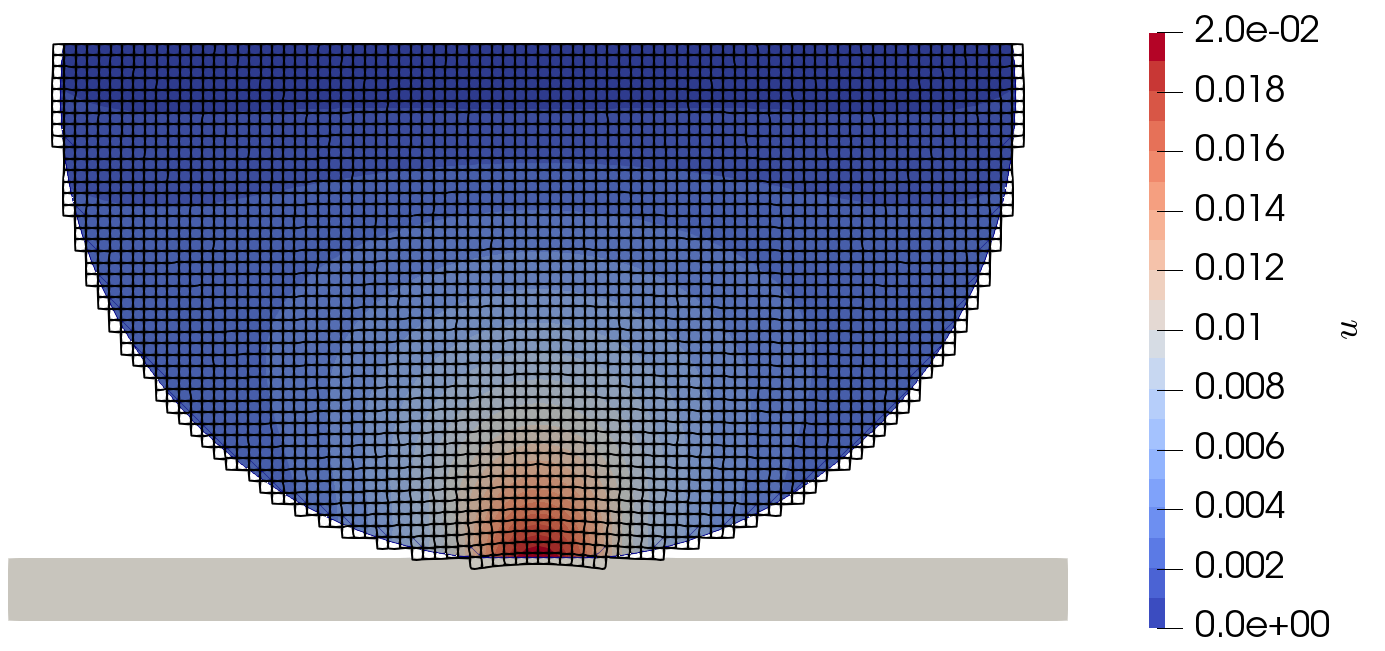}
    \caption{Setup of the Signorini's problem for \nameref{para:example1sc}.}\label{fig:singonri_1}
  \end{subfigure}
  \hfill
  \begin{subfigure}{0.49\textwidth}
    \centering
    \includegraphics[width=.9\linewidth]{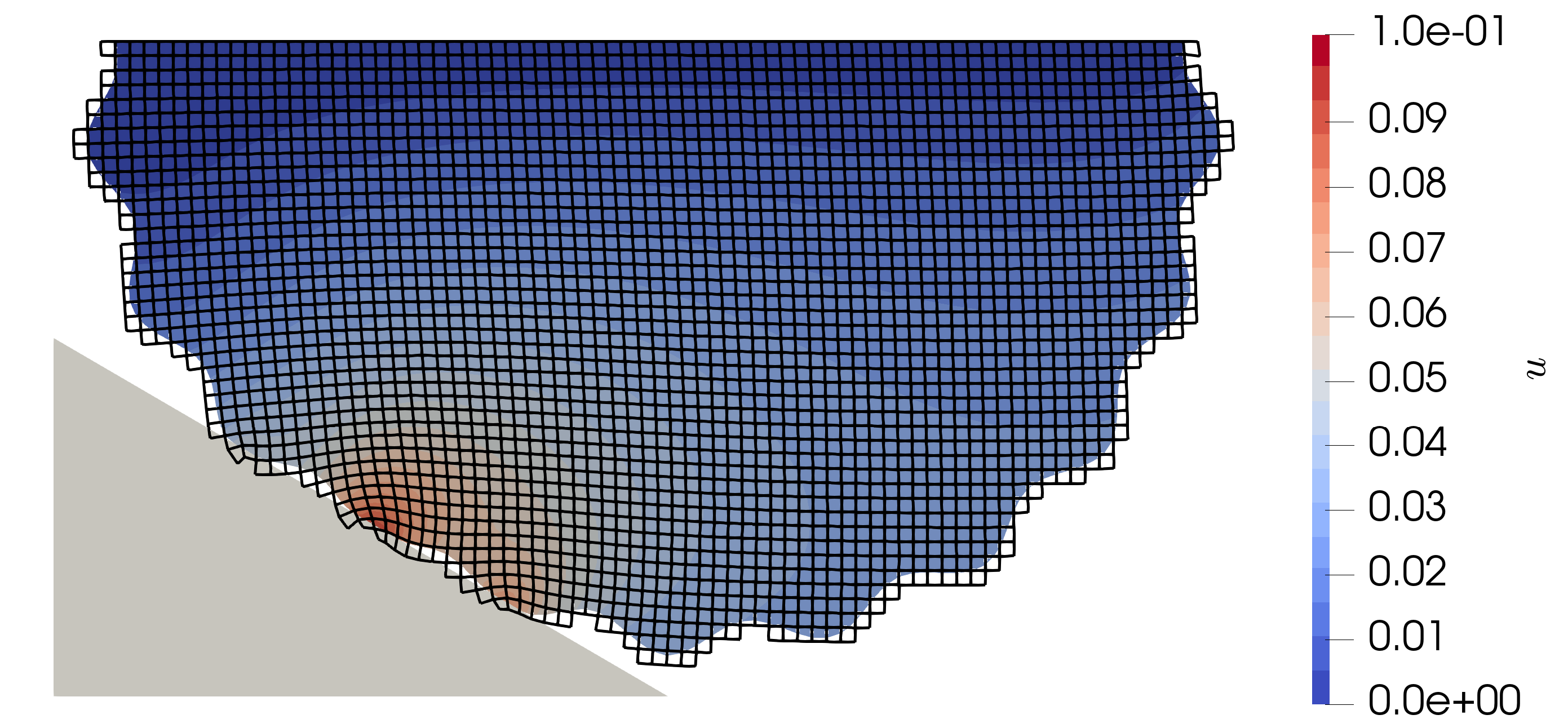}
    \caption{Setup of the Signorini's problem for \nameref{para:example2sc}.}\label{fig:singonri_2}
  \end{subfigure}
  \caption{Setup of for experiments, the object in the gray scale is the rigid obstacle. We can see the active background mesh and the displacement field.}
\end{figure*}

\paragraph{Example 1-SC}\label{para:example1sc}
For this example, we consider a rigid foundation, defined by a line $y = 0.12$.
The body $\Omega$ is pressed against the rigid foundation, and the maximum magnitude of the displacement on the body is given as $u= 0.02$.
In this experiment, we consider a semicircular domain, where the contact boundary of the domain is defined by a zero level set of a function $\Lambda_{s_1}(\matX):= r_{s_1}^2 - \|\matX - \vecc_1\|_2^2$ with radius $r_{s_1}=0.9$, and $\vecc_1$ denotes the center of the circle $(0,1)$.
The domain $\Omega$ is defined by the region where the value of the level set is positive, $\Lambda_{s_1} > 0$.
The setup of this example is depicted in Figure~\ref{fig:singonri_1}, where we see the resultant magnitude of the displacement field due to the contact with a rigid foundation.

\paragraph{Example 2-SC}\label{para:example2sc}
This example considers a non-symmetric obstacle and possible multiple contact regions.
The body $\Omega$ is considered as a semicircular domain with wavy boundaries, given as $\Lambda_{s_2}(\matX) := r^2 - 0.5(\tilde{x}^2+\tilde{y}^2)(5 + 0.3\sin(\frac{\pi}{36} + 21 \arctan(\frac{\tilde{y}}{\tilde{x}}) ))$
where, $\tilde{x}=\matX_x-c_y$ and $\tilde{y}=\matX_y-c_y$.
Here, the radius is defined as $r = 1.31111$ and the center is given as $\vecc = (10^{-5},1+10^{-5})$.
The domain $\Omega$ is defined as the region where the level set function $\Lambda_{s_2}$ has positive values.
The rigid foundation is defined as a line which passes through the points $(-1,0.7)$ and $(0.2,0)$.
In Figure~\ref{fig:singonri_2}, we can observe the setup and magnitude of the displacement field due to contact with the rigid foundation.

\subsubsection{Convergence Study of the Discretization Method}
\begin{figure*}[t] \begin{subfigure}{0.32\textwidth}
    \includegraphics[trim={6.5cm 7cm 5.2cm 7cm},clip,width=\linewidth]{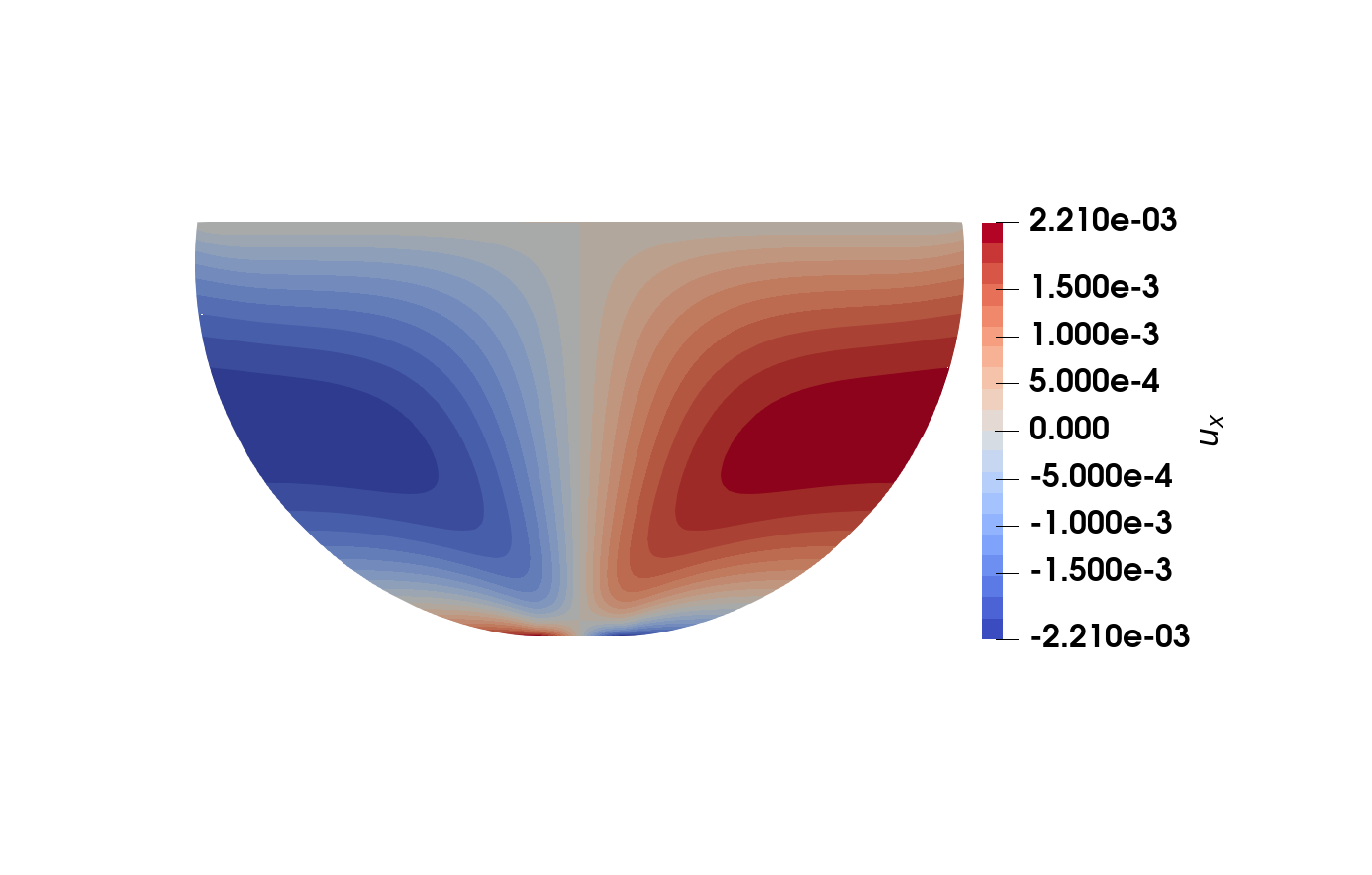}
    \caption{Displacement in X direction $\vecu_x$.}\end{subfigure}
  \hfill
  \begin{subfigure}{0.32\textwidth}
    \includegraphics[trim={6.5cm 7cm 5.2cm 7cm},clip,width=\linewidth]{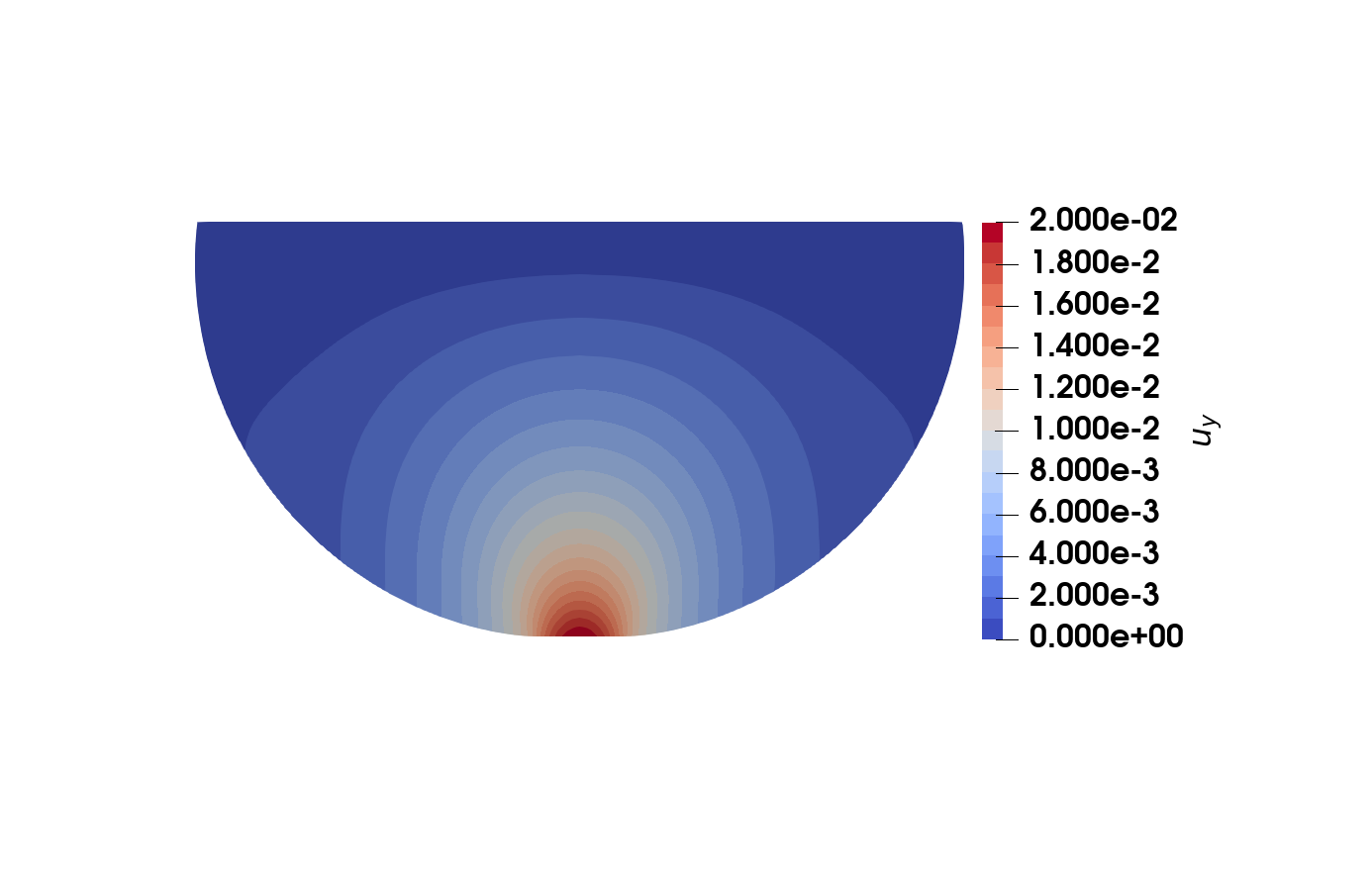}
    \caption{Displacement in Y direction $\vecu_y$.}\end{subfigure}
  \hfill
  \begin{subfigure}{0.32\textwidth}
    \includegraphics[trim={6.5cm 7cm 5.2cm 7cm},clip,width=\linewidth]{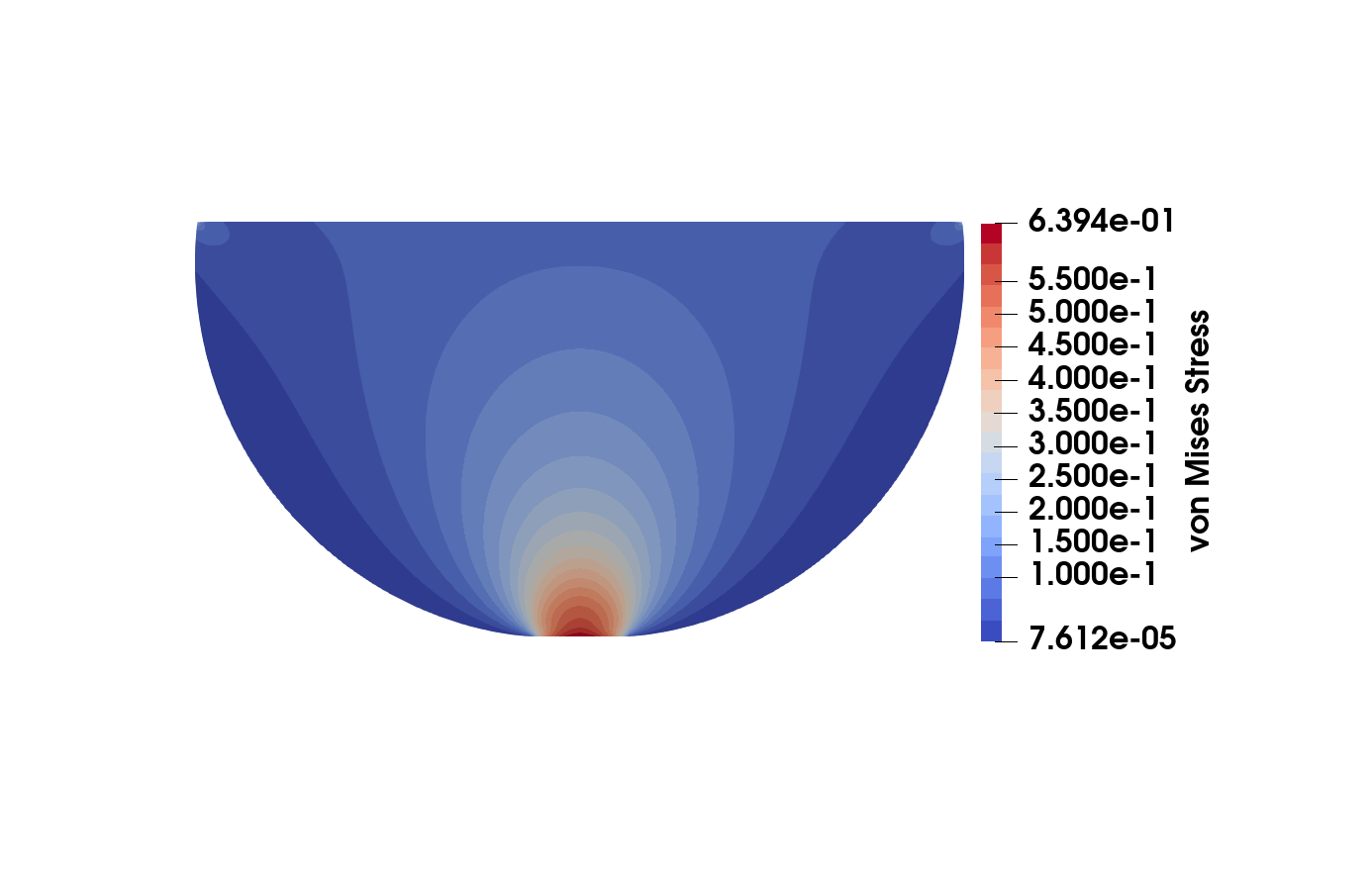}
    \caption{von Mises Stress.} \end{subfigure}
  \par\bigskip \begin{subfigure}{0.32\textwidth}
    \includegraphics[trim={6.5cm 7cm 5.2cm 7cm},clip,width=\linewidth]{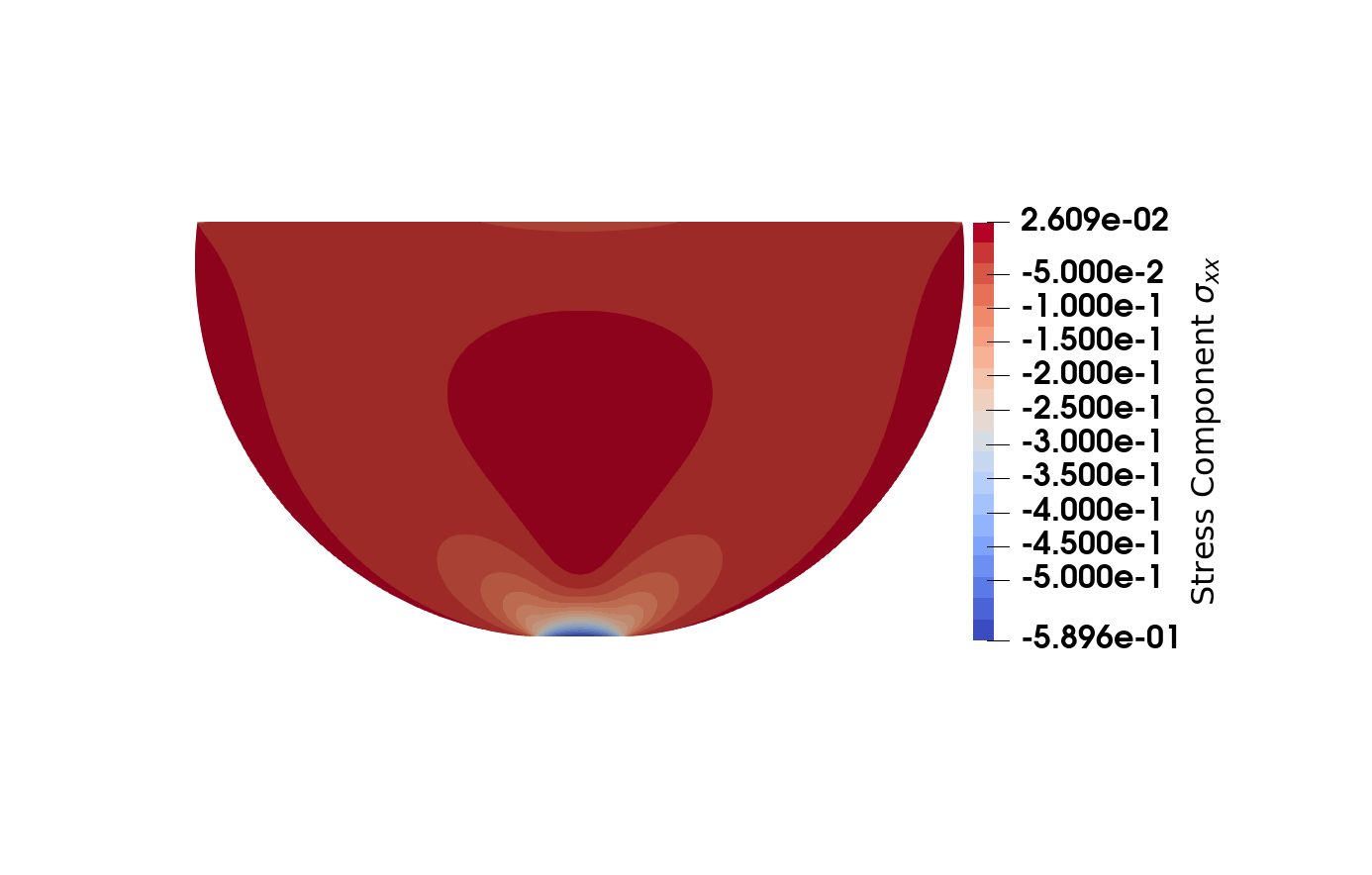}
    \caption{Stress component $\bs{\sigma}_{xx}$.} \end{subfigure}
  \hfill
  \begin{subfigure}{0.32\textwidth}
    \includegraphics[trim={6.5cm 7cm 5.2cm 7cm},clip,width=\linewidth]{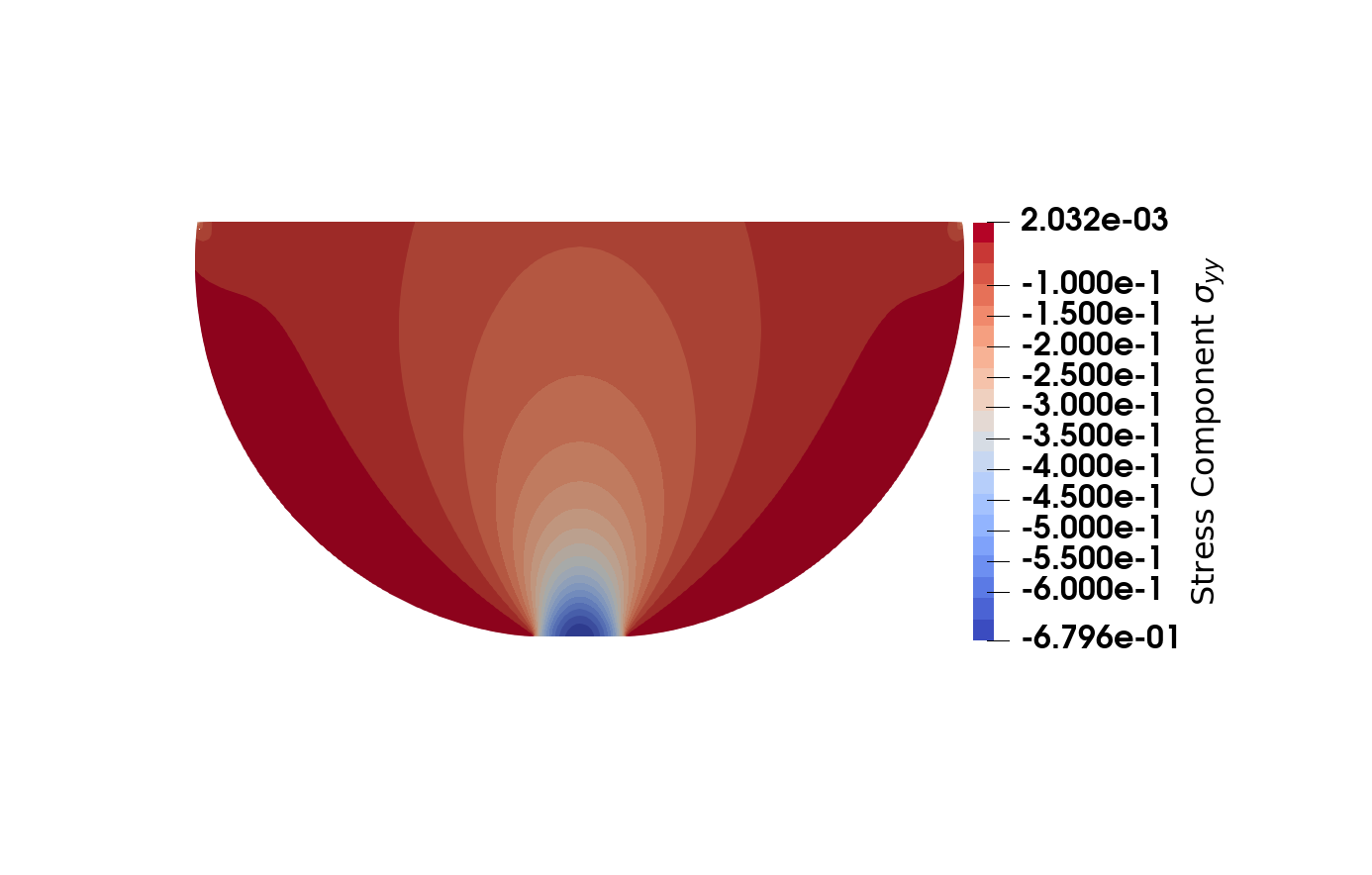}
    \caption{Stress component $\bs{\sigma}_{yy}$.} \end{subfigure}
  \hfill
  \begin{subfigure}{0.32\textwidth}
    \includegraphics[trim={6.5cm 7cm 5.2cm 7cm},clip,width=\linewidth]{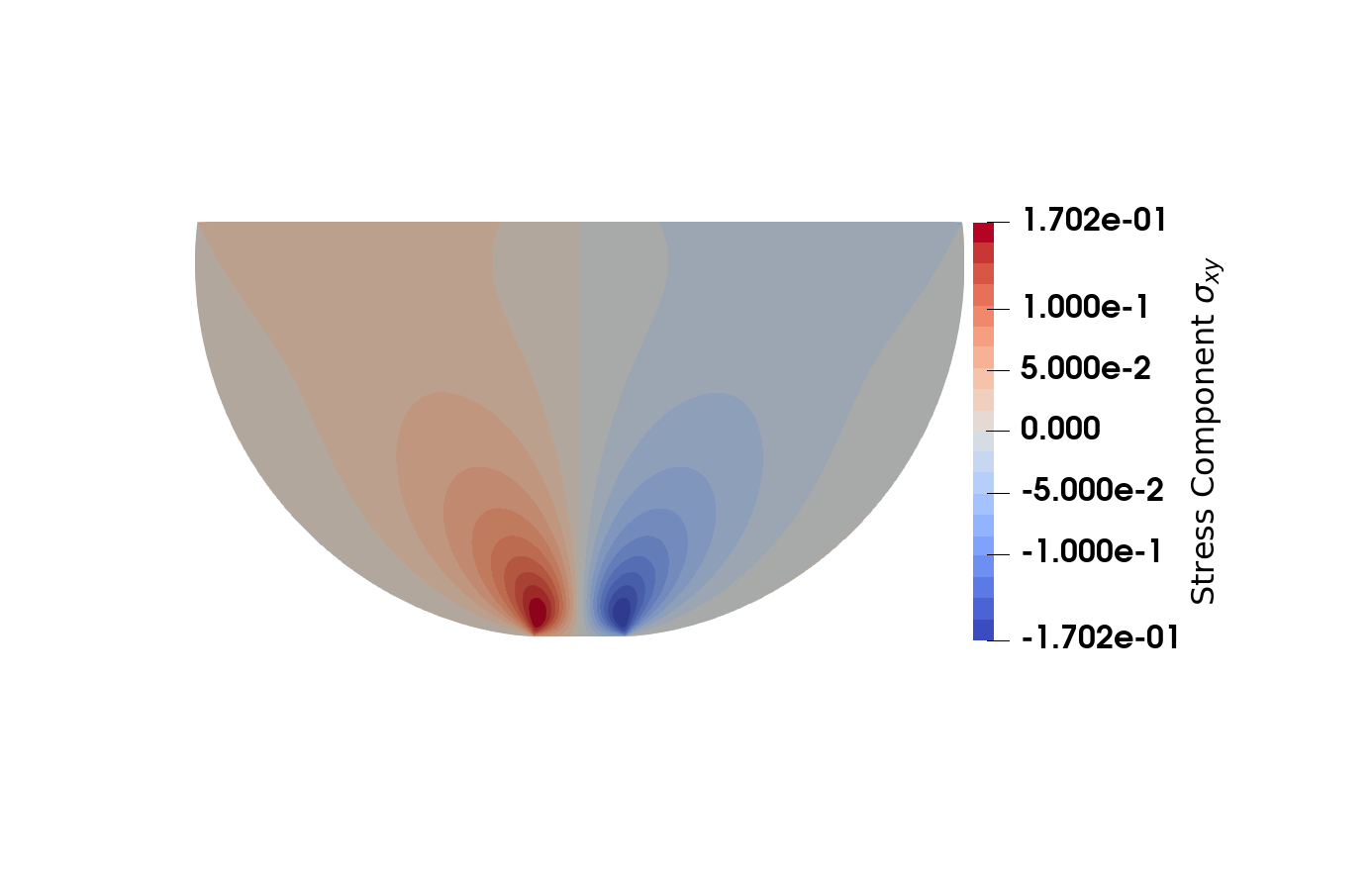}
    \caption{Stress component $\bs{\sigma}_{xy}$.}\end{subfigure}\hspace*{\fill}
  \caption{Resultant displacement field and stress field of citep \nameref{para:example1sc}.}
  \label{fig:1}
\end{figure*}

In this section, we evaluate the performance of the unfitted discretization method and Lagrange multiplier method introduced in Section~\ref{sec:model}.
We consider \nameref{para:example1sc} and \nameref{para:example2sc} as two test cases.
We use the same mesh hierarchy defined on levels $L0, L1, \ldots, L5$.
The solution computed on the mesh on the finest level $L5$ is taken as the reference solution and it is compared against the solutions on different discretization levels from $L0, L1, \ldots, L4$.
Also, we are employing the ghost penalty stabilization term in the bilinear form with the parameter $\epsilon_G = 10^{-2}$.
The resultant components of the displacement field, Cauchy stresses, and von Mises stress for \nameref{para:example1sc} are shown in Figure~\ref{fig:1}.

The discretization error of the displacement field is computed in two different norms, given as the energy norm $\|\cdot\|_{E(\Omega)}$  defined as
\[
  \|e(u_\ell)\|_{E(\Omega)}   := \|u_{\text{ref}} - u_\ell\|_{E(\Omega)} = \Bigg( \sum_{K\in \pT_\ell} \int_K \bs{\sigma}(\vecu_{\text{ref}} - \vecu_\ell) : \bs{\varepsilon}(\vecu_{\text{ref}} - \vecu_\ell) \ d\Omega \Bigg)^{\half}
\]
and the $H_1$-norm given as
\[
  \|e(u_\ell)\|_{H^1(\Omega)} := \|u_{\text{ref}} - u_\ell\|_{H^1(\Omega)} = \Bigg( \sum_{K\in \pT_\ell} \int_K (\vecu_{\text{ref}} - \vecu_\ell)^2 + (\nabla \vecu_{\text{ref}} - \nabla \vecu_\ell)^2 \ d\Omega \Bigg)^{\half}.
\]
Additionally, the discretization error in the normal stresses on the contact boundary is given in a mesh dependent norm on the interface as
\[
  \|e((\sigma_n)_\ell)\|_{H^{-\half}(\Gamma_c),h} := \|(\sigma_n)_{\text{ref}} - (\sigma_n)_\ell\|_{H^{-\half}(\Gamma_c),h} = \Bigg( \sum_{K\in \pT_{\ell,\Gamma_c}} \int_{\Gamma_K} h_K  \big( (\sigma_n)_{\text{ref}} - (\sigma_n)_\ell \big)^2 \ d \Gamma \Bigg)^{\half}.
\]

In Table~\ref{tab:discretization_error}, we show the discretization error for \nameref{para:example1sc} and \nameref{para:example2sc} in the three different norms and also compute the experimental order of convergence (EOC), $\text{EOC}_\ell = \log \big(\frac{\|e(u_{\ell-1})\|}{\|e(u_{\ell})\|}\big) / \log\big(\frac{h_{\ell-1}}{h_\ell}\big)$.
From Table~\ref{tab:discretization_error}, we can see that the EOC in the energy-norm and the $H^1$-norm have optimal convergence rate, as the error reduces by order $\pO(h)$.
While for the normal stresses on the contact interfaces, we expect the convergence rate in the discretization error to be of order $\pO(h^{\frac{3}{2}})$.
We can see from Table~\ref{tab:discretization_error}, \nameref{para:example1sc} demonstrates the optimal EOC for the normal stresses.
For \nameref{para:example2sc} we can see a better convergence rate than the optimal, this behavior can be attributed to the complex geometry, as we refine the mesh the geometry is captured more accurately the rate of convergence of the normal stresses accelerates.

\begin{table*}[t]
  \begin{subtable}{1\textwidth}
    \centering
    \begin{tabular}{|c | c c | c c | c c|}
      \hline
      $h_\ell$             & $\|e(u_\ell)\|_{H^1(\Omega)}$ & EOC$_\ell$ & $\|e(u_\ell)\|_{E(\Omega)}$ & EOC$_\ell$ & $\|e((\sigma_n)_\ell)\|_{H^{-\half}(\Gamma_c),h}$ & EOC$_\ell$ \\ \hline
      \rowcolor{Gray}
      $2.18 \cdot 10^{-2}$ & $1.28662 \cdot 10^{-3}$       & -          & $6.38102 \cdot 10^{-3}$     & -          & $3.38268\cdot 10^{-3}$                            & -          \\
      $1.09 \cdot 10^{-2}$ & $6.51964 \cdot 10^{-4}$       & 0.981      & $3.24611 \cdot 10^{-3}$     & 0.975      & $1.20316\cdot 10^{-3}$                            & 1.491      \\
      \rowcolor{Gray}
      $5.45 \cdot 10^{-3}$ & $3.29916 \cdot 10^{-4}$       & 0.983      & $1.64475 \cdot 10^{-3}$     & 0.981      & $4.98003\cdot 10^{-4}$                            & 1.273      \\
      $2.72 \cdot 10^{-3}$ & $1.64492 \cdot 10^{-4}$       & 1.004      & $8.20526 \cdot 10^{-4}$     & 1.003      & $2.06753\cdot 10^{-4}$                            & 1.268      \\
      \rowcolor{Gray}
      $1.36 \cdot 10^{-3}$ & $7.37150 \cdot 10^{-5}$       & 1.158      & $3.68012 \cdot 10^{-4}$     & 1.157      & $7.20550\cdot 10^{-5}$                            & 1.521      \\ \hline
    \end{tabular}
    \caption{\nameref{para:example1sc}.}
  \end{subtable}

  \bigskip
  \begin{subtable}{1\textwidth}
    \centering
    \begin{tabular}{|c | c c | c c | c c|}
      \hline
      $h_\ell$             & $\|e(u_\ell)\|_{H^1(\Omega)}$ & EOC$_\ell$ & $\|e(u_\ell)\|_{E(\Omega)}$ & EOC$_\ell$ & $\|e((\sigma_n)_\ell)\|_{H^{-\half}(\Gamma_c),h}$ & EOC$_\ell$ \\ \hline
      \rowcolor{Gray}
      $2.18 \cdot 10^{-2}$ & $7.94344 \cdot 10^{-3}$       & -          & $4.11494 \cdot 10^{-2}$     & -          & $3.85899\cdot 10^{-2}$                            & -          \\
      $1.09 \cdot 10^{-2}$ & $4.12698 \cdot 10^{-3}$       & 0.945      & $2.14838 \cdot 10^{-2}$     & 0.938      & $1.34525\cdot 10^{-2}$                            & 1.520      \\
      \rowcolor{Gray}
      $5.45 \cdot 10^{-3}$ & $2.08019 \cdot 10^{-3}$       & 0.988      & $1.09024 \cdot 10^{-2}$     & 0.979      & $4.34291\cdot 10^{-3}$                            & 1.631      \\
      $2.72 \cdot 10^{-3}$ & $1.02215 \cdot 10^{-3}$       & 1.025      & $5.36743 \cdot 10^{-3}$     & 1.022      & $1.46408\cdot 10^{-3}$                            & 1.569      \\
      \rowcolor{Gray}
      $1.36 \cdot 10^{-3}$ & $4.57430 \cdot 10^{-4}$       & 1.160      & $2.40582 \cdot 10^{-3}$     & 1.158      & $3.86884\cdot 10^{-4}$                            & 1.920      \\ \hline
    \end{tabular}
    \caption{\nameref{para:example2sc}.}
  \end{subtable}
  \caption{Discretization error in the displacement field at different levels in $H^1$-norm, energy norm and the normal stresses in mesh dependent norm on contact interface.}
  \label{tab:discretization_error}
\end{table*}

\begin{figure*}[t]
  \begin{subfigure}{0.49\textwidth}
    \centering
    \includegraphics{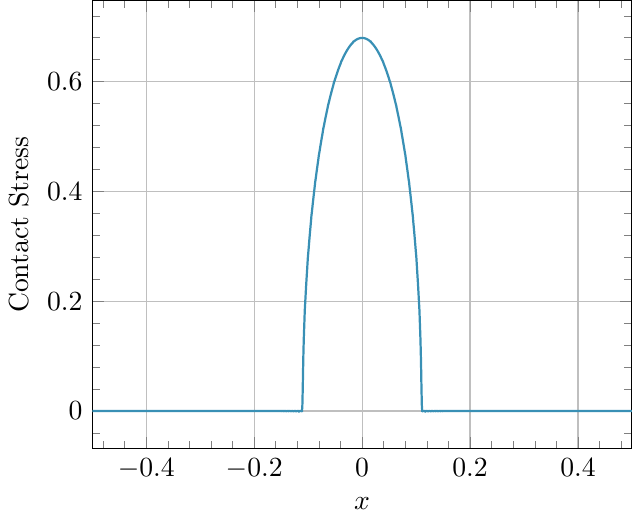}
\caption{\nameref{para:example1sc}.}\end{subfigure}
  \hfill
  \begin{subfigure}{0.49\textwidth}
    \centering
    \includegraphics{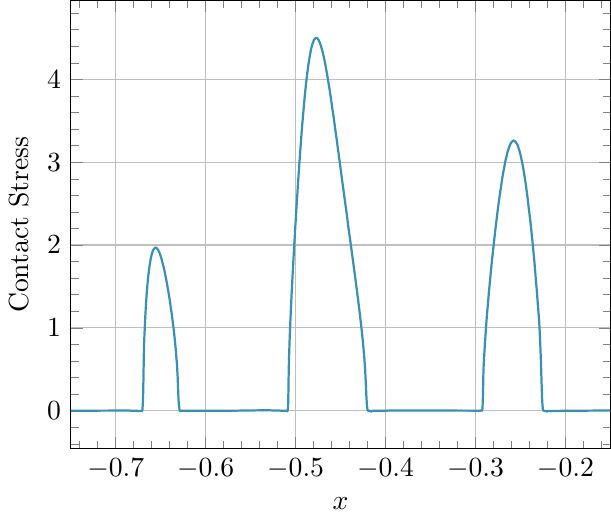}
\caption{\nameref{para:example2sc}.}\end{subfigure}
  \caption{The induced normal contact stress on the contact boundary.}
  \label{fig:contact_stresses}
\end{figure*}

In Figure~\ref{fig:contact_stresses}, we can see the resultant induced normal stresses on the contact boundary for \nameref{para:example1sc} and \nameref{para:example2sc}.
The resultant contact stresses also agree with the Hertzian contact theory.
As for \nameref{para:example1sc}, the contact stress exists only in the contact boundary and it is zero everywhere else, and the contact stress also has a parabolic shape.
For \nameref{para:example2sc}, we can also observe three distinct zones where the contact stresses can be seen.
Additionally, as the rigid foundation is non-symmetric, the resultant normal stresses are also non-symmetric.

\subsubsection{Effect of the Ghost Penalty Term}
In Section~\ref{sec:ghost_penalty}, we discussed the ghost penalty term which is used to overcome the issue of ill-conditioning.
In this section, we evaluate the robustness of this term by comparing the effect of different values of the ghost penalty parameter ($\epsilon_G$).
Additionally, we also evaluate the performance of the proposed multigrid method with respect to various values of the parameter $\epsilon_G$.

\begin{figure*}[t]
  \begin{center}
    \includegraphics{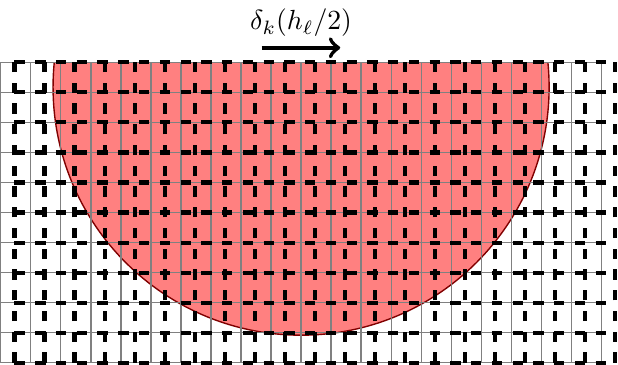}
\end{center}
  \caption{Moving the background mesh, while keeping the domain $\Omega$ fixed. The solid mesh is background mesh configuration for \nameref{para:example1sc}, while the mesh denoted with dashed line is the mesh moved in X-direction with $\delta_k (h_\ell/2)$ distance.}\label{fig:moving_mesh}
\end{figure*}

In this experiment, we consider \nameref{para:example1sc} on the predefined domain with the discussed boundary conditions.
While the domain is kept fixed, we move the background mesh $\widetilde{\pT}_h$ in the $X$-direction.
We generate a set of background meshes $\{\widetilde{\pT}^k_h\}_{k=0}^{10}$.
Here, we consider the problem defined on level $L2$, with $400$ and $200$ elements in $X$-direction and $Y$-direction, respectively.
We recall that the $\widetilde{\pT}_h$ is defined on a rectangle of dimension $[-1.09,1.09]\times[0,1.09]$.
The meshes $\widetilde{\pT}^k_h$ are given as $ \widetilde{\pT}^k_h= \widetilde{\pT}_h + \delta_k (h_{\ell}/2, 0) $, where $\delta_k = 0.1 k$ and $h_\ell = 5.45\cdot10^{-3}$.
A sketch of translated mesh configuration can be seen in Figure~\ref{fig:moving_mesh}.
We note that for this experiment, the multigrid method utilizes $V(5,5)$-cycle with the modified PGS method as a smoother on the finest level, while the symmetric Gauss-Seidel method is used as a smoother on coarser levels.
On the coarsest level, we use a direct solver.
The multigrid hierarchy is equipped with 3 levels, where $L0$ denoted the coarsest level and $L2$ denotes the finest level.

\begin{figure*}[t]
  \centering
  \tikzset{external/export next=false}
\begin{tikzpicture}[spy using outlines={rectangle,magnification=3,connect spies,size=2cm,blue}]
    \begin{groupplot}[
        group style={
group size = 3 by 1,
x descriptions at=edge bottom,
horizontal sep=30pt
          },
        xtick={0,0.2,0.4,0.6,0.8,1},
        minor x tick num=1,
scale=0.64,
        grid=major,
        xlabel={$\delta_k$},
        label style={font=\scriptsize}, tick label style={font=\scriptsize}, legend style={font=\scriptsize}
      ]
\nextgroupplot[align=left, title={{\scriptsize Condition number $\kappa(\matA)$}},ymode=log,ytick={1e5,1e6,1e7,1e8,1e9,1e10,1e11,1e12,1e13,1e14},ymin=1e5,ymax=1e14]
      \addplot [very thick,red!70!black,mark=square]table[x=epsbyh, y=gamma_g_1e-1, col sep=comma]{results/condnum_vs_gp_cut_configuration_variable.csv};
      \label{pgfplots:gamma1em1}
      \addplot [very thick,green!70!black,mark=triangle]table[x=epsbyh, y=gamma_g_1e-2, col sep=comma]{results/condnum_vs_gp_cut_configuration_variable.csv};
      \label{pgfplots:gamma1em2}
      \addplot [very thick,blue!70!black,mark=star]table[x=epsbyh, y=gamma_g_1e-3, col sep=comma]{results/condnum_vs_gp_cut_configuration_variable.csv};
      \label{pgfplots:gamma1em3}
      \addplot [very thick,black!70!black,mark=o]table[x=epsbyh, y=gamma_g_1e-4, col sep=comma]{results/condnum_vs_gp_cut_configuration_variable.csv};
      \label{pgfplots:gamma1em4}
      \addplot [very thick,yellow!70!black,mark=+]table[x=epsbyh, y=gamma_g_1, col sep=comma]{results/condnum_vs_gp_cut_configuration_variable.csv};
      \label{pgfplots:gamma1}
      \addplot [very thick,cyan!70!black,mark=diamond]table[x=epsbyh, y=gamma_g_0, col sep=comma]{results/condnum_vs_gp_cut_configuration_variable.csv};
      \label{pgfplots:gamma0}
      \coordinate (spypoint) at (axis cs:0.85, 8e5);
      \coordinate (magnifyglass) at (axis cs:0.24, 8e8);
\spy [width=2cm,height=1.5cm] on (spypoint) in node[fill=white] at (magnifyglass);      \nextgroupplot[align=left, title={{\scriptsize Number of iterations}},ytick={10,15,20,25,30,35}]
      \addplot [very thick,cyan!70!black,mark=diamond]table[x=epsbyh, y=numiter_0, col sep=comma]{results/mgiter_convrate_gp_cut_configuration_variable.csv};
      \addplot [very thick,red!70!black,mark=square]table[x=epsbyh, y=numiter_1e-1, col sep=comma]{results/mgiter_convrate_gp_cut_configuration_variable.csv};
      \addplot [very thick,green!70!black,mark=triangle]table[x=epsbyh, y=numiter_1e-2, col sep=comma]{results/mgiter_convrate_gp_cut_configuration_variable.csv};
      \addplot [very thick,blue!70!black, mark=star]table[x=epsbyh,y=numiter_1e-3,col sep=comma]{results/mgiter_convrate_gp_cut_configuration_variable.csv};
      \addplot [very thick,black!70!black,mark=o]table[x=epsbyh, y=numiter_1e-4, col sep=comma]{results/mgiter_convrate_gp_cut_configuration_variable.csv};
      \addplot [very thick,yellow!70!black,mark=+]table[x=epsbyh,y=numiter_1,col sep=comma]{results/mgiter_convrate_gp_cut_configuration_variable.csv};

\nextgroupplot[align=left, title={{\scriptsize Asymptotic convergence rate ($\rho^\ast$)}}]
      \addplot [very thick,cyan!70!black, mark=diamond] table[x=epsbyh,y=convrate_0,col sep=comma]{results/mgiter_convrate_gp_cut_configuration_variable.csv};
      \addplot [very thick,red!70!black, mark=square] table[x=epsbyh,y=convrate_1e-1,col sep=comma]{results/mgiter_convrate_gp_cut_configuration_variable.csv};
      \addplot [very thick,green!70!black, mark=triangle] table[x=epsbyh,y=convrate_1e-2,col sep=comma]{results/mgiter_convrate_gp_cut_configuration_variable.csv};
      \addplot [very thick,blue!70!black, mark=star]table[x=epsbyh,y=convrate_1e-3,col sep=comma]{results/mgiter_convrate_gp_cut_configuration_variable.csv};
      \addplot [very thick,black!70!black,mark=o]table[x=epsbyh, y=convrate_1e-4, col sep=comma]{results/mgiter_convrate_gp_cut_configuration_variable.csv};
      \addplot [very thick,yellow!70!black,mark=+]table[x=epsbyh,y=convrate_1,col sep=comma]{results/mgiter_convrate_gp_cut_configuration_variable.csv};

    \end{groupplot}
    \matrix [ draw, matrix of nodes, anchor = south east, node font=\scriptsize,
      column 1/.style={nodes={align=right,text width=.8cm}},
      column 2/.style={nodes={align=left,text width=1.0cm}},
      column 3/.style={nodes={align=right,text width=1cm}},
      column 4/.style={nodes={align=left,text width=1.5cm}},
      column 5/.style={nodes={align=right,text width=1cm}},
      column 6/.style={nodes={align=left,text width=1.5cm}},
      column 7/.style={nodes={align=right,text width=1cm}},
      column 8/.style={nodes={align=left,text width=1.5cm}},
      column 9/.style={nodes={align=right,text width=1cm}},
      column 10/.style={nodes={align=left,text width=1.5cm}},
      column 11/.style={nodes={align=right,text width=1cm}},
      column 12/.style={nodes={align=left,text width=1cm}},
    ] at ($(current axis.south east)+(-0.0,-1.5)$)
    {
      \ref{pgfplots:gamma1} & $\epsilon_G = 1 $ & \ref{pgfplots:gamma1em1} & $\epsilon_G = 10^{-1}$ & \ref{pgfplots:gamma1em2} & $\epsilon_G = 10^{-2}$ & \ref{pgfplots:gamma1em3} & $\epsilon_G = 10^{-3}$ & \ref{pgfplots:gamma1em4} & $\epsilon_G = 10^{-4}$ & \ref{pgfplots:gamma0} & $\epsilon_G= 0 $ \\
    };

  \end{tikzpicture}
  \caption{Comparing the condition number of the system matrix, number of total MG iterations and the asymptotic convergence rate of the MG method for different values of ghost penalty parameter, for \nameref{para:example1sc}.}
  \label{fig:moving_mesh_rho}
\end{figure*}
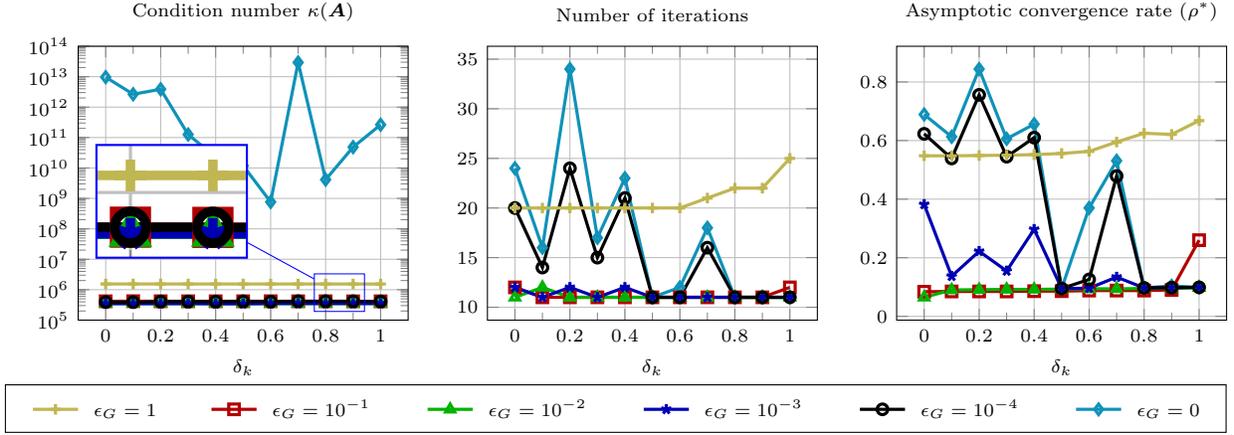

First, we compare the condition number of the stiffness matrix $\kappa(\matA)$ with respect to the different values of ghost penalty parameter.
We can see in Figure~\ref{fig:moving_mesh_rho} that the condition number of the system matrix is highly unstable if the ghost penalty term is not employed i.e., $\epsilon_G = 0$.
In the unfitted FE framework, the elements of the background mesh intersect the interface arbitrarily and on the elements with the small cuts, the gradients of the function are not bounded.
If we employ the stabilization term, the condition number of the system matrix becomes stable, and also it does not vary with respect to various cut configurations while translating the mesh.
We can observe that if smaller values of $\epsilon_G$ is chosen, e.g., $\epsilon_G=\{10^{-4}, 10^{-3},10^{-2},10^{-1}\}$, the condition number of the system matrix becomes stable.
In the next part, we compare the number of iterations of the MG method to reach the predefined tolerance criterion and compare also the asymptotic convergence rate.
We can see that, for $\epsilon_G= \{0, 10^{-4}\}$, the number iteration is not stable and varies with respect to different values of $\delta_k$.
But still, the number of iterations for $\epsilon_G=10^{-4}$ is smaller than the case without the ghost penalty term.
Similar behavior can also be witnessed in terms of the asymptotic convergence rate, as the asymptotic convergence rate $\rho^\ast$ is smaller for $\epsilon_G=10^{-4}$ than for $\epsilon_G = 0$.
For $\epsilon_G=10^{-3}$, the number of iterations is quite stable but the asymptotic convergence rate still oscillates with respect to moving mesh.
Interestingly, the number of iterations and the asymptotic convergence rate ($\rho^\ast$) is stable for $\epsilon_G =1 $, but these values are still considerably higher than the smaller values of the ghost penalty parameter.
The number of iterations and the asymptotic convergence rate are most stable for $\epsilon_G = \{10^{-2},10^{-1}\}$.
But with close observation, we can claim that for $\epsilon_G = 10^{-2}$, the condition number of the system matrix is the smallest and this also reflects in the performance of the multigrid method.
Due to this reason, the default value of the ghost penalty parameter in the previous and the next experiments is chosen as $\epsilon_G = 10^{-2}$.

\subsubsection{Performance of the Multigrid method}
In this section, the performance of the multigrid method is evaluated for increasing problem size and increasing the number of levels in the multigrid hierarchy.
Here, all the experiments are carried out on successively finer refinement levels $L1, L2, \ldots, L5$.
We employ $V(5,5)$-cycle and $W(5,5)$-cycle in the multigrid method, with the modified PGS method on the finest level and symmetric Gauss-Seidel method on the coarser levels as smoothers.
We increase the number of levels in the multigrid hierarchy with the refinement levels, i.e., MG employs $2$ levels for discretization level $L1$ while $6$ levels are employed for discretization level $L5$.

\begin{table*}[t]
  \centering
  \begin{tabular}{| c | c |  c c | c c | c c | c c | }
    \hline
         & \multirow{3}{*}{\shortstack[l]{\# levels}} & \multicolumn{4}{c|}{ \nameref{para:example1sc}}                   & \multicolumn{4}{c|}{ \nameref{para:example2sc}}                                                                                                                        \\ \cline{3-10}
         &                                            & \multicolumn{2}{c|}{ $V(5,5)$ } & \multicolumn{2}{c|}{ $W(5,5)$} & \multicolumn{2}{c|}{ $V(5,5)$} & \multicolumn{2}{c|}{ $W(5,5)$ }                                                     \\ \cline{3-10}
         &                                            & \# iter    & ($\rho^\ast$)      & \# iter    & ($\rho^\ast$)     & \# iter & ($\rho^\ast$)        & \# iter & ($\rho^\ast$) \\ \hline
    $L1$ & 2                                          & 10         & (0.072)            & 10         & (0.072)           & 11      & (0.092)              & 11      & (0.092)       \\
    \rowcolor{Gray}
    $L2$ & 3                                          & 11         & (0.061)            & 10         & (0.049)           & 12      & (0.083)              & 11      & (0.062)       \\
    $L3$ & 4                                          & 11         & (0.069)            & 10         & (0.054)           & 14      & (0.138)              & 12      & (0.097)       \\
    \rowcolor{Gray}
    $L4$ & 5                                          & 12         & (0.086)            & 11         & (0.046)           & 15      & (0.132)              & 13      & (0.086)       \\
    $L5$ & 6                                          & 15         & (0.136)            & 13         & (0.095)           & 17      & (0.182)              & 13      & (0.094)       \\ \hline
  \end{tabular}
  \caption{The number of iterations of the generalized multigrid method to reach a predefined tolerance for solving Signorini's problem.}
  \label{tab:QRMG_SC}
\end{table*}

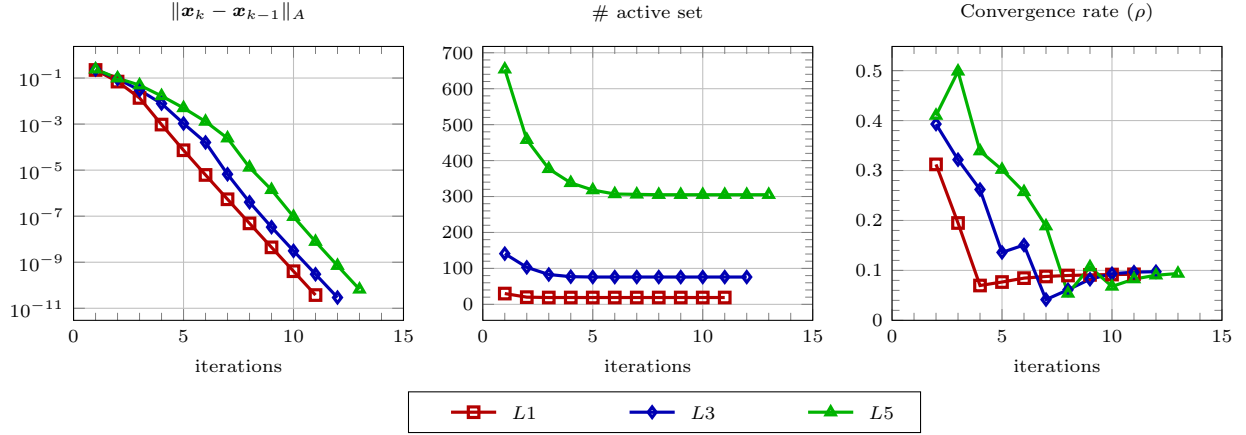
\begin{figure*}[t]
  \centering
  \tikzset{external/export next=false}
  \begin{tikzpicture}[]
    \begin{groupplot}[
        group style={
group size = 3 by 1,
x descriptions at=edge bottom,
horizontal sep=30pt,
          },
xtick={0,5,10,15},
        minor x tick num=4,
        xmin=0,
        xmax=15,
        scale=0.64,
        grid=major,
        xlabel={iterations},
        label style={font=\scriptsize}, tick label style={font=\scriptsize}, legend style={font=\scriptsize}
      ]
\nextgroupplot[align=left, ymode=log,title={\scriptsize $\|\vecx_k-\vecx_{k-1}\|_A$},ytick={1e-1,1e-3,1e-5,1e-7,1e-9,1e-11}]
      \addplot [very thick,red!70!black,mark=square]table[x=iter, y=correction_l1, col sep=comma]{results/rate_corr_activeset_sc3_W.csv};
      \label{pgfplots:L1}
\label{pgfplots:L2}
      \addplot [very thick,blue!70!black,mark=diamond]table[x=iter, y=correction_l3, col sep=comma]{results/rate_corr_activeset_sc3_W.csv};
      \label{pgfplots:L3}
\label{pgfplots:L4}
      \addplot [very thick,green!70!black,mark=triangle]table[x=iter, y=correction_l5, col sep=comma]{results/rate_corr_activeset_sc3_W.csv};
      \label{pgfplots:L5}

\nextgroupplot[align=left, title={\scriptsize \# active set}, minor y tick num=4,ytick={0,100,200,300,400,500,600,700}]
      \addplot [very thick,red!70!black,mark=square]table[x=iter, y=AS_l1, col sep=comma]{results/rate_corr_activeset_sc3_W.csv};
\addplot [very thick,blue!70!black,mark=diamond]table[x=iter, y=AS_l3, col sep=comma]{results/rate_corr_activeset_sc3_W.csv};
\addplot [very thick,green!70!black,mark=triangle]table[x=iter, y=AS_l5, col sep=comma]{results/rate_corr_activeset_sc3_W.csv};

\nextgroupplot[align=left, title={\scriptsize Convergence rate ($\rho$)}, minor y tick num=4, ymin=0,ytick={0,0.1,0.2,0.3,0.4,0.5}]
      \addplot [very thick,red!70!black,mark=square]table[x=iter, y=rate_l1, col sep=comma]{results/rate_corr_activeset_sc3_W.csv};
\addplot [very thick,blue!70!black,mark=diamond]table[x=iter, y=rate_l3, col sep=comma]{results/rate_corr_activeset_sc3_W.csv};
\addplot [very thick,green!70!black,mark=triangle]table[x=iter, y=rate_l5, col sep=comma]{results/rate_corr_activeset_sc3_W.csv};

    \end{groupplot}
    \matrix [ draw, matrix of nodes, anchor = south east, node font=\scriptsize,
      column 1/.style={nodes={align=right,text width=1cm}},
      column 2/.style={nodes={align=left,text width=1cm}},
      column 3/.style={nodes={align=centre,text width=1.5cm}},
      column 4/.style={nodes={align=right,text width=1cm}},
      column 5/.style={nodes={align=left,text width=1cm}},
      column 6/.style={nodes={align=left,text width=1.5cm}},
      column 7/.style={nodes={align=right,text width=1cm}},
      column 8/.style={nodes={align=left,text width=0.5cm}},
    ] at ($(current axis.south east)+(-4,-1.5)$)
    {
      \ref{pgfplots:L1} & $L1$ &  & \ref{pgfplots:L3} & $L3$ &  & \ref{pgfplots:L5} & $L5$ \\
    };

  \end{tikzpicture}
  \caption{The history of the correction in the energy norm, active set, the convergence rate ($\rho$) at each iteration of the MG method with $W(5,5)$-cycle for solving \nameref{para:example2sc}.}
  \label{fig:compare_MGSC}
\end{figure*}

Table~\ref{tab:QRMG_SC} illustrates the number of iterations of the generalized multigrid method to reach the termination criterion \eqref{eq:termination_cont} and the asymptotic convergence rate.
From Table~\ref{tab:QRMG_SC}, it is evident that the number of iterations does not vary significantly with the increasing problem size and an increasing number of levels in the hierarchy.
As expected, the MG method with the $W$-cycle outperforms the $V$-cycle.
While using the $V$-cycle the asymptotic convergence rate increases slightly with increasing problem size, but if we employ the $W$-cycle the asymptotic convergence rate becomes stable ($\rho^\ast < 0.1$).
In Figure~\ref{fig:compare_MGSC}, we compare the convergence history of the correction, the size of the active set, and the convergence rate at each iteration for $L1, L3, L5$, while using the MG method with $W$-cycle.
From Figure~\ref{fig:compare_MGSC}, we can see that a few initial iterations are spent for identifying the active set.
Due to this reason, we can see the convergence rate is quite high and correction in the energy norm is also reducing slowly for a few initial iterations.
But once the active set is identified, the norm of correction reduces at the same rate and the convergence rate also becomes stable.
As we increase the problem size, the size of the active set also increases, and a few more iterations are required to identify the exact active set.
Due to this reason, we observe a small increase in the number of iterations with increasing problem size.

\subsection{Two-body Contact Problem}
For the two-body contact problem, we consider two different types of embedded interfaces: a circular interface and an elliptical interface and different material parameters.

\subsubsection{Problem Description}
The experiments in this section are carried out on a structured background mesh with quadrilateral elements.
Both bodies are considered to be in contact with each other in the absence of external forces. The background mesh is given in square domain $\Omega$ in $(0,1)^2$.
We start with $50$ elements in each direction, this mesh is denoted by $\widetilde{\pT}_0$ and it is associated with level $L0$.
We create a hierarchy of meshes by uniformly refining this mesh until we have $1600$ elements in each direction, this mesh is defined as $\widetilde{\pT}_5$.
The Dirichlet boundary conditions is defined as $\vecu = (0,0)$ on $x=[0,1]$ and $y=0$, while the Neumann boundary condition is defined as $\bs{\sigma}\vecn = (0, 5)$ on  $x=[0,1]$ and $y=1$.
The body force for this example is considered to be zero.

\paragraph{Example 1-TC}\label{para:example1TC}
For this example, we consider a circular contact interface denoted as $\Gamma_c$.
The circular interface is defined as a zero level set of a function $\Lambda_c(\matX):=r_{t_1}^2 - \|\matX - \vecc_{t} \|^2_2$, with radius $r_{t_1}^2 = 3 - 2 \cdot 2^{1/2}$, and $\vecc_t $ is the center of the circle, chosen as $(0.5,0.5)$.
The circular interface decomposes the domain $\Omega$ into $\Omega_1$, where $\Lambda_c(\matX) > 0$ and $\Omega_2$ where $\Lambda_c(\matX) < 0$.
For this example, we consider two different sets of material parameters.
We choose Young's modulus as $E_1=10\,\text{MPa}$ and $E_2 = \{10\,\text{MPa}, 50\,\text{MPa}\}$ and the Poisson's ratio is chosen as $\nu_1 = \nu_2 =0.3$.

\paragraph{Example 2-TC}\label{para:example2TC}
This example considers an elliptical contact interface denoted as $\Gamma_e$.
The interface is defined as a zero level set of a function
\[
  \Lambda_e(\matX) := r_{t_2}^2 - \Bigl| \frac{x-c_x}{a} \Bigr|^2 - \Bigl| \frac{y-c_y}{b} \Bigr|^2.
\]
Here, $r_{t_2}$ denotes the radius of the ellipse, chosen as $r_{t_2}^2 = 2(3 - 2\cdot 2^{1/2})$.
The symbols $a$ and $b$ denote the major and minor axis of the ellipse, chosen as ${a=1, b=0.8}$.
Here the center of the circle is chosen as $(0.5,0.5)$.
The circular interface decomposes the domain $\Omega$ into $\Omega_1$, where $\Lambda_e(\matX) > 0$ and $\Omega_2$, where $\Lambda_e(\matX) < 0$.
For this example, we consider the same set of material parameters, as used in the previous example.
Young's modulus is chosen as $E_1=10\,\text{MPa}$, $E_2 = \{10\,\text{MPa}, 50\,\text{MPa}\}$ and the Poisson's ratio $\nu_1=\nu_2 =0.3$.

\begin{figure*}[t] \begin{subfigure}{0.32\textwidth}
    \includegraphics[width=\linewidth]{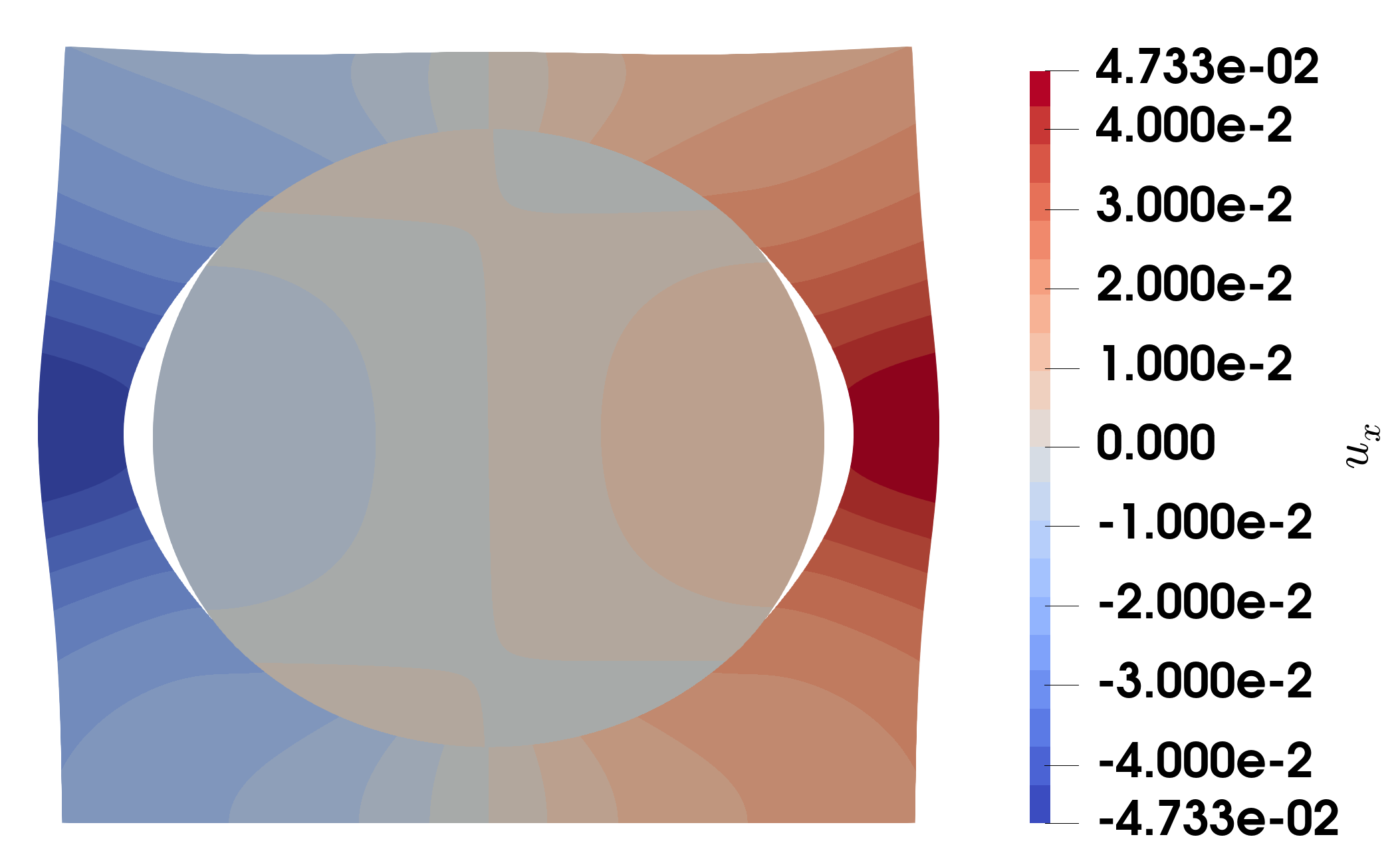}
    \caption{Displacement in $X$-direction $\vecu_x$.}\end{subfigure}\hfill
  \begin{subfigure}{0.32\textwidth}
    \includegraphics[width=\linewidth]{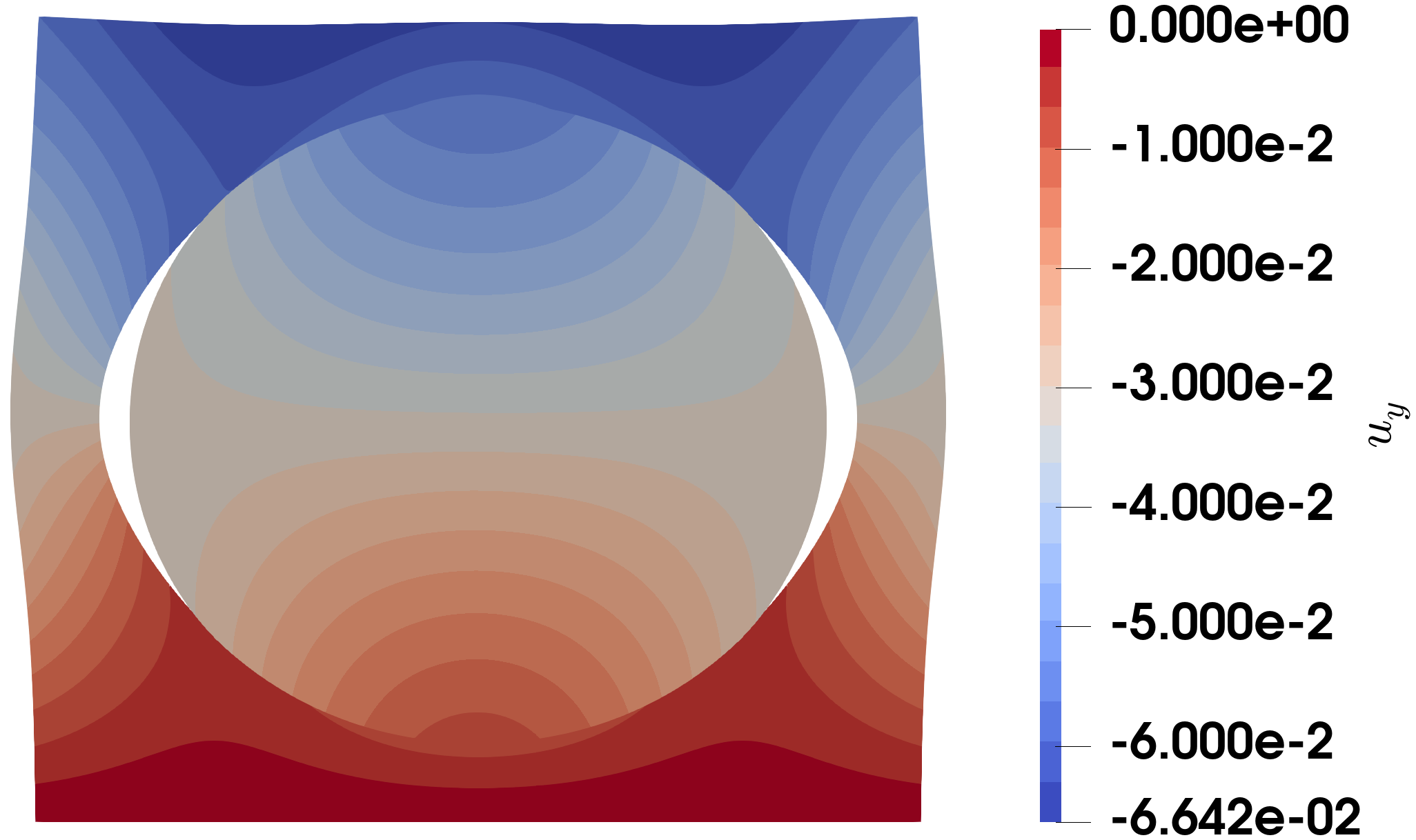}
    \caption{Displacement in $Y$-direction $\vecu_y$.}\end{subfigure}\hfill
  \begin{subfigure}{0.32\textwidth}
    \includegraphics[width=\linewidth]{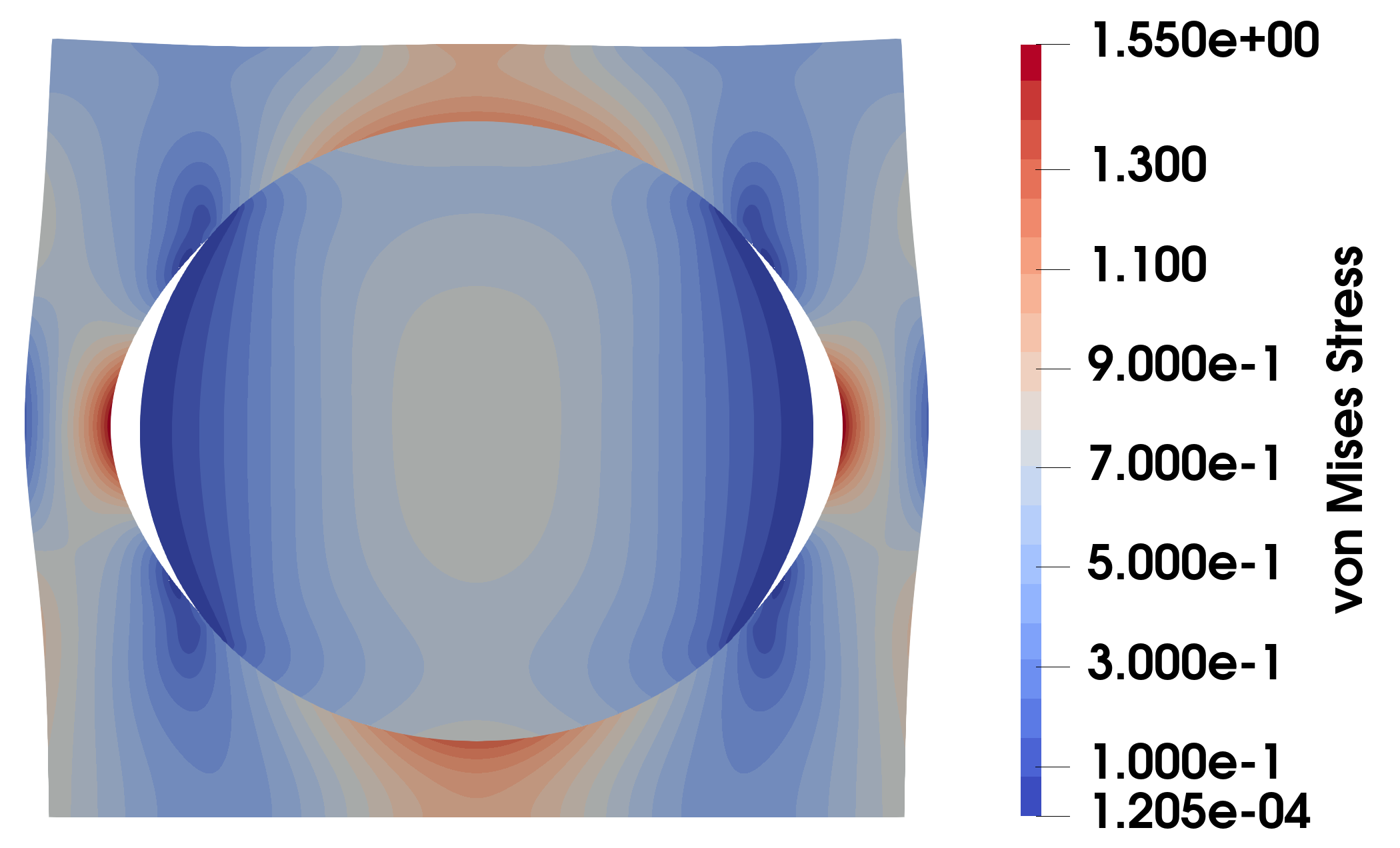}
    \caption{von Mises Stress.}\end{subfigure}
  \par\bigskip \begin{subfigure}{0.32\textwidth}
    \includegraphics[width=\linewidth]{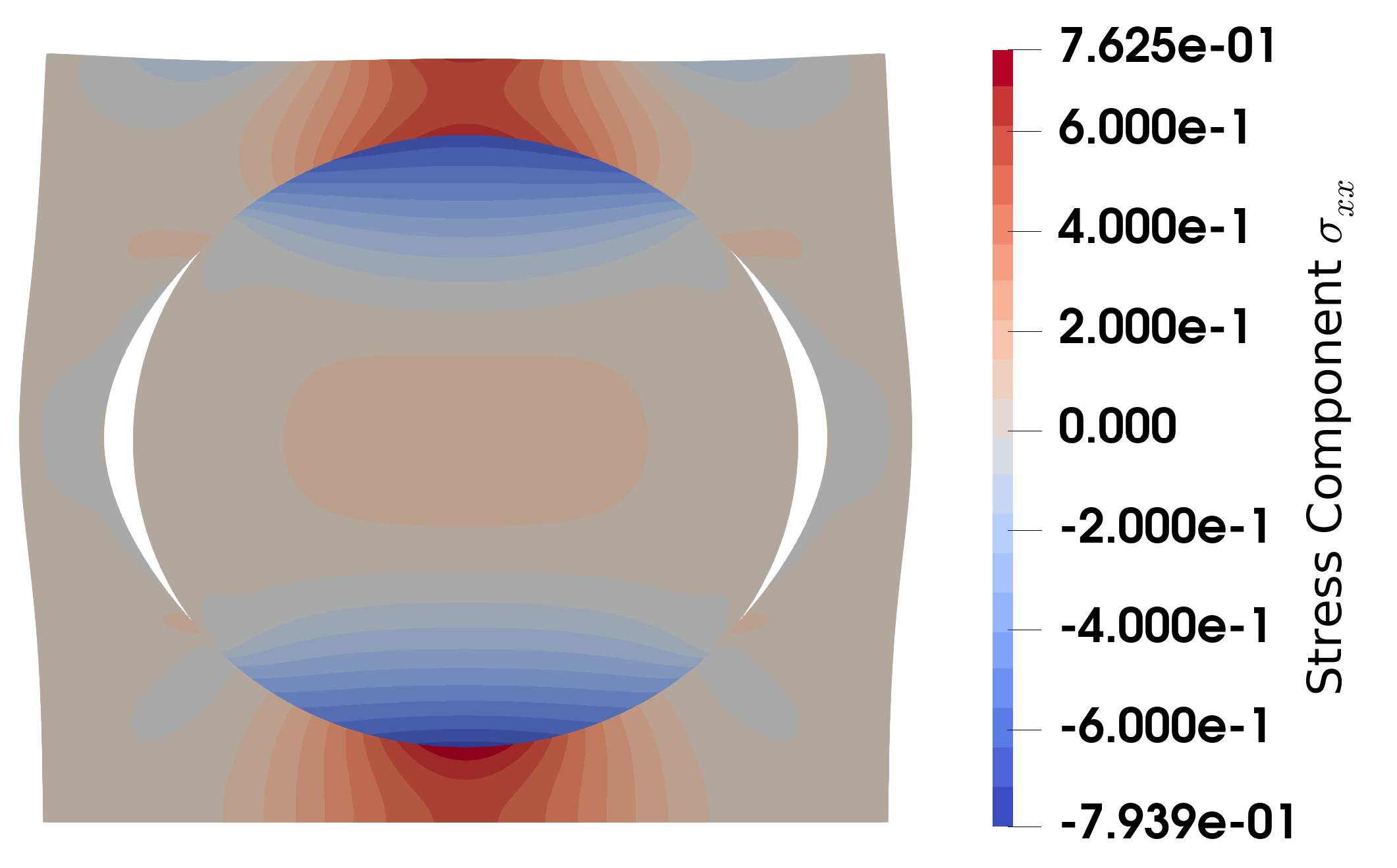}
    \caption{Stress component $\bs{\sigma}_{xx}$.} \end{subfigure}\hfill
  \begin{subfigure}{0.32\textwidth}
    \includegraphics[width=\linewidth]{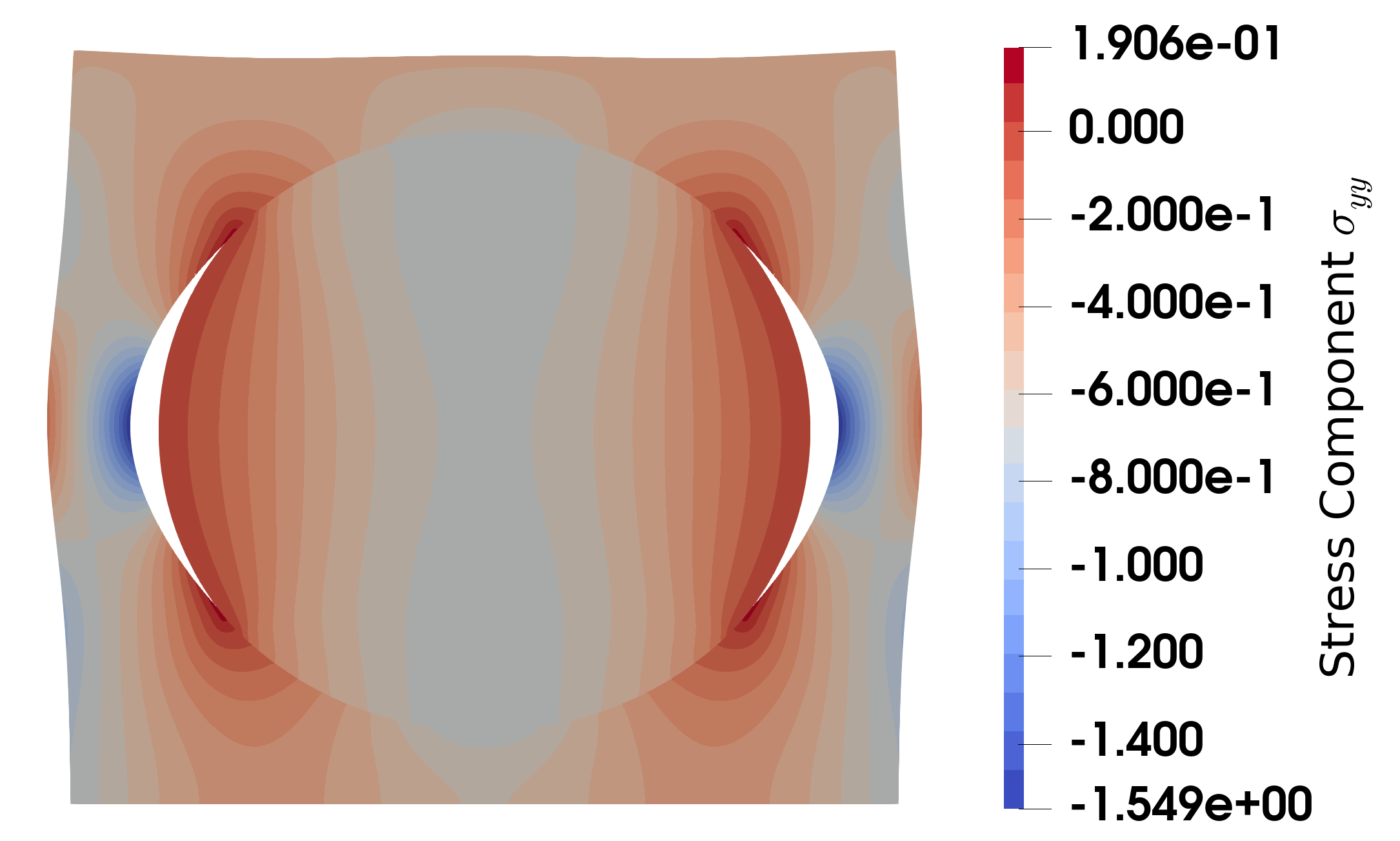}
    \caption{Stress component $\bs{\sigma}_{yy}$.} \end{subfigure}\hfill
  \begin{subfigure}{0.32\textwidth}
    \includegraphics[width=\linewidth]{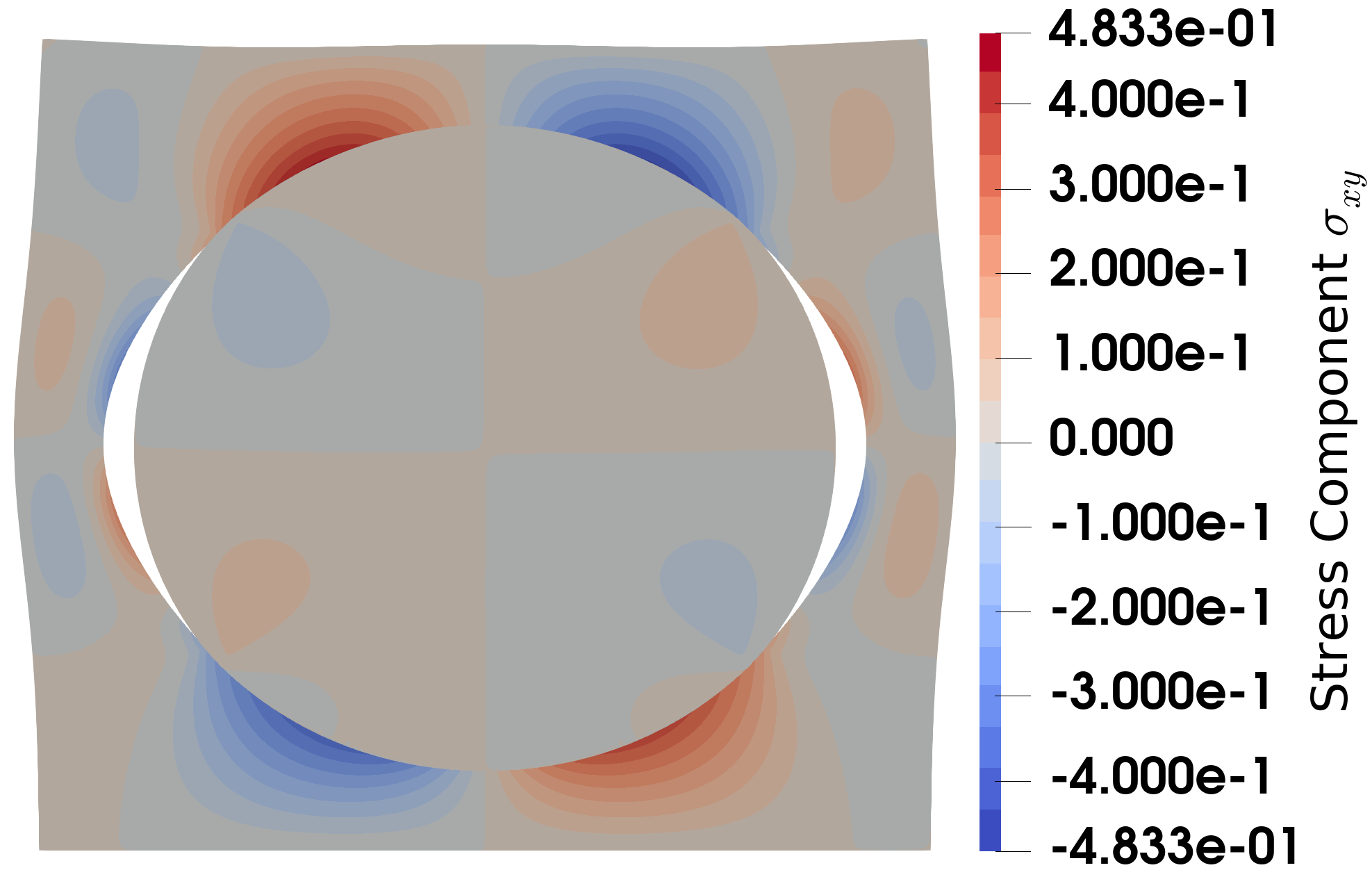}
    \caption{Stress component $\bs{\sigma}_{xy}$.} \end{subfigure}
  \caption{Resultant displacement field and stress field, as a solution of the two-body contact problem, \nameref{para:example1TC}, with Young's modulus $E_1=E_2=10\,\text{MPa}$, where the domain $\Omega_2$ is a circle.}
  \label{fig:TB_circle}
\end{figure*}
\begin{figure*}[!hbt] \begin{subfigure}{0.32\textwidth}
    \includegraphics[width=\linewidth]{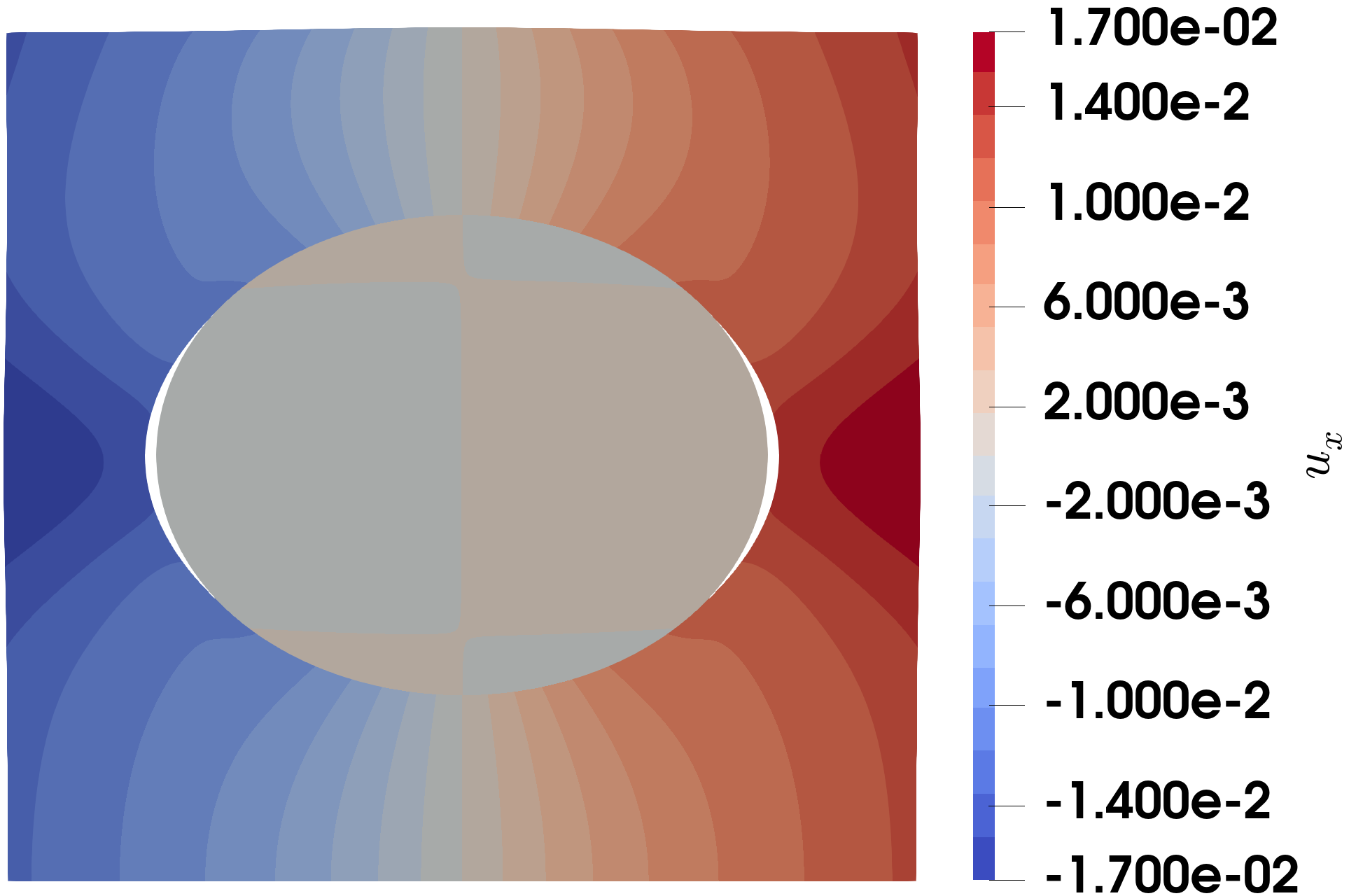}
    \caption{Displacement in $X$-direction $\vecu_x$.}\end{subfigure}\hspace*{\fill}
  \begin{subfigure}{0.32\textwidth}
    \includegraphics[width=\linewidth]{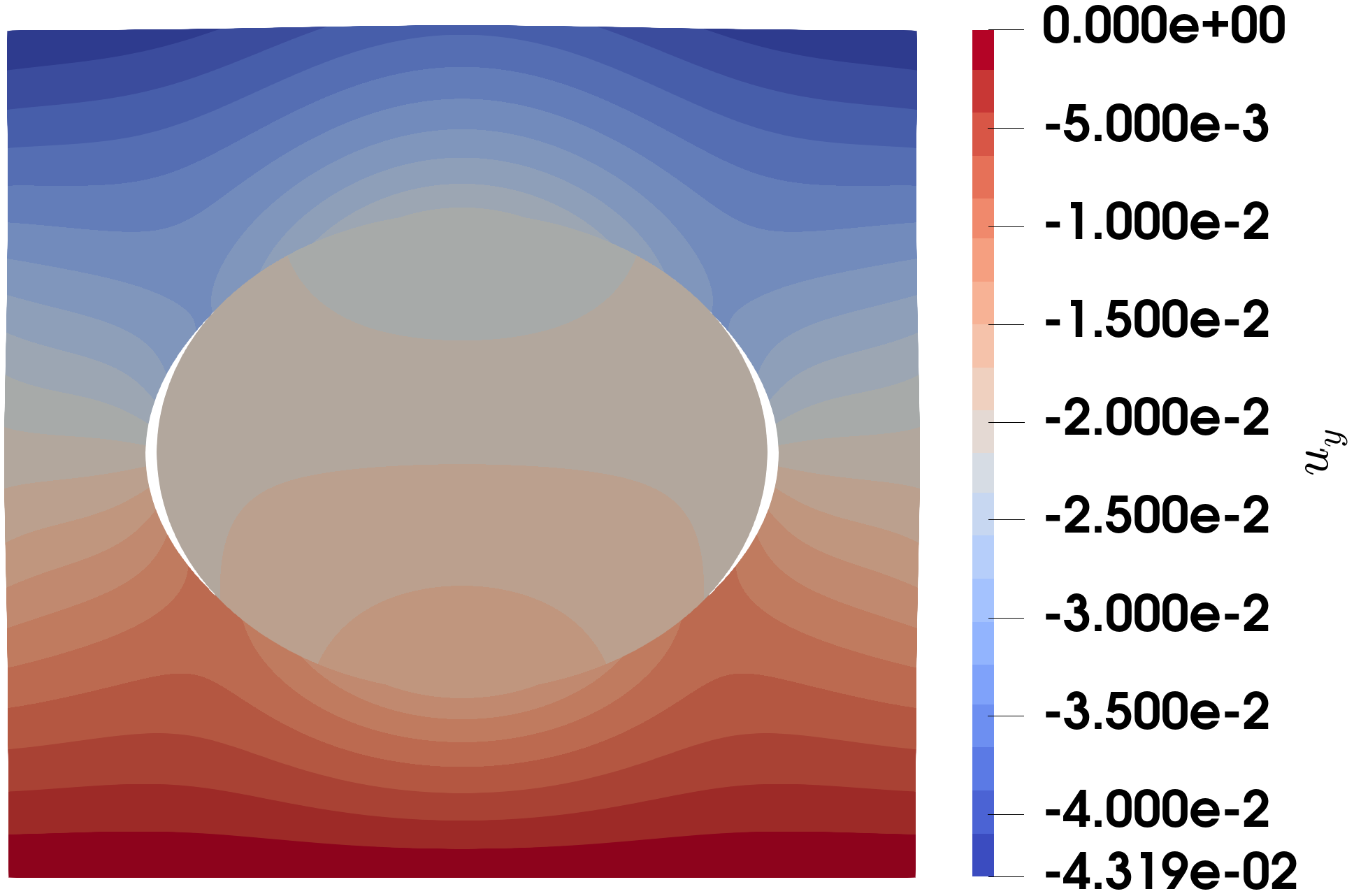}
    \caption{Displacement in $Y$-direction $\vecu_y$.} \end{subfigure}\hspace*{\fill}
  \begin{subfigure}{0.32\textwidth}
    \includegraphics[width=\linewidth]{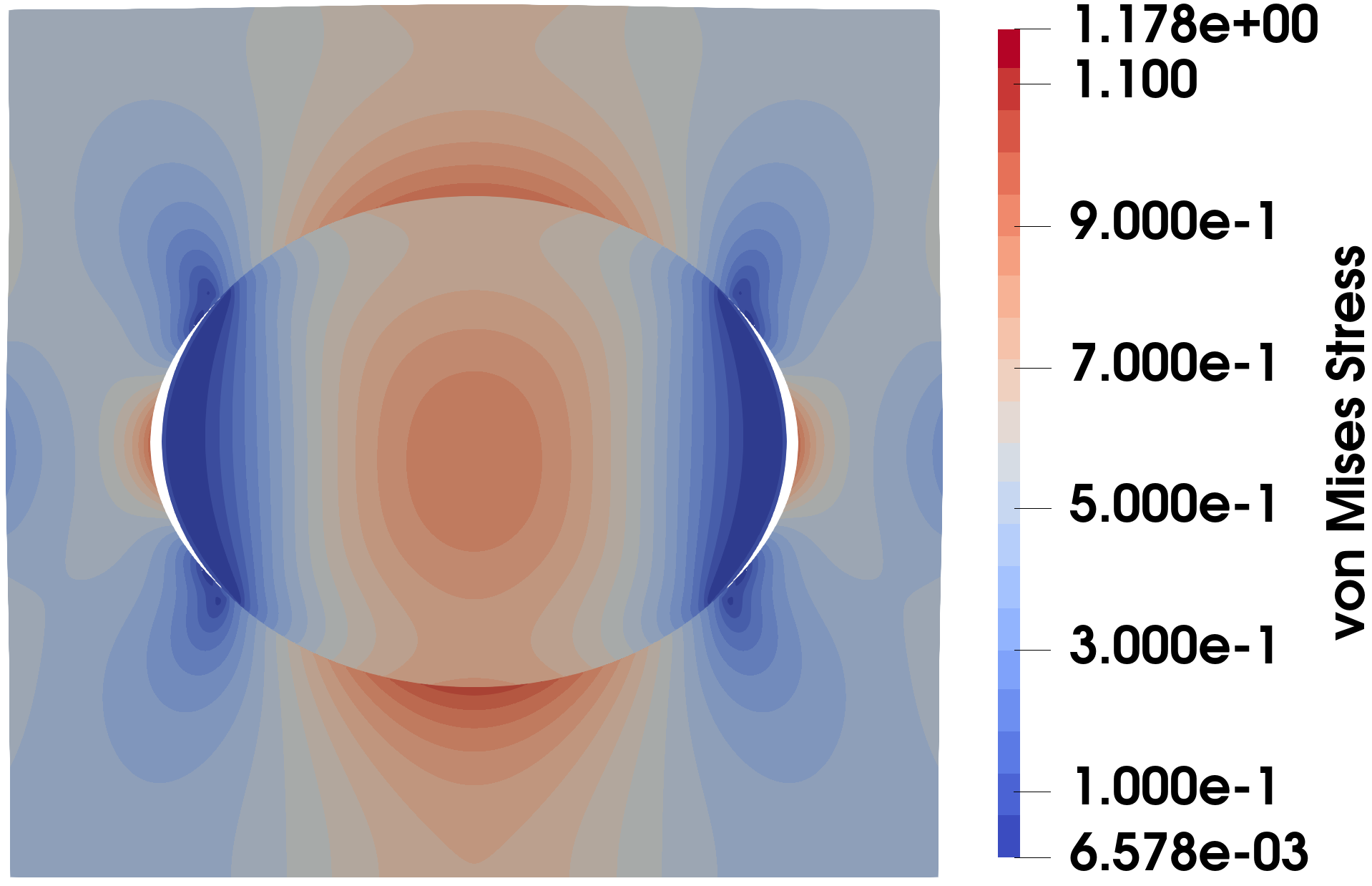}
    \caption{von Mises Stress.}\end{subfigure}

  \par\bigskip \begin{subfigure}{0.32\textwidth}
    \includegraphics[width=\linewidth]{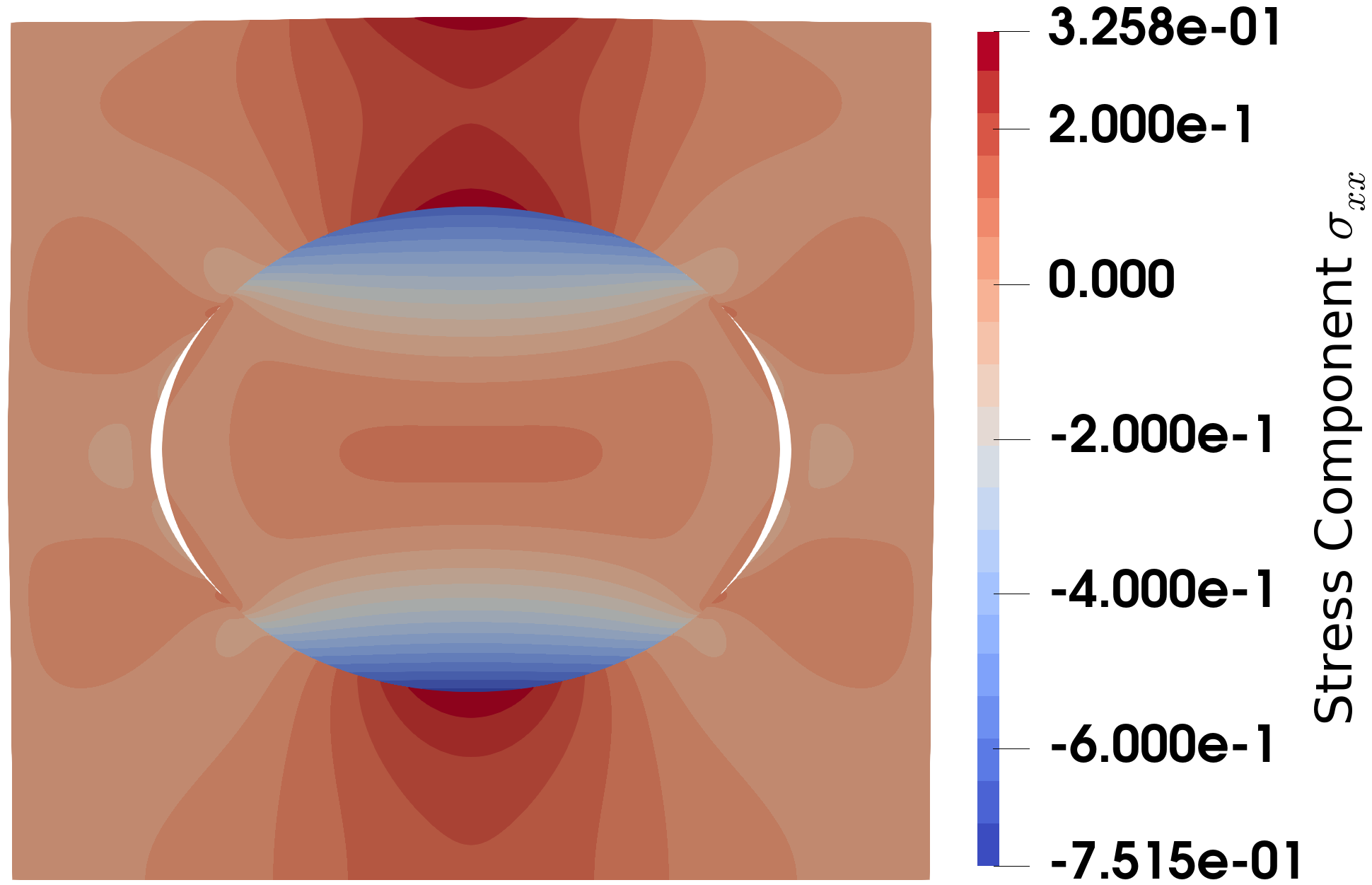}
    \caption{Stress component $\bs{\sigma}_{xx}$.}\end{subfigure}\hspace*{\fill}
  \begin{subfigure}{0.32\textwidth}
    \includegraphics[width=\linewidth]{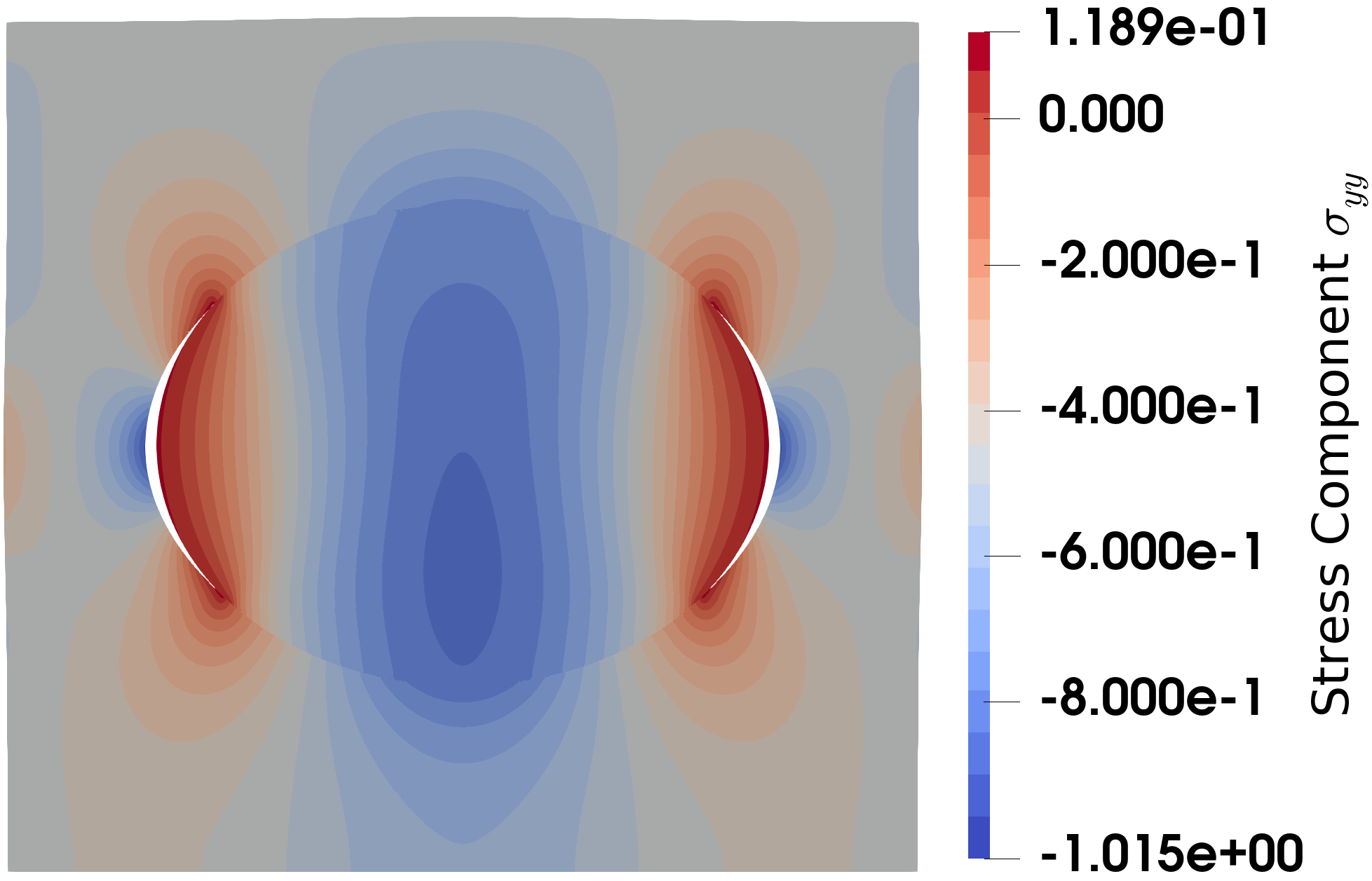}
    \caption{Stress component $\bs{\sigma}_{yy}$.}\end{subfigure}\hspace*{\fill}
  \begin{subfigure}{0.32\textwidth}
    \includegraphics[width=\linewidth]{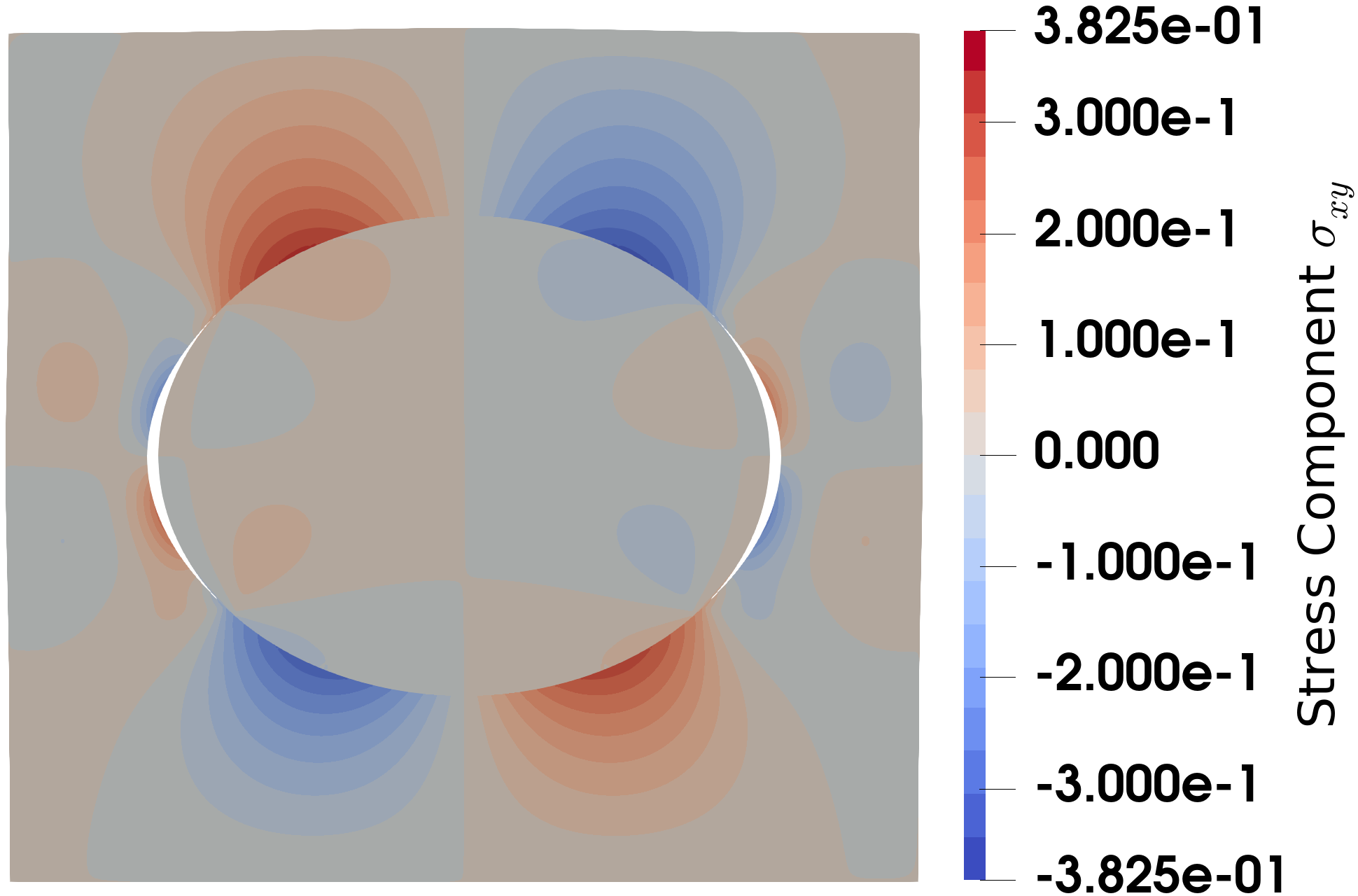}
    \caption{Stress component $\bs{\sigma}_{xy}$.}\end{subfigure}\hspace*{\fill}
  \caption{Resultant displacement field and stress field, as a solution of the two-body contact problem, \nameref{para:example2TC}, with Young's modulus $E_1=10\,\text{MPa}$ and $E_2=50\,\text{MPa}$, where the domain $\Omega_2$ is an ellipse.}
  \label{fig:TB_ellipse}
\end{figure*}

\begin{figure*}[t] \begin{subfigure}{0.23\textwidth}
    \includegraphics[trim={6.5cm 7cm 5.2cm 7cm},width=\linewidth]{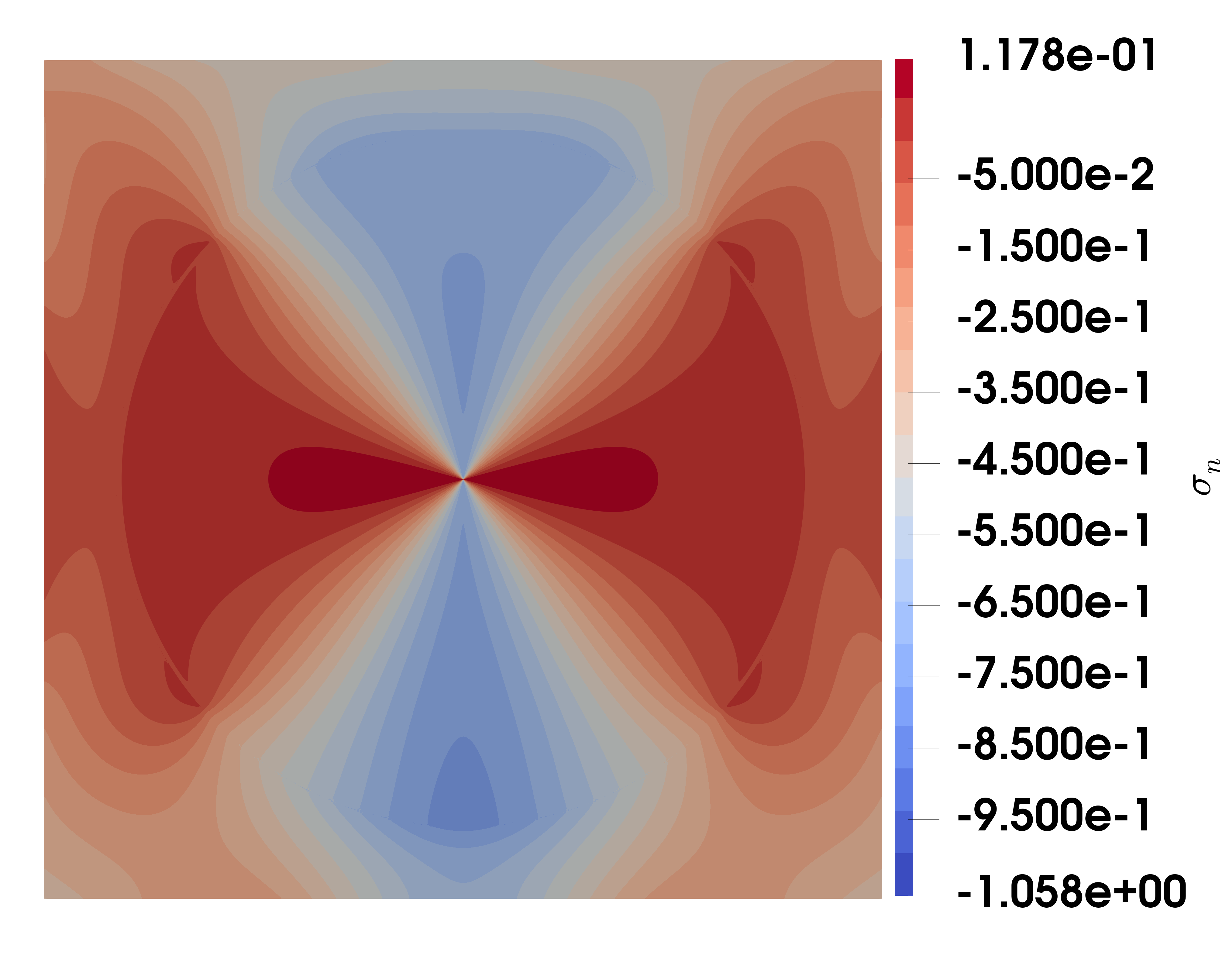}
    \caption{\nameref{para:example1TC}, \\ $E_1=10\,\text{MPa}$, $E_2=10\,\text{MPa}$.}\end{subfigure}\hspace*{\fill}
  \begin{subfigure}{0.23\textwidth}
    \includegraphics[trim={6.5cm 7cm 5.2cm 7cm},width=\linewidth]{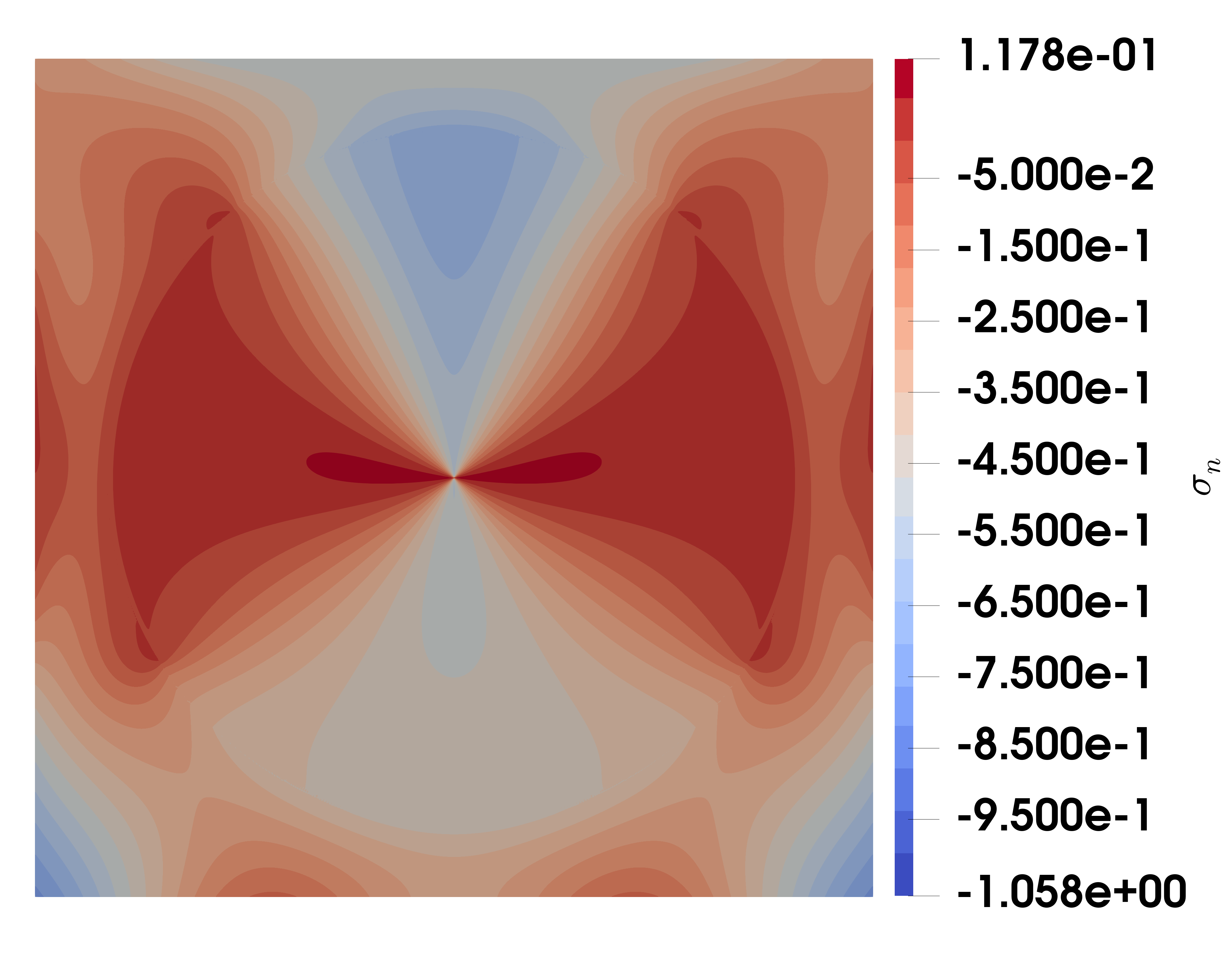}
    \caption{\nameref{para:example1TC}, \\ $E_1=10\,\text{MPa}$, $E_2=50\,\text{MPa}$.}\end{subfigure}\hspace*{\fill}
  \begin{subfigure}{0.23\textwidth}
    \includegraphics[trim={6.5cm 7cm 5.2cm 7cm},width=\linewidth]{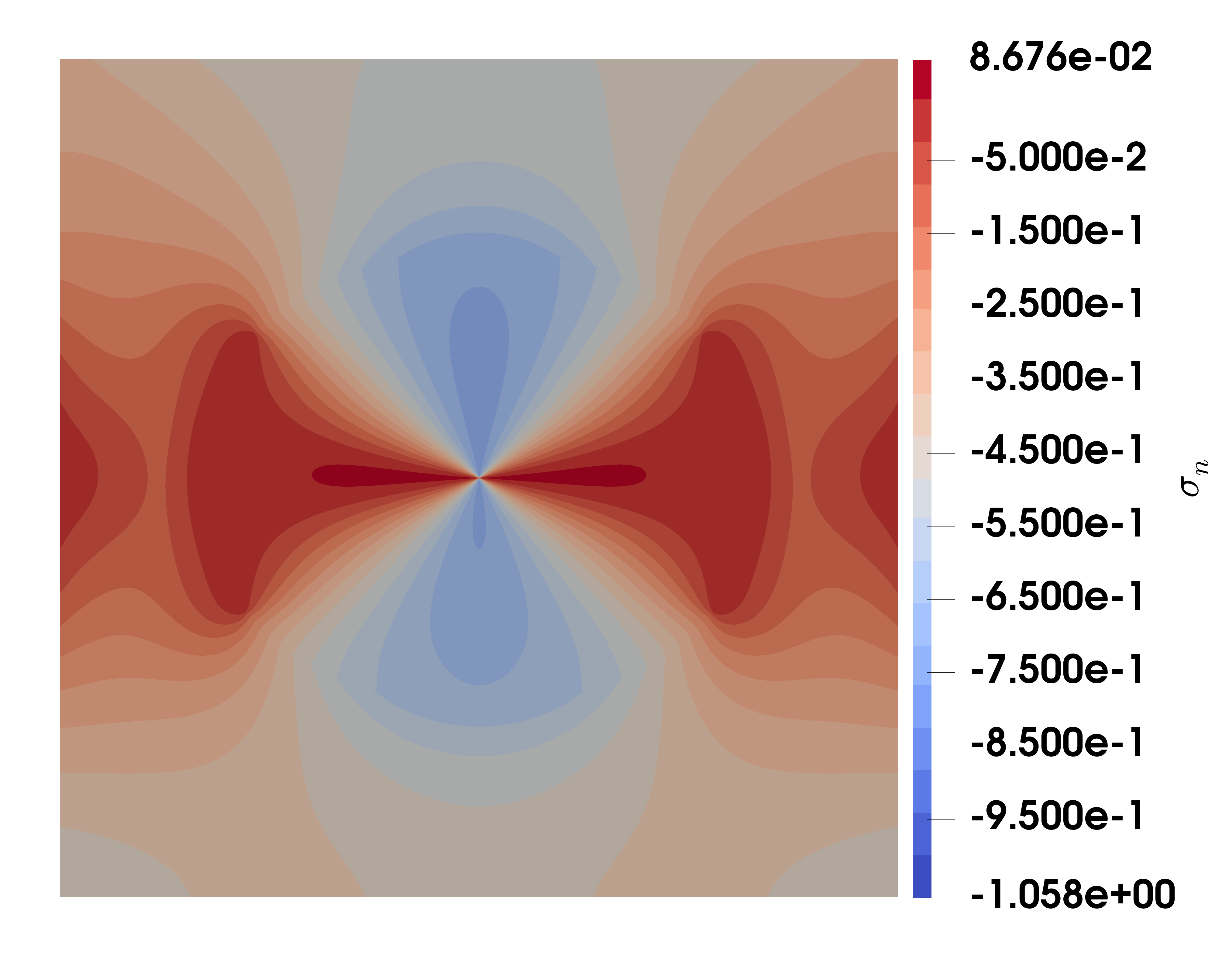}
    \caption{\nameref{para:example2TC},\\ $E_1=10\,\text{MPa}$, $E_2=10\,\text{MPa}$.}\end{subfigure}\hspace*{\fill}
  \begin{subfigure}{0.23\textwidth}
    \includegraphics[trim={6.5cm 7cm 5.2cm 7cm},width=\linewidth]{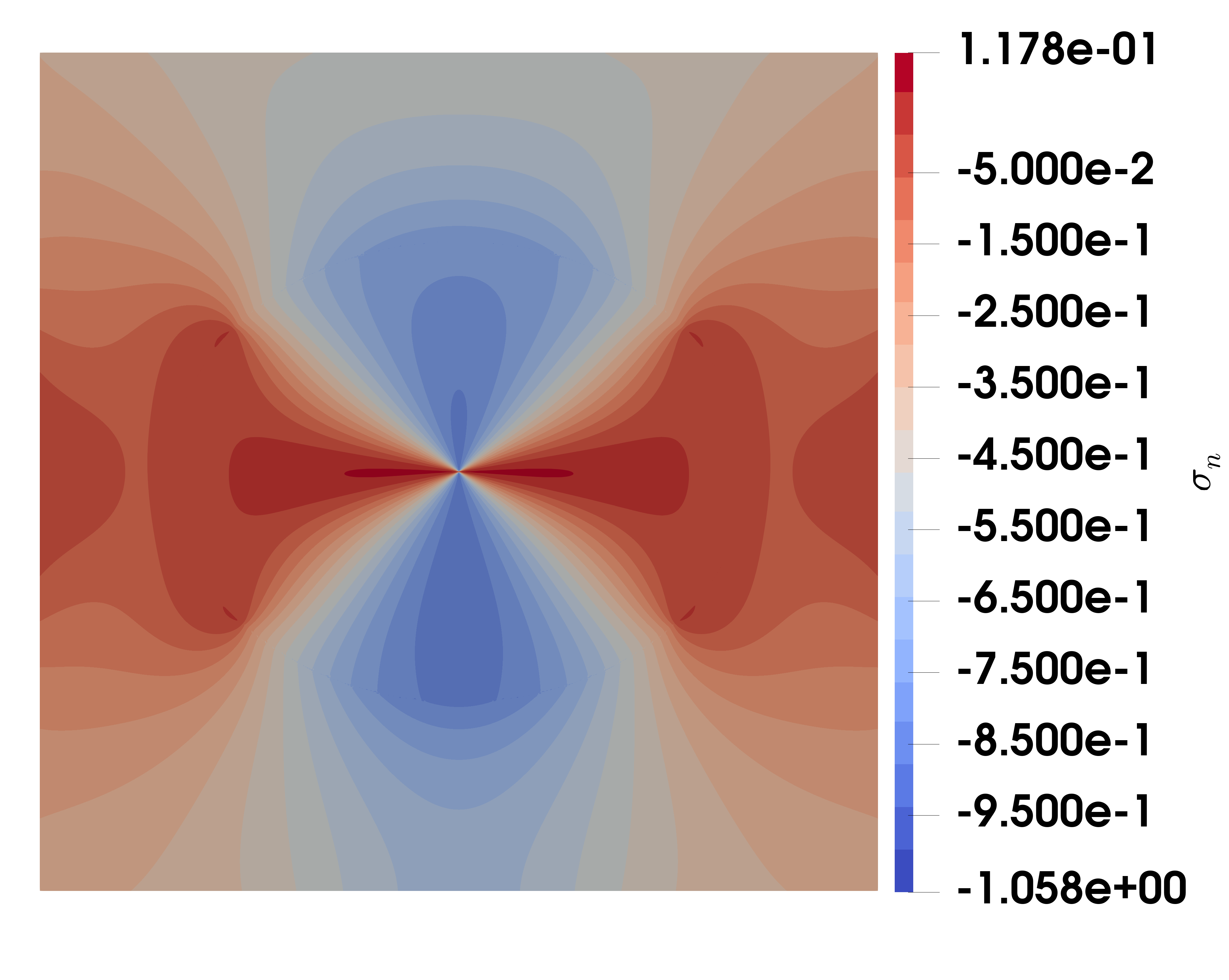}
    \caption{\nameref{para:example2TC}, \\ $E_1=10\,\text{MPa}$, $E_2=50\,\text{MPa}$.}\end{subfigure}\hspace*{\fill}
      \caption{Resultant normal stresses ($\sigma_n$) computed on nodes for example two-body contact problems in reference configuration.}
  \label{fig:TB_normal_stress}
\end{figure*}
In Figure~\ref{fig:TB_circle}, we can observe the resultant displacement field and the stress fields for two-body problem with circular interface \nameref{para:example1TC} with $E_1=E_2=10\,\text{MPa}$.
Similarly, the result of the two-body contact problem with an elliptic interface \nameref{para:example2TC} with $E_1=50\,\text{MPa}$ and $E_2=10\,\text{MPa}$ can be seen in Figure~\ref{fig:TB_ellipse}.
We can observe that the circle and the ellipse are in contact with the surrounding block on the top and bottom.
In Figure~\ref{fig:TB_normal_stress}, we can observe the resultant normal stresses computed on the whole domain for all exampled.
From Figure~\ref{fig:TB_normal_stress}, it is clear that the resultant normal stresses on the embedded interfaces are continuous and they are negative where the two bodies are in contact with each other.
\subsubsection{Performance of the Multigrid Method}
In this section, all the experiments are carried out on increasing problem size and with an increasing number of levels in the multigrid hierarchy.
We employ the multigrid method with $W(5,5)$-cycle with the modified PGS method on the finest level and symmetric Gauss-Seidel method on coarse levels as smoothers.
Table~\ref{tab:QRMG_TB} shows the number of iterations of the generalized multigrid method to reach the termination criterion \eqref{eq:termination_cont}.
We can conclude from the table, that the number of iterations does not change with an increasing number of levels in the multigrid hierarchy.
Also, in Table~\ref{tab:QRMG_SC}, we can observe the asymptotic convergence rate of the multigrid method.
We can see the difference in the asymptotic converge rates, despite the number of iterations required to reach the termination criterion is almost the same.
For the case with a homogeneous value of Young's modulus, the asymptotic convergence rate is quite low ($\rho^\ast < 0.1$).
While, for the case with different values of Young's modulus, the asymptotic convergence rate is much larger ($\rho^\ast < 0.55$) for the circular interface, whereas this value is much smaller for the elliptical interface ($\rho^\ast < 0.15$).
Even though the asymptotic convergence rate is higher for the circular interface with different Young's modulus, the number of iterations and the asymptotic convergence rate do not increase with increasing problem size.
Thus, we can conclude that the proposed generalized multigrid method is robust with respect to the number of levels, the material parameters, type of obstacle or the shape of the interface.
\begin{table*}[t]
  \centering
  \begin{tabular}{ | c | c |  c c | c c | c c | c c | }
    \hline
         & \multirow{3}{*}{\shortstack[l]{\# levels}} & \multicolumn{4}{c|}{{\nameref{para:example1TC}} } & \multicolumn{4}{c|}{{\nameref{para:example2TC}} }  \\ \cline{3-10}
         &  & \multicolumn{2}{c|}{{$E_1=10, E_2=10$} }  & \multicolumn{2}{c|}{{$E_1=10, E_2=50$}} & \multicolumn{2}{c|}{{$E_1=10, E_2=10$} } & \multicolumn{2}{c|}{{$E_1=10, E_2=50$} } \\ \cline{3-10}
&  & \# iter  & ($\rho^\ast$)    & \# iter                                     & ($\rho^\ast$)                               & \# iter & ($\rho^\ast$) & \# iter & ($\rho^\ast$) \\ \hline
    $L1$ & 2                                          & 11                                                & (0.060)     & 16    & (0.526)   & 10      & (0.078)       & 13      & (0.142)       \\
    \rowcolor{Gray}
    $L2$ & 3                                          & 11                                                & (0.055)     & 19    & (0.488)   & 11      & (0.075)       & 14      & (0.134)       \\
    $L3$ & 4                                          & 12                                                & (0.064)     & 13    & (0.165)   & 11      & (0.065)       & 14      & (0.120)       \\
    \rowcolor{Gray}
    $L4$ & 5                                          & 13                                                & (0.057)     & 18    & (0.533)   & 13      & (0.059)       & 13      & (0.087)       \\
    $L5$ & 6                                          & 14                                                & (0.073)     & 16    & (0.533)   & 13      & (0.109)       & 13      & (0.088)       \\ \hline
  \end{tabular}
  \caption{The number of iterations of the generalized multigrid method ($W(5,5$)-cycle) to reach a predefined tolerance for solving two-body contact problems.}
  \label{tab:QRMG_TB}
\end{table*}

\subsubsection{Comparison with other solution methods}
In this section, we compare the performance of the generalized MG method with the other solution strategies such as the semi-smooth Newton (SSN) method~\cite{michaelulbrich2011-07-13} and interior-point (IP) method~\cite{jorgenocedal2000-04-27}.
For this comparison, we use the MG method with $W(5,5)$-cycle and increasing number of levels in the multigrid hierarchy as discussed earlier.
We use the predictor-corrector variant of the IP method~\cite{mehrotra_implementation_1992}, hence for each iteration the linear system of equations is solved twice.
The IP method is used in the reduced form, hence at each iteration, the linear system with $nd$ unknowns has to be solved.
The SSN method for the linear inequality constraints can not be formulated in a reduced form.
Hence, we need to solve the enhanced KKT-system in each SSN iteration, where the linear system is formed as a non-symmetric system with $nd+m$ unknowns.

From Table~\ref{tab:compare_solvers}, we can see that for all the numerical experiments the multigrid method is the cheapest amongst all solution methods.
The multigrid method requires between $10$-$20$ iterations to reach the termination criterion.
For the SSN method, the number of iterations required to converge is smaller than the IP method.
It is not easy to make a comparison between the IP method and SSN method, as arising linear system of equations has a different structure.
But we can safely claim that the multigrid method is at least an order of magnitude times cheaper than the IP method and SSN method.
This is due to the reason that, despite the SSN method and IP method in some cases require fewer iterations than the multigrid method, the computational complexity of these methods per iteration is significantly higher.

\begin{table*}[t]
  \centering
  \begin{tabular}{| c |  c c c | c c c | c c c | c c  c |} \hline
         & \multicolumn{3}{c|}{{\nameref{para:example1sc}}} & \multicolumn{3}{c|}{{\nameref{para:example2sc}} } & \multicolumn{3}{c|}{{\nameref{para:example1TC}}} & \multicolumn{3}{c|}{{\nameref{para:example2TC}} }  \\ \cline{2-13}
         & \multicolumn{3}{c|}{{$E=10$} }                   & \multicolumn{3}{c|}{{$E = 10$} }                  & \multicolumn{3}{c|}{{$E_1=10, E_2=50 $} }        & \multicolumn{3}{c|}{{$E_1=10, E_2=50$}} \\ \cline{2-13}
         & GMG                                               & SSN                                               & IP                       & GMG  & SSN & IP & GMG & SSN & IP & GMG & SSN & IP \\ \hline
    $L1$ & 10                                               & 9                                                 & 17                       & 11  & 8   & 17 & 16 & 8   & 11 & 13 & 8   & 11 \\
    \rowcolor{Gray}
    $L2$ & 10                                               & 9                                                 & 17                       & 11  & 9   & 18 & 19 & 8   & 11 & 14 & 7   & 12 \\
    $L3$ & 10                                               & 11                                                & 18                       & 12  & 11  & 18 & 13 & 11  & 11 & 14 & 9   & 11 \\
    \rowcolor{Gray}
    $L4$ & 11                                               & 13                                                & 18                       & 13  & 16  & 19 & 18 & 11  & 11 & 13 & 10  & 11 \\
    $L5$ & 13                                               & 14                                                & 18                       & 13  & 15  & 19 & 16 & 12  & 10 & 14 & 12  & 11 \\ \hline
  \end{tabular}
  \caption{The number of iterations of the generalized multigrid (GMG) method, semi-smooth Newton (SSN) method and interior-point (IP) method to reach predefined tolerance.}
  \label{tab:compare_solvers}
\end{table*}

 \section{Conclusion}
In this paper, we introduced an unfitted FE discretization for Signorini's and two-body contact problem.
We utilize the vital vertex algorithm to create a stable Lagrange multiplier space which we used for discretizing the non-penetration condition.
In the numerical section, we evaluated the convergence of discretization error on two different examples of Signorini's problem, which demonstrated optimal convergence properties of the unfitted FE discretization method.
Later, we introduced a generalized multigrid method as an extension of the monotone multigrid method, which can handle linear inequality constraints.
We demonstrated the robustness and the efficiency of the multigrid method for solving Signorini's problem and two-body contact problems.

The generalized multigrid method introduced in this work can be used to solve constraint minimization problems where the number of constraints is significantly smaller than the number of unknowns.
We aim to extend this multigrid method to solve the contact problem with higher-order discretization schemes in fitted and/or unfitted FE framework.
Additionally, the extension of this method for the hyperelastic material models would also be a quite interesting pursuit.

 \section*{Acknowledgments}
The authors would like to thank the Swiss National Science Foundation for their support through the project and the Deutsche Forschungsgemeinschaft (DFG) for their support in the SPP 1962 “ Stress-Based Methods for Variational Inequalities in Solid Mechanics: Finite Element Discretization and Solution by Hierarchical Optimization [186407]”
Additionally, we would also like to gratefully acknowledge the support of Platform for Advanced Scientific Computing (PASC)  through projects FraNetG: Fracture Network Growth and FASTER: Forecasting and Assessing Seismicity and Thermal Evolution in geothermal Reservoirs.

\bibliographystyle{alpha}      \def\url#1{}

\begin{thebibliography}{MHWL12}

\bibitem[AVY19]{akula_mortex_2019}
B.~R. Akula, J.~Vignollet, and V.~A. Yastrebov.
\newblock {MorteX} method for contact along real and embedded surfaces:
  coupling {X}-{FEM} with the {Mortar} method.
\newblock {\em arXiv:1902.04000 [cs, math]}, 2019.

\bibitem[BCH{\etalchar{+}}15]{burman_cutfem:_2015}
E.~Burman, S.~Claus, P.~Hansbo, M.~G. Larson, and A.~Massing.
\newblock {CutFEM}: {Discretizing} geometry and partial differential equations.
\newblock {\em International Journal for Numerical Methods in Engineering},
  104(7):472--501, 2015.

\bibitem[BH91]{barbosa_finite_1991}
H.~J.~C. Barbosa and T.~J.~R. Hughes.
\newblock The finite element method with {Lagrange} multipliers on the
  boundary: circumventing the {Babu{\v s}ka}-{Brezzi} condition.
\newblock {\em Computer Methods in Applied Mechanics and Engineering},
  85(1):109--128, 1991.

\bibitem[BH92]{barbosa_circumventing_1992}
H.~J.~C. Barbosa and T.~J.~R. Hughes.
\newblock Circumventing the {Babu{\v s}ka}-{Brezzi} condition in mixed finite
  element approximations of elliptic variational inequalities.
\newblock {\em Computer Methods in Applied Mechanics and Engineering},
  97(2):193--210, 1992.

\bibitem[BH10]{burman_fictitious_2010}
E.~Burman and P.~Hansbo.
\newblock Fictitious domain finite element methods using cut elements: {I}. {A}
  stabilized {Lagrange} multiplier method.
\newblock {\em Computer Methods in Applied Mechanics and Engineering},
  199(41-44):2680--2686, 2010.

\bibitem[BMW09]{bechet_stable_2009}
{\'E}.~B{\'e}chet, N.~Mo{\"e}s, and B.~Wohlmuth.
\newblock A stable {Lagrange} multiplier space for stiff interface conditions
  within the extended finite element method.
\newblock {\em International Journal for Numerical Methods in Engineering},
  78(8):931--954, 2009.

\bibitem[BPM{\etalchar{+}}03]{belytschko_structured_2003}
T.~Belytschko, C.~Parimi, N.~Mo{\"e}s, N.~Sukumar, and S.~Usui.
\newblock Structured extended finite element methods for solids defined by
  implicit surfaces.
\newblock {\em International Journal for Numerical Methods in Engineering},
  56(4):609--635, 2003.

\bibitem[Bur10]{burman_ghost_2010}
E.~Burman.
\newblock Ghost penalty.
\newblock {\em Comptes Rendus Mathematique}, 348(21-22):1217--1220, 2010.

\bibitem[CH13]{chouly_nitsche-based_2013}
F.~Chouly and P.~Hild.
\newblock A {Nitsche}-{Based} {Method} for {Unilateral} {Contact} {Problems}:
  {Numerical} {Analysis}.
\newblock {\em SIAM Journal on Numerical Analysis}, 51(2):1295--1307, 2013.

\bibitem[CHR15]{chouly_symmetric_2015}
F.~Chouly, P.~Hild, and Y.~Renard.
\newblock Symmetric and non-symmetric variants of {Nitsche}{\textquoteright}s
  method for contact problems in elasticity: theory and numerical experiments.
\newblock {\em Mathematics of Computation}, 84(293):1089--1112, 2015.

\bibitem[CK18]{claus_stable_2018}
S.~Claus and P.~Kerfriden.
\newblock {A} stable and optimally convergent {LaTIn}-{CutFEM} algorithm for
  multiple unilateral contact problems.
\newblock {\em International Journal for Numerical Methods in Engineering},
  113(6):938--966, 2018.

\bibitem[CSS11]{coon_nitsche-eXtended_2011}
E.~T. Coon, B.~E. Shaw, and M.~Spiegelman.
\newblock A {Nitsche}-extended finite element method for earthquake rupture on
  complex fault systems.
\newblock {\em Computer Methods in Applied Mechanics and Engineering},
  200(41):2859--2870, 2011.

\bibitem[DF08]{dolbow_residual-free_2008}
J.~E. Dolbow and L.~P. Franca.
\newblock Residual-free bubbles for embedded {Dirichlet} problems.
\newblock {\em Computer Methods in Applied Mechanics and Engineering},
  197(45-48):3751--3759, 2008.

\bibitem[DK09]{dickopf_efficient_2009}
T.~Dickopf and R.~Krause.
\newblock Efficient simulation of multi-body contact problems on complex
  geometries: {A} flexible decomposition approach using constrained
  minimization.
\newblock {\em International Journal for Numerical Methods in Engineering},
  77(13):1834--1862, 2009.

\bibitem[DK14]{dickopf_evaluating_2014}
T.~Dickopf and R.~Krause.
\newblock Evaluating {Local} {Approximations} of the ${L^2}$-{Orthogonal}
  {Projection} {Between} {Non}-{Nested} {Finite} {Element} {Spaces}.
\newblock {\em Numerical Mathematics: Theory, Methods and Applications},
  7(3):288--316, 2014.

\bibitem[DMB01]{dolbow_extended_2001}
J.~Dolbow, N.~Mo{\"e}s, and T.~Belytschko.
\newblock An extended finite element method for modeling crack growth with
  frictional contact.
\newblock {\em Computer Methods in Applied Mechanics and Engineering},
  190(51):6825--6846, 2001.

\bibitem[FPR16]{fabre_fictitious_2016}
M.~Fabre, J.~Pousin, and Y.~Renard.
\newblock A fictitious domain method for frictionless contact problems in
  elasticity using {Nitsche}{\textquoteright}s method.
\newblock {\em The SMAI journal of computational mathematics}, 2:19--50, 2016.

\bibitem[Glo84]{rolandglowinski1984-02-14}
R.~Glowinski.
\newblock {\em Numerical Methods for Nonlinear Variational Problems}.
\newblock Springer-Verlag, Feb 1984.

\bibitem[Hac86]{wolfganghackbusch1986-01-14}
W.~Hackbusch.
\newblock {\em Multi-Grid Methods and Applications}.
\newblock Springer-Verlag, Jan 1986.

\bibitem[HAD12]{hautefeuille_robust_2012}
M.~Hautefeuille, C.~Annavarapu, and J.~E. Dolbow.
\newblock Robust imposition of {Dirichlet} boundary conditions on embedded
  surfaces.
\newblock {\em International Journal for Numerical Methods in Engineering},
  90(1):40--64, 2012.

\bibitem[HH02]{hansbo_unfitted_2002}
A.~Hansbo and P.~Hansbo.
\newblock An unfitted finite element method, based on
  {Nitsche}{\textquoteright}s method, for elliptic interface problems.
\newblock {\em Computer Methods in Applied Mechanics and Engineering},
  191(47):5537--5552, 2002.

\bibitem[HH04]{hansbo_finite_2004}
A.~Hansbo and P.~Hansbo.
\newblock A finite element method for the simulation of strong and weak
  discontinuities in solid mechanics.
\newblock {\em Computer Methods in Applied Mechanics and Engineering},
  193(33):3523--3540, 2004.

\bibitem[HR09]{haslinger_new_2009}
J.~Haslinger and Y.~Renard.
\newblock A {New} {Fictitious} {Domain} {Approach} {Inspired} by the {Extended}
  {Finite} {Element} {Method}.
\newblock {\em SIAM Journal on Numerical Analysis}, 47, 2009.

\bibitem[JD04]{ji_strategies_2004}
H.~Ji and J.~E. Dolbow.
\newblock On strategies for enforcing interfacial constraints and evaluating
  jump conditions with the extended finite element method.
\newblock {\em International Journal for Numerical Methods in Engineering},
  61(14):2508--2535, 2004.

\bibitem[KDL07]{kim_mortared_2007}
T.~Y. Kim, J.~Dolbow, and T.~Laursen.
\newblock A mortared finite element method for frictional contact on arbitrary
  interfaces.
\newblock {\em Computational Mechanics}, 39(3):223--235, 2007.

\bibitem[KK01]{kornhuber_adaptive_2001}
R.~Kornhuber and R.~Krause.
\newblock Adaptive multigrid methods for {Signorini}'s problem in linear
  elasticity.
\newblock {\em Computing and Visualization in Science}, 4(1):9--20, 2001.

\bibitem[KK19]{kothari_multigrid_2019}
H.~Kothari and R.~Krause.
\newblock A {Multigrid} {Method} for a {Nitsche}-based {Extended} {Finite}
  {Element} {Method}.
\newblock {\em arXiv:1912.00496 [cs, math]}, 2019.

\bibitem[KK21]{kothari_multigrid_2021}
H.~Kothari and R.~Krause.
\newblock Multigrid and saddle-point preconditioners for unfitted finite
  element modelling of inclusions.
\newblock In {\em 14th {WCCM}-{ECCOMAS} {Congress} 2020}, 2021.

\bibitem[KN06]{khoei_contact_2006}
A.~R. Khoei and M.~Nikbakht.
\newblock Contact friction modeling with the extended finite element method
  ({X}-{FEM}).
\newblock {\em Journal of Materials Processing Technology}, 177(1):58--62,
  2006.

\bibitem[Kor94]{kornhuber_monotone_1994}
R.~Kornhuber.
\newblock Monotone multigrid methods for elliptic variational inequalities {I}.
\newblock {\em Numerische Mathematik}, 69(2):167--184, 1994.

\bibitem[Kor96]{kornhuber_monotone_1996}
R.~Kornhuber.
\newblock Monotone multigrid methods for elliptic variational inequalities
  {II}.
\newblock {\em Numerische Mathematik}, 72(4):481--499, 1996.

\bibitem[Kra01]{krause_monotone_2001}
R.~Krause.
\newblock {\em Monotone multigrid methods for {Signorini}'s problem with
  friction}.
\newblock PhD thesis, Freie Universit\"at Berlin, Universit\"at, 2001.

\bibitem[Kra09]{krause_nonsmooth_2009}
R.~Krause.
\newblock A {Nonsmooth} {Multiscale} {Method} for {Solving} {Frictional}
  {Two}-{Body} {Contact} {Problems} in {2D} and {3D} with {Multigrid}
  {Efficiency}.
\newblock {\em SIAM Journal on Scientific Computing}, 31(2):1399--1423, 2009.

\bibitem[Lad12]{ladeveze2012nonlinear}
P.~Ladev{\`e}ze.
\newblock {\em Nonlinear computational structural mechanics: new approaches and
  non-incremental methods of calculation}.
\newblock Springer Science \& Business Media, 2012.

\bibitem[Lau13]{toda.laursen2013-03-13}
T.~A. Laursen.
\newblock {\em Computational Contact and Impact Mechanics: Fundamentals of
  Modeling Interfacial Phenomena in Nonlinear Finite Element Analysis}.
\newblock Springer Science \& Business Media, Mar 2013.

\bibitem[LB08]{liu_contact_2008}
F.~Liu and R.~I. Borja.
\newblock A contact algorithm for frictional crack propagation with the
  extended finite element method.
\newblock {\em International Journal for Numerical Methods in Engineering},
  76(10):1489--1512, 2008.

\bibitem[MBT06]{moes_imposing_2006}
N.~Mo{\"e}s, E.~B{\'e}chet, and M.~Tourbier.
\newblock Imposing {Dirichlet} boundary conditions in the extended finite
  element method.
\newblock {\em International Journal for Numerical Methods in Engineering},
  67(12):1641--1669, 2006.

\bibitem[MDB99]{moes_finite_1999}
N.~Mo{\"e}s, J.~Dolbow, and T.~Belytschko.
\newblock A finite element method for crack growth without remeshing.
\newblock {\em International Journal for Numerical Methods in Engineering},
  46(1):131--150, 1999.

\bibitem[MDH07]{mourad_bubble-stabilized_2007}
H.~M. Mourad, J.~Dolbow, and I.~Harari.
\newblock A bubble-stabilized finite element method for {Dirichlet} constraints
  on embedded interfaces.
\newblock {\em International Journal for Numerical Methods in Engineering},
  69(4):772--793, 2007.

\bibitem[Meh92]{mehrotra_implementation_1992}
S.~Mehrotra.
\newblock On the {Implementation} of a {Primal}-{Dual} {Interior} {Point}
  {Method}.
\newblock {\em SIAM Journal on Optimization}, 2(4):575--601, 1992.

\bibitem[MHWL12]{mueller-hoeppe_crack_2012}
D.~S. Mueller-Hoeppe, P.~Wriggers, and S.~Loehnert.
\newblock Crack face contact for a hexahedral-based {XFEM} formulation.
\newblock {\em Computational Mechanics}, 49(6):725--734, 2012.

\bibitem[NGM{\etalchar{+}}09]{nistor_x-fem_2009}
I.~Nistor, M.~L.~E. Guiton, P.~Massin, N.~Mo{\"e}s, and S.~G{\'e}niaut.
\newblock An {X}-{FEM} approach for large sliding contact along
  discontinuities.
\newblock {\em International Journal for Numerical Methods in Engineering},
  78(12):1407--1435, 2009.

\bibitem[NW00]{jorgenocedal2000-04-27}
J.~Nocedal and S.~Wright.
\newblock {\em Numerical Optimization}.
\newblock Springer, Apr 2000.

\bibitem[Ren13]{renard_generalized_2013}
Y.~Renard.
\newblock Generalized {Newton}{\textquoteright}s methods for the approximation
  and resolution of frictional contact problems in elasticity.
\newblock {\em Computer Methods in Applied Mechanics and Engineering},
  256:38--55, 2013.

\bibitem[SLK20]{sticko_high-order_2020}
S.~Sticko, G.~Ludvigsson, and G.~Kreiss.
\newblock High-order cut finite elements for the elastic wave equation.
\newblock {\em Advances in Computational Mathematics}, 46(3):45, 2020.

\bibitem[SMMB00]{sukumar_extended_2000}
N.~Sukumar, N.~Mo{\"e}s, B.~Moran, and T.~Belytschko.
\newblock Extended finite element method for three-dimensional crack modelling.
\newblock {\em International Journal for Numerical Methods in Engineering},
  48(11):1549--1570, 2000.

\bibitem[Ulb11]{michaelulbrich2011-07-13}
M.~Ulbrich.
\newblock {\em Semismooth Newton Methods for Variational Inequalities and
  Constrained Optimization Problems in Function Spaces (MPS-SIAM Series on
  Optimization)}.
\newblock Society for Industrial \& Applied Mathematics, Jul 2011.

\bibitem[WK03]{wohlmuth_monotone_2003}
B.~I. Wohlmuth and R.~H. Krause.
\newblock Monotone {Multigrid} {Methods} on {Nonmatching} {Grids} for
  {Nonlinear} {Multibody} {Contact} {Problems}.
\newblock {\em SIAM Journal on Scientific Computing}, 25(1):324--347, 2003.

\bibitem[Woh11]{wohlmuth_variationally_2011}
B.~Wohlmuth.
\newblock Variationally consistent discretization schemes and numerical
  algorithms for contact problems.
\newblock {\em Acta Numerica}, 20:569--734, 2011.

\bibitem[Wri06]{peterwriggers2006-06-28}
P.~Wriggers.
\newblock {\em Computational Contact Mechanics}.
\newblock Springer, Jun 2006.

\bibitem[ZKN{\etalchar{+}}16]{utopiagit}
P.~Zulian, A.~Kopani{\v c}{\'a}kov{\'a}, M.~C.~G. Nestola, A.~Fink, N.~Fadel,
  A.~Rigazzi, V.~Magri, T.~Schneider, E.~Botter, J.~Mankau, and R.~Krause.
\newblock {U}topia: {A} {C}++ embedded domain specific language for scientific
  computing. {G}it repository.
\newblock https://bitbucket.org/zulianp/utopia, 2016.

\end{thebibliography}

\newcommand{\etalchar}[1]{$^{#1}$}

\end{document}